\documentclass[pdflatex,sn-mathphys-num]{sn-jnl}% Math and Physical Sciences Numbered Reference Style
\usepackage{graphicx}%
\usepackage{multirow}%
\usepackage{amsmath,amssymb,amsfonts}%
\usepackage{amsthm}%
\usepackage{mathrsfs}%
\usepackage[title]{appendix}%
\usepackage{xcolor}%
\usepackage{textcomp}%
\usepackage{manyfoot}%
\usepackage{booktabs}%
\usepackage{algorithm}%
\usepackage{algorithmicx}%
\usepackage{algpseudocode}%
\usepackage{listings}%
\usepackage{bm}
\usepackage{mathrsfs}

\theoremstyle{thmstyleone}%
\newtheorem{theorem}{Theorem}%  meant for continuous numbers
\newtheorem{proposition}[theorem]{Proposition}% 
\newtheorem{lemma}[theorem]{Lemma}
\newtheorem{corollary}[theorem]{Corollary}

\theoremstyle{thmstyletwo}%
\newtheorem{remark}{Remark}%

\theoremstyle{thmstylethree}%
\newtheorem{definition}{Definition}%
\newtheorem{assumption}{Assumption}

\newcommand{\cB}{\mathcal{B}}

\newcommand{\cF}{\mathcal{F}}

\newcommand{\cH}{\mathcal{H}}
\newcommand{\cI}{\mathcal{I}}

\newcommand{\cL}{\mathcal{L}}

\newcommand{\cP}{\mathcal{P}}

\newcommand{\cR}{\mathcal{R}}
\newcommand{\cS}{\mathcal{S}}
\newcommand{\cT}{\mathcal{T}}

\newcommand{\cX}{\mathcal{X}}
\newcommand{\cY}{\mathcal{Y}}

\newcommand{\sigmam}{\Sigma_m^\phi}
\DeclareMathOperator*{\argmin}{argmin}

\newcommand{\NN}{\mathbb{N}}

\newcommand{\RR}{\mathbb{R}}

\newcommand{\dd}{\mathop{}\!\mathrm{d}}
\newcommand{\tr}{\mathrm{Tr}}

\begin{document}

\title[Convergence of the RFM]{A Convergence Analysis of the Random Feature Method for Linear PDEs}


\author[1,2,3]{\fnm{Jingrun} \sur{Chen}}\email{jingrunchen@ustc.edu.cn}

\author*[2,3]{\fnm{Yixin} \sur{Luo}}\email{seeing@ustc.edu.cn}

\author[1]{\fnm{Xiaojun} \sur{Peng}}\email{pengxiaojun@mail.ustc.edu.cn}

%\affil*[1]{\orgdiv{Department}, \orgname{Organization}, \orgaddress{\street{Street}, \city{City}, \postcode{100190}, \state{State}, \country{Country}}}
\affil[1]{\orgdiv{School of Mathematical Sciences}, \orgname{University of Science and Technology of China}, \orgaddress{\street{96 Jinzhai Road}, \city{Hefei}, \postcode{230026}, \state{Anhui}, \country{China}}}

\affil[2]{\orgdiv{Suzhou Institute of Advanced
Research}, \orgname{University of Science and Technology of China}, \orgaddress{\street{166 Renai Road}, \city{Suzhou}, \postcode{215123}, \state{Jiangsu}, \country{China}}}

\affil[3]{\orgdiv{School of Artificial Intelligence and Data Science}, \orgname{University of Science and Technology of China}, \orgaddress{\street{166 Renai Road}, \city{Suzhou}, \postcode{215123}, \state{Jiangsu}, \country{China}}}

%\affil[3]{\orgdiv{Department}, \orgname{Organization}, \orgaddress{\street{Street}, \city{City}, \postcode{610101}, \state{State}, \country{Country}}}


\abstract{The Random Feature Method (RFM) distinguishes itself among machine learning approaches to partial differential equations (PDEs) through its superior accuracy and computational efficiency. Nevertheless, the theoretical underpinnings of RFM remain incompletely elucidated. This paper provides a comprehensive analysis of RFM within the framework of kernel methods, focusing on two principal aspects: (i) approximation by finitely many random features in
Sobolev norms and the PDE residual seminorm, and (ii) estimation over the induced reproducing kernel Hilbert
space (RKHS) from finitely many noiseless collocation
observations. We relate both problems to kernel ridge
regression (KRR) and establish upper and lower bounds
in terms of the eigenvalues of the associated integral
operators under the respective assumptions. Combining the approximation and estimation lower bounds,
we establish a lower bound under simultaneous restrictions
on the numbers of random features and collocation observations. We further quantify how activation regularity affects
these eigenvalues, obtaining algebraic upper bounds
for finitely smooth activations and stretched-exponential
upper bounds for analytic activations. }




\maketitle

\section{Introduction}
Partial differential equations (PDEs) constitute a cornerstone in the mathematical modeling of physical and engineering systems. In recent years, the advent of deep learning techniques has facilitated the emergence of novel methodologies for solving PDEs \cite{han2017deep,han2018solving,yu2018deep,sirignano2018dgm,raissi2019physics,zang2020weak}. 
The notable efficacy of these methods stems from their mesh-free nature and their capacity to alleviate the curse of dimensionality. However, their application is constrained by two primary bottlenecks: limited accuracy and high computational expense. These limitations fundamentally originate from the nonconvex optimization landscape encountered during neural network training. 
Several approaches address this difficulty by fixing
the inner-layer parameters \cite{dong2021local, chen2022bridging,shang2023randomized}. For linear PDEs, optimizing the output-layer coefficients
then reduces to a linear least-squares problem, thereby avoiding the complexities associated with nonconvex optimization. A prominent instance is the random feature method
(RFM)~\cite{chen2022bridging}, which uses a random feature
model~\cite{rahimi2007random} to represent the approximate
solution. This model is a two-layer neural network whose
hidden-layer parameters are randomly sampled and then
held fixed. Specifically, the corresponding function
class is
\[
\sigmam
:=
\left\{
u_{\bm\alpha}:
u_{\bm\alpha}(\bm x)
=
\frac{1}{\sqrt{m}}
\sum_{i=1}^{m}
\alpha_i\phi(\bm x,\bm\theta_i),
\quad
\bm\alpha\in\mathbb R^m
\right\},
\]
where $\bm\theta_1,\ldots,\bm\theta_m$ are the sampled
feature parameters and $\bm\alpha$ contains the trainable
output coefficients.
RFM has exhibited exceptional performance in diverse numerical experiments spanning structural analysis, interface problems, and time-dependent problems \cite{chen2022bridging,luo2023random,chi2024random,sun2025}. Notably, in many cases, the achieved errors decay at an exponential or near-exponential rate as the number of random features or collocation points increases. Despite these empirical successes, a comprehensive theoretical explanation for such convergence rates remains incomplete. Closing this theoretical gap is essential both to substantiate the empirical efficacy of RFM and to extend its applicability to a wider range of challenging problems.

Throughout this paper, we consider linear boundary value problems of
the form
\[
\begin{cases}
\cL u^* = f^* & \text{in } \Omega,\\
\cB u^* = g^* & \text{on } \partial\Omega,
\end{cases}
\]
where $\cL$ is a linear differential operator, $\cB$ is a linear
boundary operator, and $u^*$ is the target solution. Our analysis of RFM distinguishes two sources of error. The
approximation error is the distance from $u^*$ to the random feature model. The estimation error arises from recovering $u^*$
from finitely many noiseless observations of $f^*$ and
$g^*$. Since random features were introduced as a computational
surrogate for kernel methods~\cite{rahimi2007random},
we analyze both problems through the induced kernels
and their associated integral operators.
For finite-feature approximation, we
construct kernels on parameter space from the Sobolev inner products of
the features and from the weighted inner products of their images
under the differential and boundary operators. For each error metric, our duality result reformulates the
regularized finite-feature approximation problem as noiseless KRR
on the parameter space, where the sampled parameters
serve as observation sites for the transformed target. Through this correspondence, the noiseless recovery lower bound and
the KRR upper bound in
\cite{long2025optimalratessaturationnoiseless}
yield the lower and upper bounds for the worst-case approximation
error in each metric. For finite collocation estimation,
we construct a spatial kernel for which the
rescaled interior and boundary observations are point values of
the transformed target. The collocation lower bound then follows
from the noiseless recovery bound in
\cite[Theorem~4.1]{long2025optimalratessaturationnoiseless}. For the collocation upper bound, we adapt the concentration
argument in the proof of the upper bound in \cite{long2025optimalratessaturationnoiseless} to account for the separate
sampling of interior and boundary points. We control each
empirical operator under its own sampling distribution
and combine the resulting interior and boundary operator
deviation bounds to estimate the population residual
of the regularized estimator over the full RKHS. Although finite-feature approximation admits a direct analysis,
we use duality to study it within the same noiseless KRR framework
as recovery from collocation observations. Our approach based on duality yields the same best approximation
error as the direct approach
(see Remark~\ref{rem:direct-dual-approximation}).
Both the approximation and collocation bounds depend on the
eigenvalues of the corresponding integral operators. To make the bounds explicit, we relate the decay of the associated
eigenvalues to the regularity of the activation function. We use
Volkov's singular-value estimates~\cite{volkov2024eigenvalue}
together with our ambient regularity lemmas to derive upper bounds
for the eigenvalues governing Sobolev approximation and PDE
collocation.

Our main results are summarized below, with precise
statements and assumptions given in Section~\ref{sec:main}.
\begin{itemize}

\item
\textbf{Approximation in Sobolev norms and the residual seminorm.}
We establish upper and lower bounds for the worst-case
approximation error over the unit ball of $\mathcal H_k$.
Under the respective assumptions of
Theorems~\ref{upper_appro} and~\ref{lower_appro}
for Sobolev approximation and
Theorems~\ref{upper_appro_residual}
and~\ref{lower_appro_residual}
for residual approximation, the bounds take the form
\[
\begin{aligned}
\tilde\mu_{m+1}^{(r)}
&\leq
\sup_{\|u^*\|_{\mathcal H_k}\leq1}
\inf_{u\in\sigmam}
\|u-u^*\|_{H^r(\Omega)}^2
\lesssim
\tilde\mu_{l_m}^{(r)},\\
\tilde\mu_{m+1}^{(\cR_\eta)}
&\leq
\sup_{\|u^*\|_{\mathcal H_k}\leq1}
\inf_{u\in\sigmam}
\|u-u^*\|_{\cR_\eta}^2
\lesssim
\tilde\mu_{l_m}^{(\cR_\eta)}.
\end{aligned}
\]
For each error measure, the upper bound holds with
probability at least $1-\delta$ under its corresponding
feature sampling distribution, whereas the lower bound
holds for every fixed admissible choice of feature
parameters.
Here, $\tilde\mu_j^{(r)}$ and
$\tilde\mu_j^{(\cR_\eta)}$ denote the $j$-th largest
positive eigenvalues of the integral operators defined
in \eqref{int_oper_appro} and
\eqref{int_oper_appro_residual}, respectively.
The index $l_m$ is specified in the respective
upper-bound theorem with \(
l_m\asymp\frac{m}{\log(m/\delta)}.
\)

\item
\textbf{Estimation from finitely many collocation observations.}
We analyze regularized recovery over the full RKHS
using the estimator
\[
\hat u_{n,\eta,\lambda}
:=
\argmin_{u\in\mathcal H_k}
\left\{
\hat{\cR}_{n,\eta}(u)
+\lambda\|u\|_{\mathcal H_k}^2
\right\},
\qquad \lambda>0,
\]
where $\hat{\cR}_{n,\eta}$ is the weighted empirical
PDE residual defined in \eqref{discre_loss}, constructed
from $n$ noiseless interior and boundary observations.
Under the sampling and regularization conditions of
Theorem~\ref{thm: Estimation_error upper},
the estimator $\hat u_{n,\eta,\lambda}$ satisfies,
with probability at least $1-\delta$,
\[
\sup_{\|u^*\|_{\mathcal H_k}\leq1}
\|\hat u_{n,\eta,\lambda}-u^*\|_{\cR_\eta}^2
\lesssim
\mu_{l_n}^{\mathrm{PDE}},
\qquad
l_n\asymp\frac{n}{\log(n/\delta)}.
\]
Here
$\mu_j^{\mathrm{PDE}}$ denotes the $j$-th largest positive
eigenvalue of the integral operator defined in
\eqref{pde_int_operator}.
The index $l_n$ is specified in the theorem.
For every fixed admissible collection of collocation
observations, we also establish
\[
\inf_{G\in\mathcal M_n^{\mathrm{PDE}}}
\sup_{\|u^*\|_{\mathcal H_k}\leq1}
\|G(S_n(u^*))-u^*\|_{\cR_\eta}^2
\gtrsim
\mu_{n+1}^{\mathrm{PDE}}.
\]
Here, $S_n(u^*)$ denotes the observations, and
$\mathcal M_n^{\mathrm{PDE}}$ is the class of reconstruction
procedures based only on these observations.

\item
\textbf{A lower bound under joint feature and observation restrictions.}
We consider reconstruction from $n$ collocation observations
with outputs restricted to a random feature model spanned
by $m$ features.
Let $E_{m,n}$ denote the smallest worst-case squared
residual error over the unit ball of $\mathcal H_k$
among all such reconstruction procedures.
For fixed positive proportions of interior and boundary
observations, we establish
\[
E_{m,n}
\gtrsim
\max\left\{
\mu_{m+1}^{\mathrm{PDE}},
\mu_{n+1}^{\mathrm{PDE}}
\right\}
\]
for every fixed admissible choice of feature parameters
and collocation points, with an implicit constant
independent of $m$ and $n$.
This bound identifies the separate limitations imposed
by the feature space and the available observations
and applies, in particular, to the RFM estimator.

\item
\textbf{Dependence of eigenvalue decay on activation regularity.}
We derive eigenvalue upper bounds for features of the form
$\phi(\bm x,\bm\theta)=\sigma(\bm k^\top\bm x+b)$.
For the PDE eigenvalues, we specialize to second-order
uniformly elliptic equations with Dirichlet boundary
conditions.
Under the finite-regularity assumptions of
Theorem~\ref{activation_fun_rate}, with
$p\geq\max\{2,r\}$, we obtain
\[
\tilde\mu_j^{(r)}
=
o\!\left(j^{-1-\frac{2(p-r)}{d}}\right),
\qquad
\mu_j^{\mathrm{PDE}}
=
o\!\left(j^{-1-\frac{2(p-2)}{d}}\right).
\]
Under the analytic assumptions of the same theorem,
there exists $c>0$ such that
\[
\tilde\mu_j^{(r)}
\lesssim e^{-c j^{1/d}},
\qquad
\mu_j^{\mathrm{PDE}}
\lesssim e^{-c j^{1/d}}.
\]
\end{itemize}

\subsection{Related Work}
This section reviews related work on random-feature
approximation, spectral analysis of kernel methods,
and kernel collocation for PDEs.

\textbf{Random-feature approximation.}
Early universal approximation results cover networks
with randomly generated hidden nodes~\cite{universal2006}.
For uniformly bounded features, a fixed target admits
a high-probability $L^2$ approximation bound of order
$O(m^{-1/2})$ if its integral representation has a
coefficient function whose ratio to the sampling density
is bounded~\cite{rahimi2008uniform}.
For $\mathrm{ReLU}^k$ features with well-distributed
fixed parameters on $\mathbb S^d$, the $H^r$ approximation
error is of order $O(m^{-1/2-(2(k-r)+1)/(2d)})$ for
$0\leq r\leq k$ and targets in
$H^{(d+2k+1)/2}(\Omega)$ on a bounded Lipschitz
domain $\Omega$~\cite[Theorem~2.2]{liu2025integral}.
For parameters sampled independently and uniformly on
$\mathbb S^d$, the corresponding high-probability $H^r$
bound incurs a logarithmic loss in the number of
features~\cite[Theorem~6.2]{liu2025integral}.
\cite{fabiani2025random} studies best $L^p$ approximation
on an interval by random projection networks with fixed hidden
parameters and smooth nonpolynomial activations. \cite{ming2025spectral} analyze RFM for one-dimensional
second-order elliptic problems and derive error estimates that
depend on the regularity of the solution. In multiple dimensions,
they establish high-probability approximation bounds for trigonometric features~\cite{ming2026spectral}. Our analysis starts from a general feature map and its induced
RKHS. Under the stated spectral and sampling assumptions, we
derive upper and lower bounds for the worst-case approximation
error over the RKHS unit ball, both in Sobolev norms and in the
PDE residual seminorm. For each error metric, the bounds are
expressed in terms of the eigenvalues of the corresponding
operator on the parameter space. This formulation applies without
prescribing a particular activation function. We then use the regularity of the activation function to obtain
explicit upper bounds on the eigenvalues.

\textbf{Spectral analysis of kernel methods.}
For random-feature approximation and kernel quadrature,
the eigenvalues of the associated integral operator yield
upper and lower bounds on the number of samples needed
to attain a prescribed accuracy. The upper bounds can
be attained by sampling from a suitable nonuniform
distribution~\cite{bach2017equivalence}.
A duality framework connects random-feature approximation
and estimation~\cite{chen2025duality}.
In our setting, the sampled feature parameters serve as
observation points in a noiseless KRR problem on parameter space.
This formulation allows us to bound the finite-feature
approximation error using KRR estimates that depend on
eigenvalue decay and target
regularity~\cite{long2025optimalratessaturationnoiseless}.
Sobolev and analytic regularity of an integral kernel
also yield decay estimates for the singular values of
the corresponding operator~\cite{volkov2024eigenvalue}.
We draw on these results to analyze approximation in
Sobolev norms and the PDE residual seminorm, recovery
from interior and boundary observations, and the dependence
of eigenvalue bounds on activation regularity.

\textbf{Kernel collocation for PDEs.}
Radial basis function collocation is a classical kernel
method for solving PDEs~\cite{kansa1,kansa2}.
Native-space interpolation theory relates approximation
errors to fill distance and kernel smoothness~\cite{rbf2}.
For second-order linear elliptic equations with Dirichlet
boundary conditions, constrained and weighted least-squares
kernel collocation methods admit Sobolev error bounds
under suitable smoothness and point-density
conditions~\cite{cheung2018h}.
For nonlinear and parametric PDEs, error analyses of
kernel and Gaussian process methods use RKHS embeddings,
PDE stability, and Sobolev interpolation
estimates~\cite{batlle2025error}.
Our analysis complements these results by deriving upper and lower bounds in terms of
operator eigenvalues for two distinct problems:
approximation by $m$ features and recovery over the
full RKHS from $n$ collocation observations.
When reconstruction is subject to both restrictions,
the worst-case squared residual error is at least the
larger of the two lower bounds.
\subsection{Notations}
We use $\NN^+$ to denote the set of positive integers and set $\NN=\{0\}\cup \NN^+$. 
We shall use $\|\cdot\|$ to denote an operator norm without specifying the input and output spaces for brevity. For any Banach spaces $\cX$ and $\cY$, let $\mathcal{B}(\cX
)$ denote the bounded linear operators on $\cX$ and $\cF(\cX;\cY)$ denote the functions $\{f:\cX \to \cY\}$. Let  $\cP(\cX)$ denote the set of probability measures over $\cX$. For any $\rho \in \cP(\cX)$, we denote by $L^p(\cX, \rho)$ the classical $L^p$ space, and often write $L^p(\rho)$ when the underlying space $\cX$ is clear from context. For any $f, g \in L^2(\rho)$, we write $\langle f, g \rangle_\rho = \int_{\cX} f(\bm{x})g(\bm{x}) \dd \rho(\bm{x})$ for the inner product, and use $\|\cdot\|_\rho$ to denote the associated norm. For a bounded domain $\Omega\subset\mathbb R^d$ and a
nonnegative integer $r$, let $H^r(\Omega)=W^{r,2}(\Omega)$
denote the Sobolev space with respect to Lebesgue measure,
equipped with the inner product
\(
\langle f,g\rangle_{H^r(\Omega)}
:=
\sum_{|\bm\beta|\leq r}
\int_\Omega
\mathrm D^{\bm\beta}f(\bm x)\,
\mathrm D^{\bm\beta}g(\bm x)\,\dd\bm x,
\)
where $\bm\beta$ is a multi-index and
$\mathrm D^{\bm\beta}$ denotes the corresponding weak
derivative. We denote the induced norm by
$\|\cdot\|_{H^r(\Omega)}$.
In particular, $H^0(\Omega)=L^2(\Omega)$. 
%Let $\dd \rho_B$ denote a prescribed probability measure on $\partial \Omega$. For any $f,g\in L^2(\rho_B)$, let $\langle f,g\rangle_{\rho_B}=\int_{\partial \Omega} f(\bm{x})g(\bm{x})\dd \rho_B(\bm{x})$.
 % Given two Banach spaces $\cA$ and $\cB$  with $\cA\subset \cB$, we denote by $\|\cA\hookrightarrow \cB\|:=\sup_{\|x\|_{\cA}=1}\|x\|_{\cB}$ the norm of the inclusion operator.
 Let $\cH_1, \cH_2$ be two Hilbert spaces, and $v\in \cH_1, u\in \cH_2$. The operator $u\otimes v: \cH_1\mapsto \cH_2$ is defined as $(u\otimes v) x = \langle v,x\rangle_{\cH_1} u$ for $x\in \cH_1$. We use $\cI$ to denote a general identity operator without specifying the domain for brevity and let $\bm{I}_d$ be the $d\times d$ identity matrix. 
 % For a bounded linear operator $A$  over a Hilbert space and $r\in\RR$, we use $(A+\lambda)^{r}$ as a shorthand for $(A+\lambda \cI)^{r}$. Let $\SS^{d-1}=\{x\in\RR^d: \|x\|_2=1\}$. 
We use big-O notations like $O(\cdot)$ and $o(\cdot)$ to hide constants and  $\tilde{O}(\cdot)$ to further hide logarithmic factors.
We also adopt the notation $a \lesssim b$ to indicate $a  =  O(b)$ and write $a\asymp b$ when both $a \lesssim b$ and $b \lesssim a$ hold.

The paper is organized as follows. Section~\ref{sec:pre} provides the necessary preliminaries for our analysis. Section~\ref{sec:problem} introduces the problem setting. Our main results are presented in Section~\ref{sec:main}, with the corresponding analysis provided in Section~\ref{Approximation and Estimation Analysis}. Conclusions follow in Section~\ref{sec:conclusions}. Auxiliary results for our proofs are provided in Appendix~\ref{math_preliminaries} and the proofs of activation-dependent eigenvalue bounds are deferred to Appendix~\ref{proof_of_thm}.
\section{Preliminaries}
\label{sec:pre}

In this section, we recall some basic definitions and results concerning reproducing kernel Hilbert spaces, kernel ridge regression, and degrees of freedom, which will be used in the subsequent sections.

\subsection{Kernel and RKHS}  

Let $\cX$ be a measurable input space, and $\rho\in \cP(\cX)$ denote the reference probability measure. 
A function $k:\cX\times\cX \mapsto \RR$ is said to be  a (positive semidefinite) kernel if 1) it is symmetric, i.e. $k(x,x') = k(x',x), \forall x,x' \in \cX$; and 2)  $\forall n \ge 1$, $x_1,\dots,x_n \in \cX$, the kernel matrix $(k(x_i,x_j))_{i,j}\in\RR^{n\times n}$ is positive semidefinite. Given a kernel $k$, there exists a unique Hilbert space $\cH$ such that $k(x,\cdot)\in \cH, \forall x\in \cX$; and $\langle f,k(x,\cdot)\rangle_\cH = f(x), \forall f\in\cH,x\in\cX$.
The second property is known as the reproducing property, and thus $\cH$ is often referred to as a reproducing kernel Hilbert space (RKHS)~\cite{aronszajn1950theory}. The RKHS can also be characterized using the associated integral operator $\cT: L^2(\cX,\rho)\mapsto L^2(\cX,\rho)$, given by
\begin{equation}\label{eqn: kernel integral}
    \cT f = \int_\cX k(\cdot,x') f(x') \dd \rho(x').
\end{equation}

Throughout this paper, we assume 
$\int k(x,x)\dd\rho(x)<\infty$, which guarantees that $\cT$ is a self-adjoint, positive, and compact operator; see \cite[Theorem 4.27]{steinwart2008support}. By the spectral theory of compact operators, we have the following  decomposition:
\begin{equation}\label{eqn: mercer's decomposition}
\cT = \sum_{j=1}^\infty \mu_j e_j\otimes e_j,
\end{equation}
 where $\{\mu_j\}_{j=1}^\infty$ are the eigenvalues in a decreasing order and $\{e_j\}_{j=1}^\infty$ are the corresponding orthonormal eigenfunctions. 
 Without loss of generality, we assume that $\{e_j\}_{j=1}^\infty$ forms a complete basis of $L^2(\cX,\rho)$. Otherwise we can always focus on the subspace spanned by the eigenfunctions associated with nonzero eigenvalues. We impose the following assumption on the
eigenvalue sequences considered in this paper.
\begin{assumption}
\label{assump:eigenvalue-decay}
Each eigenvalue sequence $\{\mu_j\}_{j\geq1}$ considered
in this paper is arranged in nonincreasing order and
satisfies one of the following conditions:
\begin{enumerate}
    \item Algebraic decay:
    \[
        \mu_j\asymp j^{-t}
        \qquad\text{for some }t>1.
    \]
    \item Stretched-exponential decay:
    \[
        \mu_j\asymp\tau^{-j^{1/d}}
        \qquad\text{for some }\tau>1
        \text{ and }d\in\mathbb N_{+},
    \]
\end{enumerate}
where $d$ is the spatial dimension, which will be specified in each setting.
\end{assumption}
We focus on these two decay regimes, for which
Lemma~\ref{lemma: N(gamma)} provides explicit estimates
of the degree of freedom.
Other eigenvalue decay regimes may be analyzed using
similar arguments.

We further extend $\cH$ to the following interpolation space with arbitrary $s \in \RR$. 
\begin{definition}[Interpolation space]\label{def: Hks-space}
For any $s\in \RR$ and $f\in L^2(\cX,\rho)$, we define
$
\|f\|_{\cH^s}=\|\cT^{-\frac{s}{2}}f\|_\rho$ and the interpolation space $\cH^{s} = \{ f: \| f \|_{\cH^{s}} < \infty \}.
$  
\end{definition}
By the spectral decomposition~\eqref{eqn: mercer's decomposition},  $\cH^s$ can  be  written in the following form
\begin{align}\label{eqn: Hks-space2}
\cH^s=\left\{\sum_{j=1}^\infty a_j\mu_j^{s/2}e_j: \sum_{j=1}^\infty a_j^2<\infty\right\},\quad \left\|\sum_{j=1}^\infty a_j\mu_j^{s/2}e_j\right\|_{\cH^s}^2=\sum_{j=1}^\infty a_j^2.
\end{align}
It is obvious that $\cH^1=\cH$, $\cH^0=L^2(\cX,\rho)$, and $\cH^s$  provides a natural interpolation between them, and beyond (i.e. $s>1$ or $s<0$). For more discussion, we refer to \cite[Section 2]{long2025optimalratessaturationnoiseless}.

\subsection{Weighted KRR}
Suppose that \(x_1,\ldots,x_n\) are i.i.d. samples drawn from a sampling distribution \(\rho^{\prime}\), possibly different from the reference distribution \(\rho\). We impose the following assumption on every pair of reference and sampling measures used
throughout the paper.
\begin{assumption}
\label{sample distribution}
 Let $\nu$ and $\nu^{\prime}$ be probability measures on the same measurable
space. The sampling measure $\nu^{\prime}$ is absolutely continuous with
respect to the reference measure $\nu$. Moreover, its
Radon--Nikodym derivative
\(
    r=\frac{\dd\nu^{\prime}}{\dd\nu}
\)
satisfies $r>0$ $\nu$-almost everywhere.   
\end{assumption}
We apply Assumption~\ref{sample distribution} to
$(\rho,\rho')$ on $\mathcal X$, with $q=\frac{\dd \rho^{\prime}}{\dd \rho}$.
For a target function $f^*\in \cH$, we consider the following KRR:
\begin{equation*}\label{eqn: krr-def}
    \hat{f}_{n,\lambda} =\argmin_{f \in \cH} \left[\frac{1}{n}\sum_{i=1}^n\frac{1}{q(x_i)} \left(f(x_i) - f^*(x_i)\right)^2 + \lambda \|f\|_{\cH}^2\right].
\end{equation*}
The weighted empirical risk is an unbiased estimate of the true risk  and standard KRR is recovered by setting $q(x)\equiv 1$. By the representer theorem \citep{scholkopf2001generalized}, $\hat{f}_{n,\lambda}$ must lie in $\mathrm{span}\{k(x_1,\cdot),\dots,k(x_n,\cdot)\}$.
% Let $\hat{\bm{k}}_n : \cX \to \mathbb{R}^n$ with $(\hat{\bm{k}}_n (\cdot))_i=\frac{1}{\sqrt{nq(x_i)}}k(x_i,\cdot)$, $\hat{\bm{y}}\in \mathbb{R}^n$ with $(\hat{\bm{y}})_i=\frac{f^*(x_i)}{\sqrt{nq(x_i)}}$, and $\bm{K}_n = \left(\frac{1}{n\sqrt{q(x_i)q(x_j)}}k(x_i,x_j)\right)_{1\leq i,j\leq n}\in \mathbb{R}^{n\times n}$. Then, the KRR solution can be expressed as: $\hat{f}_{\lambda}=\hat{\bm{k}}_n(\cdot)^{\top}\hat{\bm{b}}_{\lambda}$ with
% \begin{equation}
% \label{para_krr_sol}
%     \hat{\bm{b}}_{\lambda} = (\bm{K}_n+\lambda \bm{I}_n)^{-1} \hat{\bm{y}}.
% \end{equation}
Define
\[
\hat{\bm{k}}_n(x)
=
\left(
\frac{k(x_i,x)}{\sqrt{nq(x_i)}}
\right)_{i=1}^n,
\quad
\hat{\bm{y}}
=
\left(
\frac{f^*(x_i)}{\sqrt{nq(x_i)}}
\right)_{i=1}^n,
\quad
\bm{K}_n
=
\left(
\frac{k(x_i,x_j)}
{n\sqrt{q(x_i)q(x_j)}}
\right)_{i,j=1}^n.
\]
Then the KRR solution can be expressed as: 
\(
    \hat{f}_{n, \lambda} = \hat{\bm{k}}_{n} (\cdot)^{\top} (\bm{K}_{n} + \lambda \bm{I}_{n})^{-1} \hat{\bm{y}}. 
\)
\subsection{Degrees of Freedom}
\label{subsec:degrees-of-freedom}
To study the convergence of learning within an RKHS, the commonly used notion of complexity is the degrees of freedom (DoFs) (see, e.g., \citep{zhang2005learning,caponnetto2007optimal,blanchard2018optimal,lin2020optimal}).
To facilitate our analysis, we specialize
\cite[Definition~3.1]{long2025optimalratessaturationnoiseless}
to the case where the DoF exponent $\gamma=1$.

\begin{definition}
For $\lambda>0$ and a sampling density $q$ with respect to
$\rho$, satisfying $q>0$ $\rho$-almost everywhere, define
\begin{equation*}
N(\lambda)
:=
\operatorname{tr}\!\left[(\cT+\lambda I)^{-1}\cT\right],
\qquad
F_q(\lambda)
:=
\operatorname*{ess\,sup}_{x\in\cX}
\frac{\|(\cT+\lambda I)^{-1/2}k(x,\cdot)\|_{\cH}^2}{q(x)},
\end{equation*}
where the essential supremum is taken with respect to $\rho$.
\end{definition}
By the spectral decomposition
\eqref{eqn: mercer's decomposition},
\begin{equation}
\label{eq: dof}
N(\lambda)
=
\sum_{j=1}^{\infty}\frac{\mu_j}{\mu_j+\lambda},
\qquad
F_q(\lambda)
=
\operatorname*{ess\,sup}_{x\in\cX}
\frac1{q(x)}
\sum_{j=1}^{\infty}
\frac{\mu_j}{\mu_j+\lambda}|e_j(x)|^2.
\end{equation}
We immediately have
\[
F_q(\lambda)
\geq
\int_{\cX}
\frac{1}{q(x)}
\sum_{j=1}^{\infty}
\frac{\mu_j}{\mu_j+\lambda}|e_j(x)|^2
\,\dd\rho^{\prime}(x)
=
\sum_{j=1}^{\infty}
\frac{\mu_j}{\mu_j+\lambda}
=
N(\lambda).
\] 
One convenient choice for controlling $F_q(\lambda)$
in terms of $N(\lambda)$ is the equal-weight mixture
of the reference density $1$ and the leverage-score
density defined in
\cite[Lemma~3.5]{long2025optimalratessaturationnoiseless}
with DoF exponent $\gamma=1$:
\begin{equation}
\label{def: sample density}
q^*(k,\rho,\lambda;x)
:=
\frac12
+
\frac{1}{2N(\lambda)}
\sum_{j=1}^{\infty}
\frac{\mu_j}{\mu_j+\lambda}|e_j(x)|^2.
\end{equation}
Then $q^*(k,\rho,\lambda;\cdot)$ is a probability density
with respect to $\rho$ and satisfies
$q^*(k,\rho,\lambda;x)\geq1/2$ for $\rho$-almost every $x$. Moreover,
\begin{equation}
\label{eq: control F}
F_{q^*(k,\rho,\lambda;\cdot)}(\lambda)
=
2N(\lambda)
\operatorname*{ess\,sup}_{x\in\cX}
\left(
1-\frac{1}{2q^*(k,\rho,\lambda;x)}
\right)
\leq 2N(\lambda).
\end{equation}
For further discussion of the DoFs, see
\cite[Section~3]{long2025optimalratessaturationnoiseless}.

\section{Problem Setting}
\label{sec:problem}
\subsection{PDE Problem}
Let \(\Omega\subset\mathbb{R}^d\) be a bounded domain with Lipschitz boundary. We consider the following boundary value problem:
\begin{equation} 
\label{problem}
    \left\{\begin{aligned}
\mathcal{L}u(\bm{x})=f^*(\bm{x}), & \quad \bm{x} \in \Omega, \\
\mathcal{B}u(\bm{x})=g^*(\bm{x}), & \quad \bm{x} \in \partial \Omega,
\end{aligned}\right.
\end{equation}
where $\cL$ and $\cB$ are both linear differential operators and $f^*$ and $g^*$ are the prescribed interior and
boundary data, respectively.

Let $\rho_I\in\mathcal{P}(\Omega)$ and
$\rho_B\in\mathcal{P}(\partial\Omega)$ be reference probability measures on the interior and the boundary, respectively. Assume that $\dd\rho_I=w_I\,\dd\bm x$ and
$\dd\rho_B=w_B\,\dd S$, where $\dd\bm x$ denotes
Lebesgue measure on $\Omega$, $\dd S$ denotes surface
measure on $\partial\Omega$, and both densities are
strictly positive almost everywhere with respect
to their corresponding measures. For a prescribed penalty parameter $\eta>0$, define the population residual functional by
\begin{equation}
\label{cont_loss}
\cR_\eta(u) = \|\mathcal{L}u - f^*\|^2_{L^2(\rho_I)} + \eta \|\mathcal{B}u - g^*\|^2_{L^2(\rho_B)}. 
\end{equation}
We restrict attention to
functions whose images under $\cL$ and $\cB$ belong to
$L^2(\rho_I)$ and $L^2(\rho_B)$, respectively. Define
\begin{equation}
\label{eq: domain or residual}
D_{\cR}
:=
\left\{
f\in D(\cL)\cap D(\cB):
\cL f\in L^2(\rho_I),\;
\cB f\in L^2(\rho_B)
\right\},
\end{equation}
where $D(\cL)$ and $D(\cB)$ denote the domains of the
PDE and boundary operators, respectively. For $\eta>0$, define the residual seminorm on
$D_{\mathcal R}$ by
\[
\|v\|_{\mathcal R_\eta}
:=
\left(
\|\mathcal L v\|_{L^2(\rho_I)}^2
+
\eta\|\mathcal B v\|_{L^2(\rho_B)}^2
\right)^{1/2},
\qquad v\in D_{\mathcal R}.
\]
This seminorm is a norm if the boundary value problem
is well posed in $D_{\mathcal R}$. Let $u^*$ be the solution of \eqref{problem}. Then the population
residual can be written as
\begin{equation*}
\mathcal{R}_{\eta}(u)
=
\|\mathcal{L}(u-u^*)\|_{L^2(\rho_I)}^2
+
\eta\|\mathcal{B}(u-u^*)\|_{L^2(\rho_B)}^2
=
\|u-u^*\|_{\mathcal{R}_{\eta}}^2.
\end{equation*}
To treat the interior and boundary components in a common
space, we introduce the following definition.
\begin{definition}
\label{def:interior-boundary-space}
Let $(\Omega,\Sigma_I,\rho_I)$ and
$(\partial\Omega,\Sigma_B,\rho_B)$ be the measure spaces
associated with the interior and boundary, respectively.
Define the $\sigma$-algebra on $\bar\Omega$ by
\[
\bar\Sigma
:=
\left\{
A\subseteq\bar\Omega:
A\cap\Omega\in\Sigma_I,\;
A\cap\partial\Omega\in\Sigma_B
\right\},
\]
and define the measure $\bar\rho$ on $\bar\Sigma$ by
\[
\bar\rho(A)
:=
\rho_I(A\cap\Omega)
+
\rho_B(A\cap\partial\Omega),
\qquad A\in\bar\Sigma.
\]
Set
\[
\bar L^2(\bar\Omega)
:=
L^2(\bar\Omega,\bar\Sigma,\bar\rho).
\]
Each $f\in\bar L^2(\bar\Omega)$ is identified with a pair
$(f_I,f_B)\in L^2(\rho_I)\oplus L^2(\rho_B)$ through
\[
f(\bm x)=
\begin{cases}
f_I(\bm x), & \bm x\in\Omega,\\
f_B(\bm x), & \bm x\in\partial\Omega,
\end{cases}
\]
with equality understood almost everywhere with respect
to the corresponding measures.
For $f=(f_I,f_B)$ and $g=(g_I,g_B)$, the inner product is
\[
\langle f,g\rangle_{\bar L^2(\bar\Omega)}
:=
\langle f_I,g_I\rangle_{L^2(\rho_I)}
+
\langle f_B,g_B\rangle_{L^2(\rho_B)}.
\]
The interior and boundary components are independent,
with no trace compatibility imposed.
\end{definition}

\subsection{Weighted Random Feature Model}
Let $\Theta\subset\mathbb{R}^{\tilde d}$ be a bounded
parameter space, let
$\phi:\overline{\Omega}\times\Theta\to\mathbb{R}$
be the feature function, and let
$\pi\in\mathcal{P}(\Theta)$ be a reference probability measure.
We draw $\bm{\theta}_1,\ldots,\bm{\theta}_m$
independently from a sampling distribution
$\pi'\in\mathcal{P}(\Theta)$. We impose
Assumption~\ref{sample distribution} on the pair
$(\nu,\nu')=(\pi,\pi')$ and write
$p=\mathrm{d}\pi'/\mathrm{d}\pi$.
The associated \emph{weighted random feature model}
is defined by
\begin{equation*}
\label{rf}
\sigmam
:=
\left\{
u_{\bm{\alpha}}:
u_{\bm{\alpha}}(\bm{x})
=
\frac{1}{\sqrt{m}}
\sum_{i=1}^{m}
\frac{\alpha_i}{\sqrt{p(\bm{\theta}_i)}}
\phi(\bm{x},\bm{\theta}_i),
\quad
\bm{\alpha}\in\mathbb{R}^m
\right\}.
\end{equation*}
The weights $p(\bm{\theta}_i)^{-1/2}$ account for
sampling from $\pi'$ rather than $\pi$, ensuring that $\mathbb{E}\left[\cS_{m,r}\cS_{m,r}^*\right]=\Phi_r\Phi_r^*$, where the expectation is over the parameters; see
\eqref{def: discrete_r} and
\eqref{def: Phi_r and Phi_r*}
for the definitions of $\cS_{m,r}$ and $\Phi_r$,
respectively.

\subsection{Collocation-Based Estimator}
\label{subsec:collocation-estimator}

RFM~\cite{chen2022bridging} constructs an empirical approximation of the population
residual $\mathcal{R}_{\eta}$ defined in \eqref{cont_loss}. Let
$\rho_I\in\mathcal{P}(\Omega)$ and
$\rho_B\in\mathcal{P}(\partial\Omega)$ be the reference probability
measures associated with the interior and boundary residuals,
respectively. We then independently draw the collocation points
according to
\(\{\bm{x}_i\}_{i=1}^{n_I}
\overset{\mathrm{i.i.d.}}{\sim}\rho_I^{\prime}
\) and 
\(
\{\bm{x}_{j+n_I}\}_{j=1}^{n_B}
\overset{\mathrm{i.i.d.}}{\sim}\rho_B^{\prime}.\) Here, $\rho_I^{\prime}\in\mathcal{P}(\Omega)$ and
$\rho_B^{\prime}\in\mathcal{P}(\partial\Omega)$ are the corresponding
sampling measures.
We apply Assumption~\ref{sample distribution} to the pairs
$(\rho_I,\rho_I^{\prime})$ on $\Omega$ and
$(\rho_B,\rho_B^{\prime})$ on $\partial\Omega$, with
\(
q_I:=\frac{\dd\rho_I^{\prime}}{\dd\rho_I}
\)
and 
\(
q_B:=\frac{\dd\rho_B^{\prime}}{\dd\rho_B}.
\)
We denote the total number of collocation
points by $n=n_I+n_B$. For simplicity, we consider the global RFM without a partition of
unity. We define the weighted empirical PDE residual by
\begin{equation} 
\label{discre_loss} 
\hat{\cR}_{n,\eta} (u) = \frac{1}{n_I} \sum_{i=1}^{n_I} \frac{\left( \mathcal{L}u(\bm{x}_i) - f^*(\bm{x}_i)\right)^2}{q_I(\bm{x}_i)} +  \frac{\eta}{n_B} \sum_{j=1}^{n_B} \frac{\left( \mathcal{B}u(\bm{x}_{j+n_I}) - g^*(\bm{x}_{j+n_I})\right)^2}{q_B(\bm{x}_{j+n_I})}. 
\end{equation}
For every fixed $u$, the empirical residual is unbiased:
\(
\mathbb E[\hat{\cR}_{n,\eta}(u)]
=
\mathcal R_\eta(u),
\)
where the expectation is over the collocation points. The approximate solution obtained by RFM is 
\begin{equation}
\label{eq:RFM solution}
\hat{u}_{n,m,\eta} := \argmin_{u_{\bm \alpha}\in \sigmam} \hat{\cR}_{n,\eta} (u_{\bm \alpha}). 
\end{equation}

\subsection{RKHS Characterization via the Feature Function}
We formulate our analysis in the RKHS characterized by the feature function $\phi$.
We impose the following assumption throughout the paper.
\begin{assumption}
\label{assump:basic-feature}
The feature function $\phi$ is jointly measurable in
$(x,\theta)$ and
satisfies
\[
\phi(\bm x,\cdot)\in L^2(\pi)
\qquad\text{for every }\bm x\in\Omega.
\]
Moreover,
\(
\int_\Theta
\|\phi(\cdot,\bm\theta)\|_{L^2(\rho_I)}^2
\,\mathrm d\pi(\bm\theta)<\infty.
\)\end{assumption}
For $\bm x,\bm x'\in\Omega$, define
\(
k(\bm x,\bm x^{\prime})
:=
\langle\phi(\bm x,\cdot),\phi(\bm x^{\prime},\cdot)\rangle_{L^2(\pi)}\), 
and denote the associated RKHS by \(\mathcal H_k\).
Define the feature operator  $\Phi: L^2(\pi)\mapsto L^2(\rho_I)$ by
\begin{equation}
\label{def: Phi_origin}
   \Phi \alpha = \int_{ \Theta}\alpha(\bm \theta)\phi(\cdot,\bm\theta)\dd\pi(\bm\theta). 
\end{equation}
% so the canonical inclusion
% $J:\mathcal H\to L^2(\rho_I)$ is bounded.
% Assume that $J$ is injective, and identify $\mathcal H$
% with its image in $L^2(\rho_I)$.
The feature-space characterization
of an RKHS \cite[Theorem~4.21]{steinwart2008support} yields
\begin{equation}
\label{eq: rkhs feature charact}
\mathcal H_k
=
\operatorname{Range}(\Phi)
=
\left\{\Phi a:a\in L^2(\pi)\right\}, \qquad
\|f\|_{\mathcal H_k}
=
\min_{\substack{\alpha\in L^2(\pi)\\\Phi \alpha=f}}
\|\alpha\|_{L^2(\pi)}. \end{equation}

Although a function may admit multiple representing
coefficients, its minimum-norm coefficient is unique
and belongs to \((\ker\Phi)^\perp\). Throughout this paper, we assume that the exact solution
$u^*$ to \eqref{problem} belongs to $\mathcal H_k$. We denote its unique minimum-norm representing
coefficient by \(b^*\). Thus,
\[
u^*=\Phi b^*,
\qquad
b^*\in(\ker\Phi )^\perp,
\qquad
\|b^*\|_{L^2(\pi)}=\|u^*\|_{\mathcal H_k}.
\]
\subsection{Additional Feature Assumptions}
\label{subsec:additional-feature-assumptions}

In addition to Assumption~\ref{assump:basic-feature},
we introduce two conditions for Sobolev approximation
and PDE analysis, respectively.

\begin{assumption}[Sobolev feature regularity]
\label{assump:sobolev-feature}
For the prescribed integer $r\geq0$, the map
$\bm\theta\mapsto\phi(\cdot,\bm\theta)$ is strongly
measurable with values in $H^r(\Omega)$ and satisfies
\[
\int_\Theta
\|\phi(\cdot,\bm\theta)\|_{H^r(\Omega)}^2
\,\mathrm d\pi(\bm\theta)<\infty.
\]
\end{assumption}

Under Assumption~\ref{assump:sobolev-feature},
for every $a\in L^2(\pi)$, the integral defining
$\Phi a$ in \eqref{def: Phi_origin} is well defined
as an $H^r(\Omega)$-valued Bochner integral. Moreover, $\Phi$, viewed as an operator from $L^2(\pi)$
to $H^r(\Omega)$, is Hilbert--Schmidt. It is used for the Sobolev approximation results with the specified order $r$.

For the PDE analysis, we impose the following compatibility
condition to justify interchanging $\cL$ and $\cB$ with
the feature integral. Specifically, for every $b\in L^2(\pi)$,
it ensures that
\[
\begin{aligned}
\cL(\Phi b)
&=
\int_\Theta
\cL\phi(\cdot,\bm\theta)b(\bm\theta)
\,\dd\pi(\bm\theta)
&&\text{in }L^2(\rho_I),\\
\cB(\Phi b)
&=
\int_\Theta
\cB\phi(\cdot,\bm\theta)b(\bm\theta)
\,\dd\pi(\bm\theta)
&&\text{in }L^2(\rho_B).
\end{aligned}
\]

\begin{assumption}[Compatibility with the PDE operators]
\label{assump:feature-operator-compatibility}
There exists a Banach function space
$V\subset D_{\mathcal R}$ continuously embedded in
$C(\overline\Omega)$ such that the map
$\bm\theta\mapsto\phi(\cdot,\bm\theta)$ is strongly
measurable with values in $V$ and satisfies
\[
\int_\Theta
\|\phi(\cdot,\bm\theta)\|_V^2
\,\mathrm d\pi(\bm\theta)<\infty.
\]
Furthermore,
\[
\mathcal L:V\to L^2(\rho_I),
\qquad
\mathcal B:V\to L^2(\rho_B)
\]
are bounded linear operators.
\end{assumption}
Assumption~\ref{assump:feature-operator-compatibility}
also implies that
$\mathcal H_k\subset V$ and $\sigmam\subset V$ almost surely when the
parameters are sampled from a distribution absolutely
continuous with respect to $\pi$. It is used for residual approximation and collocation
estimation. For the collocation analysis, we impose an additional
condition on the space $V$.
\begin{assumption}
\label{assump:collocation-evaluations}
Let $V$ be the space in
Assumption~\ref{assump:feature-operator-compatibility}.
For each $v\in V$, fix pointwise representatives of
$\cL v$ on $\Omega$ and $\cB v$ on $\partial\Omega$.
For every $\bm x\in\Omega$ and $\bm y\in\partial\Omega$,
the maps
\[
v\longmapsto(\cL v)(\bm x),
\qquad
v\longmapsto(\cB v)(\bm y)
\]
are bounded linear functionals on $V$.
\end{assumption}

\section{Main Results}
Our analysis separates two sources of error: approximation
error due to finitely many random features and estimation
error due to finitely many collocation observations.
We first establish upper and lower bounds for these errors
in terms of the eigenvalues of the corresponding operators.
We then establish a lower bound on the error under
simultaneous restrictions on the numbers of random features
and collocation observations and investigate how activation regularity affects these
eigenvalues.

For a measurable positive semidefinite kernel $\tilde k$
on $\Theta$, let $\widetilde{\mathcal T}:L^2(\pi)\to L^2(\pi)$
denote the associated integral operator,
\begin{equation*}
\label{eq:parameter-kernel-operator}
(\widetilde{\mathcal T}a)(\bm\theta)
:=
\int_\Theta
\tilde k(\bm\theta,\bm\theta')a(\bm\theta')
\,\mathrm d\pi(\bm\theta').
\end{equation*}
We use matching subscripts to distinguish the kernels
and their associated integral operators.

\label{sec:main}

% To avoid confusion between the parameters introduced below,
% we summarize their roles in Table~\ref{tab:regularization-parameters}.
% \begin{table}[t]
% \centering
% \caption{Summary of parameters and their roles.}
% \label{tab:regularization-parameters}
% \small
% \renewcommand{\arraystretch}{1.15}
% \begin{tabular}{@{}p{0.15\textwidth}p{0.80\textwidth}@{}}
% \toprule
% Parameter & Role \\
% \midrule
% $\eta>0$
% & Weights the boundary residual relative to the interior residual. \\
% \addlinespace
% $\lambda>0$
% & Weights the coefficient or RKHS norm penalty
% in the corresponding regularized problem. \\
% \addlinespace
% $\lambda_m$
% & Determines the feature-sampling density and controls
% the approximation-error bound. \\
% \addlinespace
% $\lambda_n$
% & Determines the collocation-sampling densities and controls
% the estimation-error bound. \\
% \bottomrule
% \end{tabular}
% \end{table}
\subsection{Approximation Error}
We first study the best approximation of functions in
the unit ball of $\mathcal H_k$ by $m$ random features.
We establish upper and lower bounds in Sobolev norms
and in the residual seminorm.

\paragraph{Approximation error in Sobolev norm.}
For the prescribed order $r$, we impose
Assumption~\ref{assump:sobolev-feature}. 
Define the kernel on $\Theta$ by
\[
\tilde k_r(\bm\theta,\bm\theta')
:=
\langle
\phi(\cdot,\bm\theta),\phi(\cdot,\bm\theta')
\rangle_{H^r(\Omega)}.
\]
Let $\{\tilde\mu_j^{(r)}\}_{j\geq1}$ denote the positive
eigenvalues of the associated integral operator
$\widetilde{\mathcal T}_r$, arranged in nonincreasing order. We first establish an upper bound for the worst-case
approximation error in $H^r(\Omega)$.
\begin{theorem}
\label{upper_appro}
Suppose that $\{\tilde\mu_j^{(r)}\}_{j\geq1}$ satisfies
Assumption~\ref{assump:eigenvalue-decay}.
There exist constants $c_0,m_0>0$, depending only on
$\{\tilde\mu_j^{(r)}\}_{j\geq1}$ and satisfying
$c_0m_0\geq1$, such that the following holds.

For $\delta\in(0,1)$, let $m$ satisfy $\frac{m}{\log(14m/\delta)}\geq m_0$. Denote $l_m:=
\left\lfloor
\frac{c_0m}{\log(14m/\delta)}
\right\rfloor$ and $p_m(\bm \theta)=q^*(\tilde{k}_r,\pi,\tilde\mu_{l_m}^{(r)};\bm \theta)$.
If
$\bm\theta_1,\ldots,\bm\theta_m$ are sampled independently
from $p_m\,\dd\pi$, then, with probability at least $1-\delta$,
\begin{equation*}
\label{eq: upper-sobolev-effective-index}
\sup_{\|u^*\|_{\mathcal H_k}\leq1}
\inf_{u\in\sigmam}
\|u-u^*\|_{H^r(\Omega)}^2
\le 16
\tilde\mu_{l_m}^{(r)}.
\end{equation*}
\end{theorem}
We next establish a lower bound that holds for every fixed
choice of feature parameters.
\begin{theorem}
\label{lower_appro}
For every fixed collection of sampled parameters
$\bm\theta_1,\ldots,\bm\theta_m$ in the parameter space $\Theta$,
\begin{equation}
\label{eq: lower-sobolev}
\sup_{\|u^*\|_{\mathcal H_k}\leq1}
\inf_{u\in\sigmam}
\|u-u^*\|_{H^r(\Omega)}^2
\geq \tilde\mu_{m+1}^{(r)}.
\end{equation}
\end{theorem}

\paragraph{Approximation error in residual seminorm.}
We next consider approximation error measured by the residual seminorm. We impose Assumption~\ref{assump:feature-operator-compatibility}
and define the kernel
\[
\begin{aligned}
\tilde k_{\mathcal R_\eta}(\bm\theta,\bm\theta')
:={}&
\langle
\mathcal L\phi(\cdot,\bm\theta),
\mathcal L\phi(\cdot,\bm\theta')
\rangle_{L^2(\rho_I)}+\eta
\langle
\mathcal B\phi(\cdot,\bm\theta),
\mathcal B\phi(\cdot,\bm\theta')
\rangle_{L^2(\rho_B)}.
\end{aligned}
\]
Let $\{\tilde\mu_j^{(\cR_\eta)}\}_{j\geq1}$ denote
the positive eigenvalues of the associated integral
operator $\widetilde{\mathcal T}_{\cR_\eta}$,
arranged in nonincreasing order. We establish the following upper bound and lower bound for the approximation error in the residual seminorm.
\begin{theorem}
\label{upper_appro_residual}
Suppose that $\{\tilde\mu_j^{(\cR_\eta)}\}_{j\geq1}$ satisfies
Assumption~\ref{assump:eigenvalue-decay}.
There exist constants $c_0,m_0>0$, depending only on
$\{\tilde\mu_j^{(\cR_\eta)}\}_{j\geq1}$ and satisfying
$c_0m_0\geq1$, such that the following holds.

For $\delta\in(0,1)$, let $m$ satisfy
$\frac{m}{\log(14m/\delta)}\geq m_0$.
Denote $l_m:=
\left\lfloor
\frac{c_0m}{\log(14m/\delta)}
\right\rfloor$ and $p_m(\bm \theta)=q^*(\tilde{k}_{\mathcal R_\eta},\pi,\tilde\mu_{l_m}^{(\mathcal R_\eta)};\bm \theta)$. If
$\bm\theta_1,\ldots,\bm\theta_m$ are sampled independently
from $p_m\,\dd\pi$, then, with probability at least $1-\delta$,
\begin{equation}
\label{eq: upper-residual-effective-index}
\sup_{\|u^*\|_{\mathcal H_k}\leq1}
\inf_{u\in\sigmam}
\|u-u^*\|_{\cR_\eta}^2
\leq 16\tilde\mu_{l_m}^{(\cR_\eta)}.
\end{equation}
\end{theorem}

\begin{theorem}
\label{lower_appro_residual}
For every fixed collection of sampled parameters
$\bm\theta_1,\ldots,\bm\theta_m$ in the parameter space $\Theta$,
\begin{equation}
\label{eq: lower-residual}
\sup_{\|u^*\|_{\mathcal H_k}\leq1}
\inf_{u\in\sigmam}
\|u-u^*\|_{\cR_\eta}^2
\geq \tilde\mu_{m+1}^{(\cR_\eta)}.
\end{equation}
\end{theorem}
The approximation error analysis is deferred to
Section~\ref{approximation analysis}.
\subsection{Estimation Error}

We now study the estimation error associated with finitely
many collocation observations.
Throughout this subsection, we impose
Assumption~\ref{assump:feature-operator-compatibility}. Let
\[
X_n:=(\bm x_1,\ldots,\bm x_n)
\in\Omega^{n_I}\times(\partial\Omega)^{n_B}, \quad
n=n_I+n_B,
\]
be a collection of collocation points. Define the observation vector by
\[
\bm z_n(X_n,u^*)
:=
\left(
\bigl(\mathcal L u^*(\bm x_i)\bigr)_{i=1}^{n_I},
\bigl(\mathcal B u^*(\bm x_i)\bigr)_{i=n_I+1}^{n}
\right)^{\!\top}
\in\mathbb R^n.
\]
Our objective is to recover $u^*$ from the data $S_n(u^*):=(X_n,\bm z_n(X_n,u^*))$.
For $\lambda>0$, consider the regularized estimator based on the data:
\begin{equation}
\label{est_problem}
\hat u_{n,\eta,\lambda}
:=
\argmin_{u\in\mathcal H_k}
\left\{
\hat{\cR}_{n,\eta}(u)
+\lambda\|u\|_{\mathcal H_k}^2
\right\},
\end{equation}
where $\hat{\cR}_{n,\eta}$ is the weighted empirical
PDE residual defined in \eqref{discre_loss}. Here $\eta$ weights the boundary term in the empirical residual
$\hat{\mathcal R}_{n,\eta}$, whereas $\lambda$ weights the squared
RKHS norm $\|u\|_{\mathcal H_k}^2$. Let $n=n_I+n_B$ and set $\nu_I:=n/n_I$, $\nu_B:=n/n_B$, and
$\nu_+:=\max\{\nu_I,\nu_B\}$.
Define the kernel on $\Theta$ by
\[
\begin{aligned}
\tilde k_{\mathrm{PDE}}(\bm\theta,\bm\theta')
:={}&
\nu_I
\left\langle
\mathcal L\phi(\cdot,\bm\theta),
\mathcal L\phi(\cdot,\bm\theta')
\right\rangle_{L^2(\rho_I)}+\eta\nu_B
\left\langle
\mathcal B\phi(\cdot,\bm\theta),
\mathcal B\phi(\cdot,\bm\theta')
\right\rangle_{L^2(\rho_B)}.
\\
\end{aligned}
\]
Let $\{\mu_j^{\mathrm{PDE}}\}_{j\geq1}$ denote the
positive eigenvalues of the associated integral operator
$\widetilde{\cT}_{\mathrm{PDE}}$, arranged in
nonincreasing order. To specify the interior and boundary sampling densities
$q_I$ and $q_B$ in the following theorem, define
\begin{equation}
\label{eq: k_L}
k_{\cL}(\bm x,\bm x')
:=
\int_\Theta
\cL\phi(\bm x,\bm\theta)\,
\cL\phi(\bm x',\bm\theta)\,
\dd\pi(\bm\theta),
\qquad \bm x,\bm x'\in\Omega,
\end{equation}
and
\begin{equation}
\label{eq: k_B}
k_{\cB}(\bm x,\bm x')
:=
\int_\Theta
\cB\phi(\bm x,\bm\theta)\,
\cB\phi(\bm x',\bm\theta)\,
\dd\pi(\bm\theta),
\qquad \bm x,\bm x'\in\partial\Omega.
\end{equation}

We establish an upper bound for the worst-case
estimation error in the residual seminorm.

\begin{theorem}
\label{thm: Estimation_error upper}
Fix $\eta>0$ and the constants $\nu_I$ and $\nu_B$.
Suppose that the eigenvalues $\{\mu_j^{\mathrm{PDE}}\}_{j\geq1}$ satisfy one of the
eigenvalue decay conditions in
Lemma~\ref{lemma: N(gamma)}.
There exist constants $c_0,n_0>0$, depending only on $\{\mu_j^{\mathrm{PDE}}\}_{j\geq1}$
and the fixed sampling ratios $\nu_I,\nu_B$, and satisfying $c_0n_0\geq1$, such that the
following holds.

For $\delta\in(0,1)$, let $n$ satisfy
$\frac{n}{\log(28n/\delta)}\geq n_0$. Denote \(
l_n:=
\left\lfloor
\frac{c_0n}{\log(28n/\delta)}
\right\rfloor
\), $q_I= q^*(k_\cL,\rho_I,\mu_{l_n}^{\mathrm{PDE}};\bm x)$ and $q_B= q^*(k_\cB,\rho_B,\mu_{l_n}^{\mathrm{PDE}}/\eta;\bm x)$.
If $\bm x_1,\ldots,\bm x_{n_I}$ are sampled independently
from $q_I\,\mathrm d\rho_I$ and
$\bm x_{n_I+1},\ldots,\bm x_n$ are sampled independently
from $q_B\,\mathrm d\rho_B$, with the two collections
independent, then, with probability at least $1-\delta$,
the estimator $\hat u_{n,\eta,\lambda}$ defined in
\eqref{est_problem} satisfies the following bound
simultaneously for every $0<\lambda\leq\mu_{l_n}^{\mathrm{PDE}}$:
\[
\sup_{\|u^*\|_{\mathcal H_k}\leq1}
\|\hat u_{n,\eta,\lambda}-u^*\|_{\cR_\eta}^2
\leq 4\mu_{l_n}^{\mathrm{PDE}}.
\]
\end{theorem}

We next establish a lower bound for every fixed
collection of collocation points.
Define the class of all estimators
based on the collocation points and the corresponding
PDE observations by
\begin{equation}
\label{M_n pde}
\mathcal M_n^{\mathrm{PDE}}
:=
\left\{
G:
\bigl(\Omega^{n_I}\times(\partial\Omega)^{n_B}\bigr)
\times\mathbb R^n
\to D_{\mathcal R}
\right\}. 
\end{equation}

\begin{theorem}
\label{thm: Estimation_error lower}
For every fixed collection $X_n$ of collocation
points, let $S_n(u^*)=(X_n,\bm z_n(X_n,u^*))$. Then we have
\begin{equation*}
\label{eq:collocation-lower}
\begin{aligned}
&\inf_{G\in\mathcal M_n^{\mathrm{PDE}}}
\sup_{\|u^*\|_{\mathcal H_k}\leq1}
\|G(S_n(u^*))-u^*\|^2_{\cR_\eta}
\geq
\frac{\mu_{n+1}^{\mathrm{PDE}}}{\nu_+}.
\end{aligned}
\end{equation*}
\end{theorem}
We now consider the combined restrictions imposed by
finitely many random features and collocation observations. Consider the class of estimators in which the outputs are restricted to the random feature model
\begin{equation}
\label{M_mn}
\mathcal M_{m,n}
:=
\left\{
G:
\bigl(\Omega^{n_I}\times(\partial\Omega)^{n_B}\bigr)
\times\mathbb R^n\to\sigmam\right\}.
\end{equation}

Define the optimal worst-case squared residual error by
\begin{equation}
\label{eq: joint-residual-error}
E_{m,n}
:=
\inf_{G\in\mathcal M_{m,n}}
\sup_{\|u^*\|_{\mathcal H_k}\leq1}
\|G(S_n(u^*))-u^*\|_{\cR_\eta}^2.
\end{equation}
\begin{corollary}
\label{thm: joint-residual-lower}
For every fixed admissible choice of feature parameters
$\{\bm\theta_i\}_{i=1}^m$ and collocation points $X_n$,
we have
\begin{equation}
\label{eq: joint-residual-lower}
 E_{m,n}
\geq
\frac{1}{\nu_+}
\max\left\{
\mu_{m+1}^{\mathrm{PDE}},
\mu_{n+1}^{\mathrm{PDE}}
\right\}.
\end{equation}
\end{corollary}
The RFM solution in \eqref{eq:RFM solution}, viewed as
a function of $S_n(u^*)$, defines an estimator
$G\in\mathcal M_{m,n}$, with a fixed selection rule
based only on these data if the minimizer is not unique.
Therefore, Corollary~\ref{thm: joint-residual-lower} yields
\[
\sup_{\|u^*\|_{\mathcal H_k}\leq1}
\|\hat u_{n,m,\eta}-u^*\|_{\cR_\eta}^2
\geq E_{m,n}
\geq
\frac{1}{\nu_+}
\max\left\{
\mu_{m+1}^{\mathrm{PDE}},
\mu_{n+1}^{\mathrm{PDE}}
\right\}.
\]
The estimation error analysis is deferred to
Section~\ref{sec: estimation}.

\begin{remark}
\label{rem: residual-to-sobolev}
The estimation error analyzed in this section
is measured in the residual seminorm. Upper bounds in other norms can be
obtained from stability estimates for the underlying boundary
value problem. For example, let $\Omega\subset\mathbb R^d$ be a bounded domain
with $C^\infty$ boundary, and consider
\begin{equation}
\label{example PDE}
  \mathcal L u
=
-\sum_{i,j=1}^d a^{ij}(\bm{x})\,\partial_{ij}u
+\sum_{i=1}^d b^i(\bm{x})\,\partial_i u
+c(\bm{x})u,
\qquad
\mathcal B u=\operatorname{Tr}u,  
\end{equation}
where the coefficients $a^{ij},b^i,c\in
C^\infty(\overline{\Omega})$ and satisfy
\[
\sum_{i,j=1}^d a^{ij}(\bm{x})\,\xi_i\xi_j
\geq \alpha|\bm{\xi}|^2,
\qquad
\bm{x}\in\overline{\Omega},\quad
\bm{\xi}\in\mathbb R^d,
\]
for some $\alpha>0$. Assume that the homogeneous Dirichlet
problem
\[
\mathcal L v=0 \quad\text{in }\Omega,
\qquad
\mathcal B v=0 \quad\text{on }\partial\Omega
\]
has only the trivial solution. The stability estimate in \cite[Appendix~B.1, Eq.~(B2)]{shin2023error} gives
\[
\|v\|_{H^{1/2}(\Omega)}^2
\lesssim
\|\mathcal L v\|_{L^2(\Omega)}^2
+
\|\mathcal B v\|_{L^2(\partial\Omega)}^2,
\qquad v\in H^2(\Omega).
\]
Suppose that $\rho_I$ and $\rho_B$ have densities bounded below
by positive constants with respect to Lebesgue measure on
$\Omega$ and surface measure on $\partial\Omega$, respectively.
For an exact solution $u^*\in H^2(\Omega)$ satisfying
$\mathcal L u^*=f^*$ and $\mathcal B u^*=g^*$, every
$u\in H^2(\Omega)$ therefore satisfies
\[
\|u-u^*\|_{H^{1/2}(\Omega)}^2
\lesssim
\|\mathcal L u-f^*\|_{L^2(\rho_I)}^2
+
\|\mathcal B u-g^*\|_{L^2(\rho_B)}^2.
\]
Consequently, upper bounds for the residual seminorm yield
upper bounds of the same order for the squared
$H^{1/2}(\Omega)$ error, and hence for the squared
$L^2(\Omega)$ error.
\end{remark}
%\subsection{Adaptation to the Regularity of Activation Functions}
\subsection{Eigenvalue Decay Induced by Activation Regularity}
To obtain explicit convergence rates, we
estimate the eigenvalue decay of the integral operator associated
with the chosen feature function. We consider the commonly used
features
\begin{equation}
\label{feature}
\phi(\bm{x},\bm{\theta})
=
\sigma(\bm{k}^{\top}\bm{x}+b),
\end{equation}
where $\sigma:\mathbb R\to\mathbb R$ is the activation function
and $\bm{\theta}=(\bm{k},b)$ belongs to a compact parameter set
$\Theta\subset\mathbb R^d\times\mathbb R$.
The residual approximation and collocation estimation
bounds in
Theorems~\ref{upper_appro_residual}--\ref{thm: Estimation_error lower}
apply to general linear differential and boundary
operators under their respective stated assumptions.
For fixed proportions of interior and boundary observations,
\eqref{mu_PDE mu_R_eta} implies
that $\tilde\mu_j^{(\cR_\eta)}\asymp\mu_j^{\mathrm{PDE}}$
uniformly in $j$. It therefore suffices to estimate
the PDE eigenvalues.
To estimate the PDE eigenvalues, we specialize to the
elliptic boundary value problem in \eqref{example PDE}. The same approach applies to other linear equations.
The following theorem bounds the Sobolev and PDE
eigenvalues in terms of the regularity of $\sigma$.
\begin{theorem}
\label{activation_fun_rate}
Assume that $\rho_I$ is a probability measure on $\Omega$
with density $w_I\in L^\infty(\Omega)$ with respect to
Lebesgue measure.
Let $\rho_B$ be a probability measure on $\partial\Omega$. Let $\mathcal L$ be the elliptic operator in
\eqref{example PDE}, with
$a^{ij},b^i,c\in L^\infty(\Omega)$, and let
$\mathcal B u=u|_{\partial\Omega}$.
Fix $\eta>0$ and the constants $\nu_I$ and $\nu_B$.
For $\bm\theta=(\bm k,b)\in\Theta$, set
\(
\phi(\bm x,\bm\theta)
=\sigma(\bm k^\top\bm x+b),
\)
and choose
\(
M>
\max_{\bm x\in\overline\Omega,\;\bm\theta\in\Theta}
|\bm k^\top\bm x+b|.
\)

% For a fixed integer $r\geq0$, let
% $\{\tilde\mu_j^{(r)}\}_{j\geq1}$ denote the eigenvalues of $\widetilde{\mathcal T}_r=
% \Phi^*_{r}\Phi_{r}$. And let $\{\mu_j^{\mathrm{PDE}}\}_{j\geq1}$ denote the
% eigenvalues of
% \(
% \widetilde{\mathcal T}_{\mathrm{PDE}}
% :=
% \Phi_{\mathrm{PDE}}^*\Phi_{\mathrm{PDE}}
% =
% \nu_I\Phi_{\cL}^*\Phi_{\cL}
% +\lambda\nu_B\Phi_{\cB}^*\Phi_{\cB},
% \).
\begin{enumerate}
\item[(i)]
If
$\sigma\in H^p([-M,M])\cap C^{p-1}([-M,M])$
for an integer $p\geq\max\{2,r\}$, then
\[
\tilde\mu_j^{(r)}
=o\!\left(j^{-\frac{2(p-r)}{d}-1}\right),
\qquad
\mu_j^{\mathrm{PDE}}
=o\!\left(j^{-\frac{2(p-2)}{d}-1}\right).
\]
\item[(ii)]
If $\sigma$ is real analytic on an open interval
containing $[-M,M]$, then there exist constants
$C>0$ and $\tau\in(0,1)$
such that
\[
\tilde\mu_j^{(r)}
\leq C\tau^{2j^{1/d}},
\qquad
\mu_j^{\mathrm{PDE}}
\leq C\tau^{2j^{1/d}}.
\]
\end{enumerate}
\end{theorem}

\begin{remark}
Under the assumptions of Theorem~\ref{activation_fun_rate},
if $\sigma\in C^\infty([-M,M])$, then, for every fixed
integer $r\geq0$ and every $t>0$,
\[
\tilde\mu_m^{(r)}=o(m^{-t}),
\qquad
\mu_n^{\mathrm{PDE}}=o(n^{-t}).
\]
Thus, both eigenvalue sequences decay super-algebraically.
Smoothness alone does not imply the stretched-exponential
bounds in the analytic case.
\end{remark}

For the rectified power activation
$\mathrm{ReLU}^k(t):=(\max\{t,0\})^k$,
where $k\geq1$ is an integer, we obtain the following
sharper upper bounds.
\begin{theorem}
\label{thm: reluk-spectral-decay}
Let $\sigma=\mathrm{ReLU}^k$ for an integer $k\geq1$.
For each integer $0\leq r\leq k$, the eigenvalues
$\{\tilde\mu_j^{(r)}\}_{j\geq1}$ satisfy
\begin{equation}
\label{eq: reluk-eigenvalue-decay}
\tilde\mu_j^{(r)}
\leq Cj^{-1-\frac{2(k-r)+1}{d}},\qquad j\geq1.
\end{equation}
For $k\geq2$, consider the second-order elliptic operator
$\mathcal L$ and the Dirichlet boundary operator
$\mathcal B u=u|_{\partial\Omega}$ specified in
Theorem~\ref{activation_fun_rate}.
The eigenvalues $\{\mu_j^{\mathrm{PDE}}\}_{j\geq1}$ satisfy
\[
\mu_j^{\mathrm{PDE}}
\leq Cj^{-1-\frac{2k-3}{d}},\qquad j\geq1.
\]
\end{theorem}
The eigenvalue bound in
Theorem~\ref{thm: reluk-spectral-decay}, together with
Theorem~\ref{upper_appro}, yields the following
approximation rate.
For each integer $0\leq r\leq k$ and fixed
$\delta\in(0,1)$, a sampling density can be chosen such
that, for all sufficiently large $m$, with probability
at least $1-\delta$,
\[
\sup_{\|u^*\|_{\mathcal H_k}\leq1}
\inf_{u\in\sigmam}
\|u-u^*\|_{H^r(\Omega)}
\leq
C_\delta
\left(\frac{m}{\log m}\right)^{
-\frac12-\frac{2(k-r)+1}{2d}}.
\]
Unlike \cite{liu2025integral}, our analysis allows
general compact parameter sets rather than requiring
the feature parameters to lie on the unit sphere.
When $\Theta=\mathbb S^d$ and $\pi$ is uniform,
the RKHS $\mathcal H_k$ coincides with
$H^{(d+2k+1)/2}(\Omega)$ up to equivalent norms;
see \cite[Theorem~2.3]{liu2025integral}.
Our bound therefore recovers the approximation rate
in \cite[Theorem~6.2]{liu2025integral}, possibly under
a different sampling density, and extends this rate
to more general compact parameter sets.
The proofs of Theorems~\ref{activation_fun_rate}
and~\ref{thm: reluk-spectral-decay} are deferred to
Appendices~\ref{proof: activation_fun_rate}
and~\ref{proof: thm: reluk-spectral-decay}, respectively.

\section{Approximation and Estimation Analysis}
\label{Approximation and Estimation Analysis}
This section develops the analysis underlying the
approximation and estimation bounds in
Section~\ref{sec:main}. We first develop a duality framework in Hilbert spaces
that establishes an isometric correspondence between function
spaces on the computational and parameter domains,
yielding an identity between the corresponding
approximation errors.
We then use this framework to derive approximation bounds
in the Sobolev norm and residual seminorm through kernel ridge
regression on the parameter domain.
Finally, we analyze estimation errors and combine the approximation and estimation lower bounds
to quantify the joint limitations imposed by finitely
many random features and collocation observations.
\subsection{A Duality Framework in Hilbert Spaces}
\label{subsec: framework}
To analyze finite-feature approximation in different error
norms, we develop a duality framework in Hilbert
spaces.
The standard operator-theoretic framework for kernel methods is
typically formulated in the $L^2$ setting, as in \cite[Appendix~B.2]{long2025optimalratessaturationnoiseless}.
To accommodate the norms arising in PDEs, we extend this
framework to general Hilbert spaces, which provides the foundation
for our subsequent analysis. Let $\mathcal X$ be a nonempty set, and let the ambient space
$B$ be a separable Hilbert space of real-valued
functions on $\mathcal X$. The choice of $B$ will depend on the error metric under
consideration. Let $\Theta$ be the parameter space equipped with the
reference distribution $\pi$. To ensure that the associated feature operator is Hilbert--Schmidt, we impose
the following condition.
\begin{assumption}
\label{assump:general-feature}
The map
\(
    \bm{\theta}\mapsto\phi(\cdot,\bm{\theta})
\)
is strongly measurable as a $B$-valued map and satisfies
\[
    \int_{\Theta}
    \|\phi(\cdot,\bm{\theta})\|_{B}^2
    \dd\pi(\bm{\theta})<\infty.
\]
Equivalently, $\phi\in L^2(\pi;B)$
in the Bochner sense.
\end{assumption}
Define the feature operator $\Phi:L^2(\pi)\to B$ by
\begin{equation}
\label{eq: abstract feature operator}
    (\Phi\alpha)(\bm x)
    = \langle \phi(\bm x,\cdot),\alpha\rangle_{L^2(\pi)},
    \qquad \alpha\in L^2(\pi).
\end{equation}
Its adjoint $\Phi^*:B\to L^2(\pi)$ is given by
\begin{equation}
\label{eq: abstract adjoint feature operator}
    (\Phi^*f)(\bm\theta)
    = \langle \phi(\cdot,\bm\theta),f\rangle_B,
    \qquad f\in B.
\end{equation}
Assumption~\ref{assump:general-feature} implies that $\Phi^*$ is
Hilbert--Schmidt. Indeed, for any orthonormal basis
$\{e_j\}_{j=1}^{\infty}$ of $B$, Tonelli's theorem and Parseval's
identity yield
\[
    \|\Phi^*\|_{\mathrm{HS}}^2
    =
    \sum_{j=1}^{\infty}\|\Phi^* e_j\|_{L^2(\pi)}^2
    =
    \int_{\Theta}
    \|\phi(\cdot,\bm\theta)\|_{B}^2
    \dd\pi(\bm\theta)
    <\infty.
\]
Hence, $\Phi^*$ is compact and admits a singular value decomposition:
\begin{equation}
\label{eqn: svd-Phi}
\Phi^* = \sum_{j=1}^{\infty} 
\mu_j^{1/2}\widetilde{e}_j\otimes e_j,
\end{equation}
where $\{\widetilde e_j\}_{j=1}^{\infty}$ and
$\{e_j\}_{j=1}^{\infty}$ are the corresponding orthonormal bases of $L^2(\pi)$ and $ B$, respectively.
Consequently, the operators
\(
\mathcal T
:=
\Phi\Phi^*:B\to B
\)
and
\(
\widetilde{\mathcal T}
:=
\Phi^*\Phi:L^2(\pi)\to L^2(\pi)
\)
are self-adjoint, positive, and compact. Their spectral
decompositions are given by
\[
\mathcal T
=
\sum_{j=1}^{\infty}
\mu_j e_j\otimes e_j,
\qquad
\widetilde{\mathcal T}
=
\sum_{j=1}^{\infty}
\mu_j\widetilde e_j\otimes\widetilde e_j.
\] We can therefore define the interpolation spaces on $\mathcal X$ and
$\Theta$, analogously to \eqref{eqn: Hks-space2}, as follows:
\begin{equation}
\label{eqn: Hks-space-x}
\mathcal H^s
=
\left\{
\sum_{j=1}^{\infty}a_j\mu_j^{s/2}e_j:
\sum_{j=1}^{\infty}a_j^2<\infty
\right\},
\qquad
\left\|
\sum_{j=1}^{\infty}a_j\mu_j^{s/2}e_j
\right\|_{\mathcal H^s}^2
=
\sum_{j=1}^{\infty}a_j^2,
\end{equation}
and
\begin{equation}
\label{eqn: Hks-space-theta}
\widetilde{\mathcal H}^s
=
\left\{
\sum_{j=1}^{\infty}a_j\mu_j^{s/2}\widetilde e_j:
\sum_{j=1}^{\infty}a_j^2<\infty
\right\},
\qquad
\left\|
\sum_{j=1}^{\infty}a_j\mu_j^{s/2}\widetilde e_j
\right\|_{\widetilde{\mathcal H}^s}^2
=
\sum_{j=1}^{\infty}a_j^2.
\end{equation}
By definition, \(\mathcal H^0=B\) and \(\widetilde{\mathcal H}^0=L^2(\pi)\). The following lemma shows that $\Phi^*$ induces an isometric
isomorphism between $\mathcal H^{s-1}$ and
$\widetilde{\mathcal H}^{s}$.
\begin{lemma}
\label{lemma RKHS equivalence}
Let $s\geq1$. For any $g\in \widetilde{\mathcal{H}}^s$, there exists a unique $u\in \cH^{s-1}$ such that $\Phi^* u=g$ and $\|g\|_{\widetilde{\mathcal{H}}^s} = \|u\|_{\cH^{s-1}}$. 

\end{lemma}
\begin{proof}
Following the argument in the proof of \cite[Proposition B.2]{long2025optimalratessaturationnoiseless}, let 
$$u_g=\sum_{j=1}^{\infty}\mu_j^{-1/2} \langle g,\widetilde{e}_j\rangle_{L^2(\pi)}e_j.$$ It follows from the singular value decomposition of $\Phi^*$ as in \eqref{eqn: svd-Phi} that 
\[
\Phi^* u_g = \sum_{j=1}^{\infty} \langle g,\widetilde{e}_j\rangle_{L^2(\pi)}\widetilde{e}_j =g,
\]
and 
\[
\|u_g\|^2_{\cH^{s-1}}=\sum_{j=1}^{\infty}\mu_j^{-(s-1)}\mu_j^{-1} \langle g,\widetilde{e}_j\rangle^2_{L^2(\pi)} = \|g\|^2_{\widetilde{\mathcal{H}}^s}.
\]
It suffices to justify the uniqueness. Suppose $u=\sum_{j=1}^{\infty} u_je_j \in \cH^{s-1}  $ such that $\Phi^* u=\sum_{j=1}^{\infty} \mu_j^{1/2} u_j \widetilde{e}_j=g$. Since $\mu_j>0$ for all $j\ge 1$, we must have $ u_j= \mu_j^{-1/2}\langle g,\widetilde{e}_j\rangle_{L^2(\pi)}$ for all $j$. Consequently, $u=u_g$.
\end{proof}
The following lemma relates the approximation errors in
the two spaces, provided that both error norms
are finite.
\begin{lemma}
\label{lemma: approximate duality revised}
Let $s\geq 1$ and $u^*\in\mathcal H^{s-1}$, and set
$g^*=\Phi^*u^*\in\widetilde{\mathcal H}^{s}$.
For any $\hat u\in B$, define $\hat g=\Phi^*\hat u$.
Then, for every $p\in[0,s]$, we have
\begin{equation}
\label{eq: approximation error identity}
\|\hat g-g^*\|_{\widetilde{\mathcal H}^{p}}^2
=
\|\hat u-u^*\|_{\mathcal H^{p-1}}^2
\end{equation}
whenever \(\hat u-u^*\in\mathcal H^{p-1}\) and
\(\hat g-g^*\in\widetilde{\mathcal H}^{p}\).

\end{lemma}

\begin{proof}
The inclusion $g^*\in\widetilde{\mathcal H}^{s}$ follows
from Lemma~\ref{lemma RKHS equivalence}.
Set
\[
b:=\hat u-u^*,
\qquad
h:=\hat g-g^*=\Phi^* b.
\]
Write $b=\sum_{j=1}^{\infty}b_j e_j$, with convergence
in $\mathcal H^{p-1}$, and define
\[
b_N:=\sum_{j=1}^N b_j e_j,
\qquad
h_N:=\Phi^* b_N
=\sum_{j=1}^N\mu_j^{1/2}b_j\widetilde e_j.
\]
By the spectral representation of $\Phi^*$, $h_N$ is the
spectral truncation of $h$. It holds that 
\[
b_N\to b\quad\text{in }\mathcal H^{p-1},
\qquad
h_N\to h\quad\text{in }\widetilde{\mathcal H}^{p}.
\]
For each $N$, 
\begin{align*}
\|h_N\|_{\widetilde{\mathcal H}^{p}}^2
&=
\|\widetilde{\mathcal T}^{-p/2}\Phi^*b_N\|_{L^2(\pi)}^2=
\left\langle
b_N,\Phi\widetilde{\mathcal T}^{-p}\Phi^*b_N
\right\rangle_B\\
&\stackrel{(i)}{=}
\left\langle b_N,\mathcal T^{1-p}b_N\right\rangle_B=
\|\mathcal T^{-(p-1)/2}b_N\|_B^2
=
\|b_N\|_{\mathcal H^{p-1}}^2,
\end{align*}
where $(i)$ follows from the identity
\[
\Phi\widetilde{\mathcal T}^{-p}\Phi^* e_j
=\mu_j^{1/2}\Phi\widetilde{\mathcal T}^{-p}\widetilde e_j=
\mu_j^{1/2}\mu_j^{-p}\Phi\widetilde e_j
=
\mu_j^{1-p}e_j
=
\mathcal T^{1-p}e_j,
\]
which holds for each $j$.
Passing to the limit as $N\to\infty$ proves
\eqref{eq: approximation error identity}.
\end{proof}
\begin{remark}
Although Lemma~\ref{lemma: approximate duality revised}
is stated for general $s$ and $p$, the subsequent analyses
of the approximation error in both the Sobolev and residual
norms use only the case $s=2$ and $p=1$.
\end{remark}
The following lemma shows that the interpolation space
$\cH^1$ does not depend on the choice of
the ambient space $B$.
\begin{lemma}
\label{lem:ambient-invariance}
Let $B_1$ and $B_2$ be real Hilbert spaces continuously
embedded in a common Hilbert space $G$, and identify them
with their images in $G$.
Suppose that the same feature map
$\theta\mapsto\phi(\cdot,\theta)$ belongs to
$L^2(\pi;B_i)$ for $i=1,2$.
Define $\Phi_i^*:B_i\to L^2(\pi)$ by
\[
(\Phi_i^*f)(\bm\theta)
=
\langle\phi(\cdot,\bm\theta),f\rangle_{B_i},
\]
and let $\Phi_i:L^2(\pi)\to B_i$ be its adjoint.
Let $\mathcal H_i^1$ be the interpolation space of order $1$ associated with $\Phi_i\Phi_i^*$, as defined in \eqref{eqn: Hks-space-x}.
Then
\[
\mathcal H_i^1=\operatorname{Range}(\Phi_i),
\qquad
\|f\|_{\mathcal H_i^1}
=
\min_{\substack{\alpha\in L^2(\pi)\\\Phi_i\alpha=f}}
\|\alpha\|_{L^2(\pi)}.
\]
Moreover, as subspaces of $G$,
\[
\mathcal H_1^1=\mathcal H_2^1,
\qquad
\|f\|_{\mathcal H_1^1}=\|f\|_{\mathcal H_2^1}.
\]
\end{lemma}

\begin{proof}
By the singular value decomposition of $\Phi_i$,
\begin{equation}
\label{eq: character of inter space}
\mathcal H_i^1=\operatorname{Range}(\Phi_i),
\qquad
\|f\|_{\mathcal H_i^1}
=
\min_{\substack{\alpha\in L^2(\pi)\\\Phi_i\alpha=f}}
\|\alpha\|_{L^2(\pi)},
\qquad i=1,2. 
\end{equation}

Since the embeddings $B_i\hookrightarrow G$ are continuous,
the Bochner integral representations give, in $G$,
\[
\Phi_1\alpha
=
\int_\Theta
\phi(\cdot,\bm\theta)\alpha(\bm\theta)\,\dd\pi(\bm\theta)
=
\Phi_2\alpha
\qquad
\text{for every }\alpha\in L^2(\pi).
\]
Consequently,
\[
\mathcal H_1^1
=
\operatorname{Range}(\Phi_1)
=
\operatorname{Range}(\Phi_2)
=
\mathcal H_2^1.
\]
For any $f\in\mathcal H_1^1=\mathcal H_2^1$,
\[
\|f\|_{\mathcal H_1^1}
=
\min_{\substack{\alpha\in L^2(\pi)\\\Phi_1\alpha=f}}
\|\alpha\|_{L^2(\pi)}
=
\min_{\substack{\alpha\in L^2(\pi)\\\Phi_2\alpha=f}}
\|\alpha\|_{L^2(\pi)}
=
\|f\|_{\mathcal H_2^1}.
\]
\end{proof}

\begin{remark}[Connection with RKHS]
In the setting of Lemma~\ref{lem:ambient-invariance},
suppose additionally that $B_1$ and $B_2$ consist of
real-valued functions on $\mathcal X$.
% and that
% \[
% (\Phi_i\alpha)(\bm x)
% =
% \langle\alpha,\phi(\bm x,\cdot)\rangle_{L^2(\pi)},
% \qquad i=1,2,
% \]
% for every $\alpha\in L^2(\pi)$ and $\bm x\in\mathcal X$,
% where $\phi(\bm x,\cdot)\in L^2(\pi)$.
% Define
% \[
% k(\bm x,\bm x')
% =
% \langle\phi(\bm x,\cdot),\phi(\bm x',\cdot)\rangle_{L^2(\pi)}.
% \]
The operators $\Phi_i$, $i=1,2$, are defined in the same
way as $\Phi$ in \eqref{def: Phi_origin}.
Moreover, the characterization of the interpolation space
of order $1$ in \eqref{eq: character of inter space}
coincides with that of the RKHS $\cH_k$
in \eqref{eq: rkhs feature charact}.
Consequently, 
\[
\mathcal H_1^1=\mathcal H_k=\mathcal H_2^1,
\qquad
\|f\|_{\mathcal H_1^1}
=
\|f\|_{\mathcal H_k}
=
\|f\|_{\mathcal H_2^1}.
\]
\end{remark}
% We focus on the approximation error of the random feature model, optimal error of the continuous loss function when using the random feature model and the estimation error over the whole RKHS based on the discrete loss function using $n$ data $\{(x_i,f(x_i))\}_{i=1}^{n_I}\cup \{(x_j,g(x_j))\}_{j=1+n_I}^{n_I+n_B}$.

% We investigate two sources of error in the random feature method:
% (i) the approximation error of the random feature model, measured in
% both the Sobolev norm and the residual norm; and
% (ii) the estimation error over the entire RKHS $\cH$, arising from
% the weighted empirical residual defined in \eqref{discre_loss} and
% evaluated using the data samples
% ${(\bm{x}_i,f^*(\bm{x}_i))}_{i=1}^{n_I}
% \cup
% {(\bm{x}_j,g^*(\bm{x}_j))}_{j=1+n_I}^{n_I+n_B}$.
\subsection{Approximation Error Analysis}
\label{approximation analysis}
In this subsection, we analyze finite-feature approximation
in Sobolev norms and the residual seminorm using the duality
framework developed in Section~\ref{subsec: framework}.
The key idea is to interpret feature selection as sampling
on the parameter space: the adjoint feature operator maps
features to kernel sections, while their inner products
with the target become point evaluations of a transformed
target. This correspondence yields a KRR representation of the regularized approximator and
an exact identity between the corresponding errors. We first study the approximation error of the random feature model in Sobolev norms.
For $r\in \mathbb{N}$ and $u^*\in \cH_k$, we consider the following linear least-squares problem 
\begin{equation*}
    \min_{u_{\bm\alpha}\in \sigmam}\|u_{\bm\alpha}-u^*\|^2_{H^r(\Omega)}.
\end{equation*}
The above optimization problem is equivalent to the following finite-dimensional optimization problem over the coefficient vector \(\bm{\alpha}\in\mathbb{R}^m\):
\begin{align*}
\min_{\bm{\alpha} \in \mathbb{R}^m} \left(\bm{\alpha}^{\top} \bm{K}_m \bm{\alpha} -2\bm{y}^{\top} \bm{\alpha} +c\right), \label{def y}
\end{align*}
where $c=\|u^*\|_{H^r(\Omega)}^2$ and 
\begin{equation}
\label{def: K_m y}
\bm{K}_m=\left(\frac{\left\langle \phi(\cdot,\bm{\theta}_i), \phi(\cdot,\bm{\theta}_j)\right\rangle_{H^r(\Omega)}}{m\sqrt{p(\bm \theta_i)p(\bm \theta_j)}}\right)_{i,j=1}^m
, \quad \bm{y}=\left(\frac{\left\langle\phi(\cdot,\bm{\theta}_i), u^*(\cdot)\right\rangle_{H^r(\Omega)}}{\sqrt{mp(\bm\theta_i)}}\right)_{i=1}^m.
\end{equation}
Since the Gram matrix \(\bm K_m\) may be singular, the coefficient vector representing the best approximation may be non-unique. To obtain a unique coefficient representation and facilitate the subsequent error analysis, we introduce a regularization parameter \(\lambda>0\) and consider
\begin{equation}
\label{approximate func}
    \hat{u}^{(r)}_{m,\lambda} = \argmin_{u_{\bm\alpha}\in  \sigmam} \left(\|u_{\bm\alpha}-u^*\|^2_{H^r(\Omega)} +  \lambda \|\bm\alpha\|_{2}^2\right).
\end{equation}
Define $\hat{\bm{\phi}}_m : \Omega \to \mathbb{R}^m $ with $(\hat{\bm{\phi}}_m(\cdot))_i =\frac{\phi(\cdot,\bm{\theta}_i)}{\sqrt{mp(\bm\theta_i)}} $. Then $\hat{u}^{(r)}_{m,\lambda} = \hat{\bm{\phi}}_m(\cdot)^{\top} \hat{\bm{\alpha}}^{(r)}_{m,\lambda}$, where 
\begin{equation}
    \label{rf_coef}
    \hat{\bm{\alpha}}^{(r)}_{m,\lambda} = \argmin_{\bm{\alpha} \in \mathbb{R}^m}  \left(\bm{\alpha}^{\top} (\bm{K}_m+ \lambda \bm{I}_m) \bm{\alpha} -2\bm{y}^{\top} \bm{\alpha} +c\right)
    = (\bm{K}_m+\lambda \bm{I}_m)^{-1}\bm{y}.   
\end{equation}
We next express the approximant in operator form. Let $\cS_{m,r}: \RR^m \to H^r(\Omega)$ be the discrete operator and let
$\cS^\ast_{m,r}: H^r(\Omega) \to \RR^m $ denote its adjoint. For $\bm{a}\in\RR^m$ and $f\in H^r(\Omega)$, these operators
are given by
\begin{equation}
\label{def: discrete_r}
\begin{aligned}
\cS_{m,r}\bm{a}
= \frac{ 1}{\sqrt{m}}\sum_{i=1}^m \frac{a_i}{\sqrt{p(\bm \theta_i)}}\phi(\cdot,\bm{\theta}_i), \quad 
(\cS_{m,r}^*f)_j
&= \frac{ \langle \phi(\cdot,\bm{\theta}_j),f(\cdot)\rangle_{H^r(\Omega)}}{\sqrt{mp(\bm\theta_j)}},
   \quad j=1,\ldots,m.\\
\end{aligned}
\end{equation}
One can easily verify that 
\begin{equation}
\label{appro func operator form}
\hat{u}^{(r)}_{m,\lambda} = \cS_{m,r} (\cS^*_{m,r}\cS_{m,r} +\lambda \bm{I}_m )^{-1} \cS^*_{m,r} u^*.
\end{equation}
The corresponding continuous operator $\Phi_r: L^2(\pi) \to H^r(\Omega)$ and its adjoint 
$\Phi_r^*$ are obtained by
specializing \eqref{eq: abstract feature operator} and
\eqref{eq: abstract adjoint feature operator}, respectively,
to $B=H^r(\Omega)$. Explicitly, for $\alpha\in L^2(\pi)$ and $f\in H^r(\Omega)$,  \begin{equation}
\label{def: Phi_r and Phi_r*}
(\Phi_r\alpha)(\bm{x})
= \langle \phi(\bm{x},\cdot),\alpha\rangle_{L^2(\pi)}, \quad 
(\Phi_r^* f)(\bm{\theta})
= \langle \phi(\cdot,\bm{\theta}),f\rangle_{H^r(\Omega)}.
\end{equation}

To connect this formulation with kernel ridge regression
(KRR), we introduce the kernel on the parameter space $\Theta$:
\begin{equation}
\label{new_k}
\tilde{k}_r(\bm{\theta},\bm{\theta}^{\prime})
= \left\langle
    \phi(\cdot,\bm{\theta}),
    \phi(\cdot,\bm{\theta}^{\prime})
  \right\rangle_{H^r(\Omega)}.
\end{equation}
The corresponding integral operator $\tilde{\mathcal{T}}_r$ is given by
\begin{equation}
\label{int_oper_appro}
    \tilde{\mathcal{T}}_r g = \int_{\Theta} \tilde{k}_r (\cdot,\bm{\theta}^{\prime}) g(\bm{\theta}^{\prime}) \dd \pi(\bm{\theta}^{\prime})=\Phi^*_{r}\Phi_{r}g.
\end{equation}
% Assume that the feature function $\phi$ is sufficiently smooth, then it is reasonable that 
% \[
% \sup_{\Theta}  \tilde{k}_r(\theta,\theta) < \infty,
% \]
% then
Its spectral decomposition is
\begin{equation}\label{mercer's decomposition_dual}
\tilde{\mathcal{T}}_r = \sum_{j=1}^\infty \tilde{\mu}^{(r)}_j \tilde{e}^{(r)}_j\otimes \tilde{e}^{(r)}_j,
\end{equation}
where $\{\tilde{\mu}^{(r)}_j\}_{j=1}^\infty$ are the eigenvalues in a decreasing order and $\{\tilde{e}^{(r)}_j\}_{j=1}^\infty$ are the corresponding orthonormal eigenfunctions in $L^2(\pi)$. Let $\cH_{(r)}^s$ and $\widetilde{\mathcal H}_{(r)}^s$ denote the
interpolation spaces defined in \eqref{eqn: Hks-space-x}
and \eqref{eqn: Hks-space-theta}, respectively, with
$B=H^r(\Omega)$. The RKHS associated with $\tilde k_r$ is isometrically
identified with $\widetilde{\mathcal H}_{(r)}^1$.

% Note that
% \begin{equation*}
% \begin{aligned}
%         y_i=\frac{1}{\sqrt{m}}\left\langle\phi(\cdot,\theta_i), u^*(\cdot)\right\rangle_{H^r(\Omega)} &= \frac{1}{\sqrt{m}}\int_{\Omega} \sum_{|\beta|\leq r} D^{\beta}_x \phi(x,\theta_i)\cdot D^{\beta}_x \left(\int_{\Theta}\phi(x,\theta) b^*(\theta) \dd \pi(\theta)\right)\dd x \\
%         &= \frac{1}{\sqrt{m}}\int_{\Omega} \sum_{|\beta|\leq r} D^{\beta}_x \phi(x,\theta_i)\cdot  \int_{\Theta} D^{\beta}_x \phi(x,\theta) b^*(\theta) \dd \pi(\theta)\dd x \\
%         & \stackrel{(i)}{=}\frac{1}{\sqrt{m}}  \int_{\Theta} 
%  \left(\int_{\Omega} \sum_{|\beta|\leq r} D^{\beta}_x \phi(x,\theta_i)\cdot D^{\beta}_x \phi(x,\theta)\dd x\right) b^*(\theta)  \dd \pi(\theta) \\
%         &= \frac{1}{\sqrt{m}} \tilde{\mathcal{T}}_r b^* (\theta_i) ,
% \end{aligned}
% \end{equation*}
% where $(i)$ uses the Fubini's theorem.
% Let $g^*= \tilde{\mathcal{T}}_r b^*$, the vector $y$ can be interpreted as the evaluation of $g^*$ at the parameters $\{\theta_i\}_{i=1}^m$. Note that $g^* \in \tilde{\mathcal{H}}_r$, 

Now consider the following KRR on $\Theta$ with input points $\{\bm{\theta}_i\}_{i=1}^{m}\overset{\text{i.i.d.}}{\sim} \pi^{\prime}$ and target function $g^*(\bm{\theta})= \left\langle\phi(\cdot,\bm{\theta}), u^*(\cdot)\right\rangle_{H^r(\Omega)}=(\Phi^*_{r} u^*)(\bm{\theta})$. It follows that $g^*\in\widetilde{\mathcal H}_{(r)}^1$. The KRR estimator is obtained by solving:
\begin{equation}\label{eqn: krr-0}
\hat{g}^{(r)}_{m,\lambda} = \argmin_{g\in \widetilde{\mathcal{H}}^1_{(r)}} \left(\frac{1}{m}\sum_{i=1}^m \frac{(g(\bm{\theta}_i)-g^*(\bm \theta_i))^2}{p(\bm \theta_i)} + \lambda \|g\|_{\widetilde{\mathcal{H}}^1_{(r)}}^2\right).
\end{equation}  
Define $\hat{\bm k}_m:\Theta\to\mathbb R^m$ and
$\tilde{\bm K}_m\in\mathbb R^{m\times m}$ by
\[
(\hat{\bm k}_m(\cdot))_i
=\frac{\tilde k_r(\bm\theta_i,\cdot)}{\sqrt{mp(\bm\theta_i)}},
\qquad
\tilde{\bm K}_m
=\left(
\frac{\tilde k_r(\bm\theta_i,\bm\theta_j)}
{m\sqrt{p(\bm\theta_i)p(\bm\theta_j)}}
\right)_{i,j=1}^m.
\] The KRR solution can then be written as $\hat{g}^{(r)}_{m,\lambda}=\hat{\bm{k}}_m(\cdot)^{\top}\hat{\bm{b}}^{(r)}_{m,\lambda}$ with 
\begin{equation}
\label{para_krr_sol1}
    \hat{\bm{b}}^{(r)}_{m,\lambda} = (\tilde{\bm{K}}_m+\lambda \bm{I}_m)^{-1} \bm{y},
\end{equation}
where $\bm y$ is defined in \eqref{def: K_m y}. 
By \eqref{new_k}, $\tilde{\bm K}_m=\bm K_m$, so
$\hat{\bm b}^{(r)}_{m,\lambda}=\hat{\bm\alpha}^{(r)}_{m,\lambda}$. Moreover, for any $\hat{\bm{b}}\in \mathbb{R}^m$,
\begin{align*}
    \hat{\bm{k}}_m(\bm{\theta})^{\top}\hat{\bm{b}} &= \frac{1}{\sqrt{m}}\sum_{i=1}^m \frac{\hat{b}_i}{\sqrt{p(\bm \theta_i)}}  \left\langle \phi(\cdot,\bm{\theta}), \phi(\cdot,\bm{\theta}_i)\right\rangle_{H^r(\Omega)} \\
    &= \left\langle \phi(\cdot,\bm{\theta}), \frac{1}{\sqrt{m}}\sum_{i=1}^m \frac{\hat{b}_i}{\sqrt{p(\bm \theta_i)}} \phi(\cdot,\bm{\theta}_i)\right\rangle_{H^r(\Omega)} \\
    &=\left\langle \phi(\cdot,\bm{\theta}), \cS_{m,r} \hat{\bm{b}}\right\rangle_{H^r(\Omega)} =(\Phi_r^* \cS_{m,r} \hat{\bm{b}})(\theta).
\end{align*}
Taking $\hat{\bm b}=\hat{\bm b}^{(r)}_{m,\lambda}$ yields
\begin{equation}
\label{express g_lambda}
\hat g^{(r)}_{m,\lambda}
=
\Phi_r^*\cS_{m,r}
(\cS_{m,r}^*\cS_{m,r}+\lambda\bm I_m)^{-1}
\cS_{m,r}^*u^*
=
\Phi_r^*\hat u^{(r)}_{m,\lambda}.
\end{equation}

Let $N^{(r)}$ and $F_q^{(r)}$ denote the DoFs obtained
from \eqref{eq: dof} by replacing $(\mu_j,e_j)$ with
$(\tilde\mu_j^{(r)},\tilde e_j^{(r)})$ and taking $\pi$ as the reference measure on $\Theta$.
We now combine the duality identity with the KRR error
bound in \cite[Theorem~4.3]{long2025optimalratessaturationnoiseless} to prove Theorem~\ref{upper_appro}.
% Suppose that
% $\{\tilde\mu_j^{(r)}\}_{j=1}^{\infty}$ satisfies one of the three
% eigenvalue decay conditions in Lemma~\ref{lemma: N(gamma)}.  Then there exist constants $c_0,m_0>0$ such that
% \[
% N^{(r)}\!\left(
% \tilde\mu_{\lfloor c_0t\rfloor}^{(r)}
% \right)
% \leq \frac{t}{20}
% \qquad \text{for all } t\geq m_0,
% \]
% and increase $m_0$, if necessary, so that $c_0m_0\geq1$.

% For $\delta\in(0,1)$ and each positive integer $m$ satisfying
% \begin{equation}
% \label{eq: sobolev-sample-size}
% \frac{m}{\log(14m/\delta)}\geq m_0,
% \end{equation}
% define
% \begin{equation}
% \label{eq: density sobolev}
% l_m:=
% \left\lfloor
% \frac{c_0m}{\log(14m/\delta)}
% \right\rfloor,
% \qquad
% \gamma_m:=\tilde\mu_{l_m}^{(r)},
% \qquad
% q_m:=\bar q_{1,\gamma_m}^{*,(r)}.   
% \end{equation}

\begin{proof}[Proof of Theorem \ref{upper_appro}]
By Lemma~\ref{lemma: N(gamma)},
$N^{(r)}(\tilde\mu_k^{(r)})\leq Ck$ for all $k\geq1$. Choosing $c_0\leq(20C)^{-1}$ and $m_0$ sufficiently
large so that $c_0m_0\geq1$, we obtain
\[
N^{(r)}\!\left(\tilde\mu_{\lfloor c_0t\rfloor}^{(r)}\right)
\leq C\lfloor c_0t\rfloor
\leq \frac{t}{20},
\qquad t\geq m_0.
\]
Write $L_m:=\log(14m/\delta)$ and $\lambda_m:=\tilde\mu_{\lfloor c_0m/L_m\rfloor}^{(r)}$. For $m/L_m\ge m_0$,
\[
N^{(r)}(\lambda_m)
=
N^{(r)}
\!\left(\tilde\mu_{\lfloor c_0m/L_m\rfloor}^{(r)}\right)
\leq\frac{m}{20L_m}.
\]
By \eqref{eq: control F}, we have
\[
F_{p_m}^{(r)}(\lambda_m)
\leq2N^{(r)}(\lambda_m)
\leq\frac{m}{10L_m}.
\]
Since $L_m>1$, it follows that $F_{p_m}^{(r)}(\lambda_m)
\leq m$ and
\[
5F_{p_m}^{(r)}(\lambda_m)
\max\left\{
1,\log\frac{14F_{p_m}^{(r)}(\lambda_m)}{\delta}
\right\}
\leq5F_{p_m}^{(r)}(\lambda_m)L_m
\leq\frac m2 < m.
\]
Thus the critical penalty
$\lambda_m^{\mathrm{crit}}$ associated with $p_m$,
as defined in
\cite[Definition~4.2]{long2025optimalratessaturationnoiseless},
satisfies
\[
\lambda_m^{\mathrm{crit}}\leq\lambda_m.
\]
Fix $\lambda\in(0,\lambda_m^{\mathrm{crit}}]$.
Applying \cite[Theorem~4.3]{long2025optimalratessaturationnoiseless}
with $s=2$ and $p=1$ yields, with probability at least
$1-\delta$,
\[
\sup_{\|g^*\|_{\widetilde{\mathcal H}_{(r)}^2}\leq1}
\|\hat g_{m,\lambda}-g^*\|_{
\widetilde{\mathcal H}_{(r)}^1}^2
\leq16\lambda_m^{\mathrm{crit}}.
\]
Since $g^*=\Phi_r^*u^*$ and $\hat g^{(r)}_{m,\lambda}=\Phi_r^*\hat u_{m,\lambda}^{(r)}$, 
Lemma~\ref{lemma: approximate duality revised} gives
\begin{align*}
\sup_{\|u^*\|_{\mathcal H_k}\leq1}
\inf_{u\in\sigmam}
\|u-u^*\|_{H^r(\Omega)}^2
&\leq
\sup_{\|u^*\|_{\mathcal H_k}\leq1}
\|\hat u_{m,\lambda}^{(r)}-u^*\|_{H^r(\Omega)}^2\\
&=
\sup_{\|g^*\|_{\widetilde{\mathcal H}_{(r)}^2}\leq1}
\|\hat g^{(r)}_{m,\lambda}-g^*\|_{\widetilde{\mathcal H}_{(r)}^1}^2\\
&\leq16\lambda_m^{\mathrm{crit}}
\leq16\lambda_m
=16\tilde\mu_{\lfloor c_0m/L_m\rfloor}^{(r)}.
\end{align*}
\end{proof}

\begin{remark}[Direct and dual representations of the approximation error]
\label{rem:direct-dual-approximation}
Let $u^*=\Phi_r b^*$ with $b^*\in L^2(\pi)$, and set
\[
\mathcal T_{m,r}:=\cS_{m,r}\cS_{m,r}^*,
\qquad
R_{m,\lambda}:=(\mathcal T_{m,r}+\lambda I)^{-1}.
\]
By \eqref{appro func operator form}, the direct approximation error is
$\hat u^{(r)}_{m,\lambda}-u^*=-\lambda R_{m,\lambda}u^*$.
Since $R_{m,\lambda}$ is self-adjoint on $H^r(\Omega)$,
the operator $A_{m,\lambda}:=-\lambda R_{m,\lambda}\Phi_r$
from $L^2(\pi)$ to $H^r(\Omega)$ has adjoint
$A_{m,\lambda}^*=-\lambda\Phi_r^*R_{m,\lambda}$.
Thus, by \eqref{express g_lambda} and $g^*=\Phi_r^*u^*$,
\[
\hat u^{(r)}_{m,\lambda}-u^*
= A_{m,\lambda}b^*,
\qquad
\hat g^{(r)}_{m,\lambda}-g^*
= \Phi_r^*(\hat u^{(r)}_{m,\lambda}-u^*)
= A_{m,\lambda}^*u^*.
\]
By Lemma~\ref{lemma: approximate duality revised}, these two
representations give the same error in their respective norms,
\[
\|\hat u^{(r)}_{m,\lambda}-u^*\|_{H^r(\Omega)}
=
\|\hat g^{(r)}_{m,\lambda}-g^*\|_{\widetilde{\mathcal H}_{(r)}^1}.
\]
This identity also gives an operator-norm interpretation of the
worst-case error. We have
\[
\|\hat g^{(r)}_{m,\lambda}-g^*\|_{\widetilde{\mathcal H}_{(r)}^1}
=
\|A_{m,\lambda}b^*\|_{H^r(\Omega)}
=
\sup_{\|v\|_{H^r(\Omega)}\leq 1}
\bigl|\langle b^*,A_{m,\lambda}^*v\rangle_{L^2(\pi)}\bigr|.
\]
By the feature representation of the RKHS unit ball, taking
the supremum over $\|u^*\|_{\mathcal H_k}\leq 1$ is equivalent
to taking the supremum over $\|b^*\|_{L^2(\pi)}\leq 1$.
Consequently,
\[
\begin{aligned}
\sup_{\|u^*\|_{\mathcal H_k}\leq 1}
\|\hat u^{(r)}_{m,\lambda}-u^*\|_{H^r(\Omega)}
&= \sup_{\|b\|_{L^2(\pi)}\leq 1}
   \|A_{m,\lambda}b\|_{H^r(\Omega)} \\
&= \|A_{m,\lambda}\|_{L^2(\pi)\to H^r(\Omega)}, \\[1ex]
\sup_{\|g^*\|_{\widetilde{\mathcal H}_{(r)}^2}\leq 1}
\|\hat g^{(r)}_{m,\lambda}-g^*\|_{\widetilde{\mathcal H}_{(r)}^1}
&= \sup_{\substack{\|b\|_{L^2(\pi)}\leq 1\\
                   \|v\|_{H^r(\Omega)}\leq 1}}
   \bigl|\langle b,A_{m,\lambda}^*v\rangle_{L^2(\pi)}\bigr| \\
&= \sup_{\|v\|_{H^r(\Omega)}\leq 1}
   \|A_{m,\lambda}^*v\|_{L^2(\pi)} \\
&= \|A_{m,\lambda}^*\|_{H^r(\Omega)\to L^2(\pi)}.
\end{aligned}
\]
Since $\|A_{m,\lambda}^*\|=\|A_{m,\lambda}\|$, the bound
on $\|A_{m,\lambda}^*\|$ which is established in KRR analysis also bounds the operator norm that
governs the direct approximation error. Thus, both analyses
control the same worst-case error and yield the same convergence
rate under the stated sampling assumptions. We use the dual representation because it places finite-feature
approximation and collocation estimation in a common noiseless
KRR framework. 
\end{remark}
We next prove the lower bound on the worst-case
Sobolev approximation error stated in
Theorem~\ref{lower_appro}.
\begin{proof}[Proof of Theorem \ref{lower_appro}]
Define
\[
V_m
:=
\operatorname{span}
\{\tilde k_r(\cdot,\bm\theta_i):i=1,\ldots,m\}
\subset\widetilde{\mathcal H}_{(r)}^1.
\]
Since for any $\bm \theta\in \Theta$,
\[
\Phi_r^*(\phi(\cdot,\bm\theta_i))(\bm \theta)
=
\tilde k_r(\bm \theta,\bm\theta_i),
\]
we have $\Phi_r^*\sigmam=V_m$.
By
Lemma~\ref{lemma: approximate duality revised},
\begin{align*}
\sup_{\|u^*\|_{\mathcal H_k}\leq1}
\inf_{u\in\sigmam}
\|u-u^*\|_{H^r(\Omega)}^2=
\sup_{\|g^*\|_{\widetilde{\mathcal H}_{(r)}^2}\leq1}
\inf_{h\in V_m}
\|h-g^*\|_{\widetilde{\mathcal H}_{(r)}^1}^2.
\end{align*}
Let
\[
S_m:=(\bm\theta_1,\ldots,\bm\theta_m),
\qquad
g(S_m):=(g(\bm\theta_1),\ldots,g(\bm\theta_m))^\top,
\]
and define
\[
\mathcal M_m
:=
\left\{
\mathcal Q:
\mathcal Qg=G(S_m,g(S_m)),\
G:\Theta^m\times\mathbb R^m\to\mathbb R^\Theta
\right\}.
\]
Let $\mathcal Q_0$ denote the orthogonal projection onto
$V_m$ in $\widetilde{\mathcal H}_{(r)}^1$. Then
\[
\sup_{\|g^*\|_{\widetilde{\mathcal H}_{(r)}^2}\leq1}
\inf_{h\in V_m}
\|h-g^*\|_{\widetilde{\mathcal H}_{(r)}^1}^2
=
\sup_{\|g^*\|_{\widetilde{\mathcal H}_{(r)}^2}\leq1}
\|\mathcal Q_0g^*-g^*\|_{\widetilde{\mathcal H}_{(r)}^1}^2.
\]
By the reproducing property, $\mathcal Q_0g^*$ depends
only on $S_m$ and $g^*(S_m)$, so
$\mathcal Q_0\in\mathcal M_m$. Then
\begin{align*}
\sup_{\|g^*\|_{\widetilde{\mathcal H}_{(r)}^2}\leq1}
\|\mathcal Q_0g^*-g^*\|_{\widetilde{\mathcal H}_{(r)}^1}^2&\geq
\inf_{\mathcal Q\in\mathcal M_m}
\sup_{\|g^*\|_{\widetilde{\mathcal H}_{(r)}^2}\leq1}
\|\mathcal Qg^*-g^*\|_{\widetilde{\mathcal H}_{(r)}^1}^2\\
&\geq\tilde\mu_{m+1}^{(r)},
\end{align*}
where the last inequality follows from the case
$s=2$, $p=1$ of
\cite[Theorem~4.1]{long2025optimalratessaturationnoiseless}.
\end{proof}

\begin{remark}[A geometric interpretation of the Sobolev approximation lower bounds]
Write $\mu_j=\tilde\mu_j^{(r)}$. In the equivalent KRR problem on
the parameter space, the target unit ball, expressed in coordinates
orthonormal in the error space $\widetilde{\mathcal H}_{(r)}^1$, is
the ellipsoid
\[
\mathcal E:=
\left\{
\sum_{j\geq1}\sqrt{\mu_j}\,b_j\psi_j:
\sum_{j\geq1}|b_j|^2\leq1
\right\},
\]
where $\psi_j$ is the $j$th principal direction. Retaining the first
$m$ directions gives worst-case squared error
\[
\sup_{\sum_j|b_j|^2\leq1}\sum_{j>m}\mu_j|b_j|^2=\mu_{m+1},
\]
with equality attained at the target $\sqrt{\mu_{m+1}}\psi_{m+1}$. 
This choice of directions is optimal. Indeed, for any trial space
$V$ of dimension at most $m$, the span of $\psi_1,\ldots,\psi_{m+1}$
contains a unit vector $w$ orthogonal to $V$. Since the first $m+1$
semiaxes all have length at least $\sqrt{\mu_{m+1}}$, the target
$g^*=\sqrt{\mu_{m+1}}\,w$ belongs to $\mathcal E$ and satisfies
\[
\inf_{v\in V}\|g^*-v\|_{\widetilde{\mathcal H}_{(r)}^1}^2
=\|g^*\|_{\widetilde{\mathcal H}_{(r)}^1}^2
=\mu_{m+1}.
\]
The space generated by $m$ sampled features has dimension at most
$m$, so its worst-case squared approximation error is at least
$\mu_{m+1}$.
\end{remark}

We next study approximation in the residual seminorm
$\|\cdot\|_{\cR_\eta}$.
For the exact solution $u^*\in\mathcal H_k$ of
\eqref{problem}, the squared approximation error is
\begin{equation}
\label{eq: approximation residual}
\inf_{u_{\bm\alpha}\in\sigmam}
\|u_{\bm\alpha}-u^*\|_{\cR_\eta}^2
=
\inf_{u_{\bm\alpha}\in\sigmam}
\left(
\|\cL u_{\bm\alpha}-f^*\|_{L^2(\rho_I)}^2
+\eta\|\cB u_{\bm\alpha}-g^*\|_{L^2(\rho_B)}^2
\right).
\end{equation}
For $\eta>0$, define the residual operator
$P_{\cR_\eta}:D_{\cR}\to\bar L^2(\bar\Omega)$ by
\begin{equation}
\label{def residual_mapping}
(P_{\cR_\eta}f)(\bm x)
=
\begin{cases}
\cL f(\bm x),
& \bm x\in\Omega,\\
\sqrt{\eta}\,\cB f(\bm x),
& \bm x\in\partial\Omega,
\end{cases}
\end{equation}
with equality understood almost everywhere with respect
to the corresponding measures. Define the residual
feature function by
\begin{equation}
\label{def residual_feature}
\phi_{\cR_\eta}(\bm x,\bm\theta)
:=
\bigl(P_{\cR_\eta}\phi(\cdot,\bm\theta)\bigr)(\bm x)
=
\begin{cases}
\cL\phi(\bm x,\bm\theta),
& \bm x\in\Omega,\\
\sqrt{\eta}\,\cB\phi(\bm x,\bm\theta),
& \bm x\in\partial\Omega.
\end{cases}
\end{equation}
Assumption~\ref{assump:feature-operator-compatibility}
implies that the map $\bm\theta\mapsto
\phi_{\cR_\eta}(\cdot,\bm\theta)
$ belongs to $L^2(\pi;\bar L^2(\bar\Omega))$. Thus, Assumption~\ref{assump:general-feature}
holds with $B=\bar L^2(\bar\Omega)$ and feature function
$\phi_{\cR_\eta}$.
The optimization problem \eqref{eq: approximation residual} is equivalent to the following finite-dimensional optimization problem over the coefficient vector \(\bm{\alpha}\in\mathbb{R}^m\):
\begin{align*}
\inf_{\bm{\alpha} \in \mathbb{R}^m} \left(\bm{\alpha}^{\top} \bm{K}^{\cR_\eta}_m \bm{\alpha} -2(\bm{y}^{\cR_\eta})^{\top} \bm{\alpha} +c^{\cR_\eta}\right), 
\end{align*}
where $c^{\cR_\eta}=\|u^*\|_{\cR_\eta}^2$ and 
\begin{align*}
\bm{K}^{\cR_\eta}_m &=\left(\frac{\left\langle \cL\phi(\cdot,\bm{\theta}_i), \cL\phi(\cdot,\bm{\theta}_j)\right\rangle_{L^2(\rho_I)}+\eta\langle\cB\phi(\cdot,\bm{\theta}_i), \cB\phi(\cdot,\bm{\theta}_j)\rangle_{L^2(\rho_B)}}{m\sqrt{p(\bm \theta_i)p(\bm \theta_j)}}\right)_{i,j=1}^m
 \\
 &=\left(\frac{\langle\phi_{\cR_{\eta}}(\cdot,\bm \theta_i),\phi_{\cR_{\eta}}(\cdot,\bm \theta_j)\rangle_{\bar L^2(\bar\Omega)}}{m\sqrt{p(\bm \theta_i)p(\bm \theta_j)}}\right)_{i,j=1}^m,\\
\bm{y}^{\cR_\eta}&=\left(\frac{\left\langle\cL\phi(\cdot,\bm{\theta}_i), f^*(\cdot)\right\rangle_{L^2(\rho_I)}+\eta\langle \cB \phi(\cdot,\bm \theta_i),g^*\rangle_{L^2(\rho_B)}}{\sqrt{mp(\bm\theta_i)}}\right)_{i=1}^m \\
&=\left(\frac{\langle\phi_{\cR_{\eta}}(\cdot,\bm \theta_i),P_{\cR_{\eta}}u^*(\cdot)\rangle_{\bar L^2(\bar\Omega)}}{\sqrt{mp(\bm\theta_i)}}\right)_{i=1}^m.
\end{align*}
As in \eqref{approximate func}, we consider
\begin{equation*}
\label{approximate func residual}
    \hat{u}^{\cR_\eta}_{m,\lambda} = \argmin_{u_{\bm\alpha}\in  \sigmam} \left(\|u_{\bm\alpha}-u^*\|^2_{\cR_\eta} +  \lambda \|\bm\alpha\|_{2}^2\right).
\end{equation*}
Let $\hat{\bm{\phi}}_m : \Omega \mapsto \mathbb{R}^m $ with $(\hat{\bm{\phi}}_m(\cdot))_i =\frac{\phi(\cdot,\bm{\theta}_i)}{\sqrt{mp(\bm\theta_i)}} $. Then $\hat{u}^{\cR_\eta}_{m,\lambda} = \hat{\bm{\phi}}_m(\cdot)^{\top} \hat{\bm{\alpha}}^{\cR_\eta}_{m,\lambda}$, where the coefficient vector $\hat{\bm{\alpha}}^{\cR_\eta}_{m,\lambda}$ is determined by
\begin{equation*}
    \label{rf_coef_residual}
    \hat{\bm{\alpha}}^{\cR_\eta}_{m,\lambda} = \argmin_{\bm{\alpha} \in \mathbb{R}^m}  \left(\bm{\alpha}^{\top} (\bm{K}^{\cR_\eta}_m+ \lambda \bm{I}_m) \bm{\alpha} -2(\bm{y}^{\cR_\eta})^{\top} \bm{\alpha} +c^{\cR_\eta}\right)
    = (\bm{K}^{\cR_\eta}_m+\lambda \bm{I}_m)^{-1}\bm{y}^{\cR_\eta}.   
\end{equation*}
% To facilitate the subsequent analysis within the framework
% of Section~\ref{subsec: framework}, we introduce a feature map
% $\phi_{\cR_\eta}:\bar\Omega\to L^2(\pi)$
% that incorporates the PDE and boundary operators.
% \begin{equation}
% \label{def residual_feature}
% \phi_{\cR_\eta}(\bm{x},\bm{\theta}) =  
% \begin{cases}
% \cL\phi(\bm{x},\bm{\theta}), & \bm{x} \in \Omega,\\
% \sqrt{\eta}\cB\phi(\bm{x},\bm{\theta}),  & \bm{x} \in \partial \Omega.
% \end{cases}
% \end{equation}
We express the approximant in operator form. Let $\cS_{m,\cR_\eta}:\RR^m\to\bar L^2(\bar\Omega)$
be the discrete operator, and let
$\cS_{m,\cR_\eta}^*:\bar L^2(\bar\Omega)\to\RR^m$
denote its adjoint.
For $\bm a\in\RR^m$ and $f\in\bar L^2(\bar\Omega)$,
these operators are given by
\begin{equation}
\label{def: discrete_residual}
\begin{aligned}
\cS_{m,\cR_\eta}\bm a
&=
\frac{1}{\sqrt m}\sum_{i=1}^m
\frac{a_i}{\sqrt{p(\bm\theta_i)}}
\phi_{\cR_\eta}(\cdot,\bm\theta_i),\\
(\cS_{m,\cR_\eta}^*f)_j
&=
\frac{
\langle\phi_{\cR_\eta}(\cdot,\bm\theta_j),f\rangle_
{\bar L^2(\bar\Omega)}
}{
\sqrt{mp(\bm\theta_j)}
},
\qquad j=1,\ldots,m.
\end{aligned}
\end{equation}
It follows that
\begin{equation}
\label{appro func operator form residual}
P_{\cR_\eta}\hat u^{\cR_\eta}_{m,\lambda}
=
\cS_{m,\cR_\eta}
\bigl(
\cS_{m,\cR_\eta}^*\cS_{m,\cR_\eta}
+\lambda\bm I_m
\bigr)^{-1}
\cS_{m,\cR_\eta}^*P_{\cR_\eta}u^*.
\end{equation}
The corresponding operator $\Phi_{\cR_\eta}: L^2(\pi) \to \bar L^2(\bar\Omega)$ and its adjoint 
$\Phi_{\cR_\eta}^*$ are obtained by
specializing \eqref{eq: abstract feature operator} and
\eqref{eq: abstract adjoint feature operator}, respectively,
to $B= \bar L^2(\bar\Omega)$ with feature function
$\phi_{\cR_\eta}$. Explicitly, for $\alpha\in L^2(\pi)$ and $f\in \bar L^2(\bar\Omega)$,  \begin{equation*}
\label{def: Phi_re and Phi_re*}
(\Phi_{\cR_\eta}\alpha)(\bm{x})
= \langle \phi_{\cR_\eta}(\bm{x},\cdot),\alpha\rangle_{L^2(\pi)}, \quad 
(\Phi_{\cR_\eta}^* f)(\bm{\theta})
= \langle \phi_{\cR_\eta}(\cdot,\bm{\theta}),f\rangle_{\bar L^2(\bar\Omega)}.
\end{equation*}
To connect this formulation with KRR on $\Theta$, we introduce the kernel on the parameter space $\Theta$:
\begin{equation*}
\label{new_k_residual}
\begin{aligned}
\tilde k_{\cR_\eta}(\bm\theta,\bm\theta')
&=
\left\langle
\phi_{\cR_\eta}(\cdot,\bm\theta),
\phi_{\cR_\eta}(\cdot,\bm\theta')
\right\rangle_{\bar L^2(\bar\Omega)}
\\
&=
\left\langle
\cL\phi(\cdot,\bm\theta),
\cL\phi(\cdot,\bm\theta')
\right\rangle_{L^2(\rho_I)}
+\eta
\left\langle
\cB\phi(\cdot,\bm\theta),
\cB\phi(\cdot,\bm\theta')
\right\rangle_{L^2(\rho_B)}.
\end{aligned}
\end{equation*}
Denote by $\tilde{\mathcal H}_{\cR_\eta}$ the RKHS
associated with $\tilde k_{\cR_\eta}$.
The corresponding integral operator
$\tilde{\mathcal T}_{\cR_\eta}$ is given by
\begin{equation}
\label{int_oper_appro_residual}
\tilde{\mathcal T}_{\cR_\eta}g
=
\int_\Theta
\tilde k_{\cR_\eta}(\cdot,\bm\theta')
g(\bm\theta')\,\dd\pi(\bm\theta')
=
\Phi_{\cR_\eta}^*\Phi_{\cR_\eta}g.
\end{equation}
Its spectral decomposition is
\begin{equation*}
\label{mercer's decomposition_dual_residual}
\tilde{\mathcal T}_{\cR_\eta}
=
\sum_{j=1}^\infty
\tilde\mu_j^{(\cR_\eta)}
\tilde e_j^{(\cR_\eta)}
\otimes
\tilde e_j^{(\cR_\eta)},
\end{equation*}
where $\{\tilde\mu_j^{(\cR_\eta)}\}_{j=1}^\infty$
are the positive eigenvalues in nonincreasing order and
$\{\tilde e_j^{(\cR_\eta)}\}_{j=1}^\infty$
are the corresponding orthonormal eigenfunctions in
$L^2(\pi)$. For the target function $\langle\phi_{\cR_{\eta}}(\cdot,\bm \theta),P_{\cR_{\eta}}u^*(\cdot)\rangle_{\bar L^2(\bar\Omega)}=(\Phi_{\cR_\eta}^*P_{\cR_\eta}u^*)(\bm \theta)$, the corresponding KRR estimator admits the operator
representation
\begin{equation}
\label{express g_residual}
\begin{aligned}
\hat g^{\cR_\eta}_{m,\lambda}
&=
\Phi_{\cR_\eta}^*\cS_{m,\cR_\eta}
\bigl(
\cS_{m,\cR_\eta}^*\cS_{m,\cR_\eta}
+\lambda\bm I_m
\bigr)^{-1}
\cS_{m,\cR_\eta}^*P_{\cR_\eta}u^*\\
&=
\Phi_{\cR_\eta}^*
P_{\cR_\eta}\hat u^{\cR_\eta}_{m,\lambda}.
\end{aligned}
\end{equation}
Let $\cH_{\cR_\eta}^s$ and
$\widetilde{\mathcal H}_{\cR_\eta}^s$ denote the
interpolation spaces defined in \eqref{eqn: Hks-space-x}
and \eqref{eqn: Hks-space-theta}, respectively, with
$B=\bar L^2(\bar\Omega)$ and feature function
$\phi_{\cR_\eta}$. 
The minimum-coefficient characterization identifies
$\cH_{\cR_\eta}^1$ with $P_{\cR_\eta}\mathcal H_k$,
equipped with the norm
\[
\|h\|_{\cH_{\cR_\eta}^1}
=
\min_{\substack{\alpha\in L^2(\pi)\\ \Phi_{\cR_{\eta}}\alpha=h}}
\|\alpha\|_{L^2(\pi)}
=
\min_{\substack{v\in\mathcal H_k\\P_{\cR_\eta}v=h}}
\|v\|_{\mathcal H_k}.
\]
Thus, $P_{\cR_\eta}$ maps the unit ball of $\mathcal H_k$
onto that of $\cH_{\cR_\eta}^1$. For $h^*\in\cH_{\cR_\eta}^1$, let $g^*(\bm\theta)
:=\Phi_{\cR_{\eta}}^*h^*$.
By Lemma~\ref{lemma RKHS equivalence}, the map
$h^*\mapsto g^*$ is an isometric isomorphism from
$\cH_{\cR_\eta}^1$ onto
$\widetilde{\mathcal H}_{\cR_\eta}^{2}$.
In particular,
\(
\|g^*\|_{\widetilde{\mathcal H}_{\cR_\eta}^{2}}
=
\|h^*\|_{\cH_{\cR_\eta}^1}.
\)
Lemma~\ref{lemma: approximate duality revised}
with $s=2$ and $p=1$ gives
\[
\begin{aligned}
\sup_{\|u^*\|_{\mathcal H_k}\leq1}
\inf_{u\in\sigmam}
\|u-u^*\|_{\cR_\eta}^2&=
\sup_{\|h^*\|_{\cH_{\cR_\eta}^1}\leq1}
\inf_{h\in P_{\cR_\eta}\sigmam}
\|h-h^*\|_{\bar L^2(\bar\Omega)}^2\\
&\quad=
\sup_{\|g^*\|_{\widetilde{\mathcal H}_{\cR_\eta}^{2}}\leq1}
\inf_{\bm a\in\mathbb R^m}
\left\|
\sum_{i=1}^m
a_i\tilde k_{\cR_\eta}(\cdot,\bm\theta_i)-g^*
\right\|_{\widetilde{\mathcal H}_{\cR_\eta}^{1}}^2.
\end{aligned}
\]
Consequently, the upper and lower bounds for the approximation error
in the residual seminorm follow by the same arguments
as in the proofs of Theorems~\ref{upper_appro} and
\ref{lower_appro}, respectively.

\subsection{Estimation Error Analysis}
\label{sec: estimation}
We now study the estimation error associated with
a finite number of collocation observations. For $\lambda>0$, consider the regularized estimator defined in \eqref{est_problem}. By the feature-space characterization
of $\cH_k$ in \eqref{eq: rkhs feature charact}, 
\begin{equation}
\label{est_problem feature}
\hat u_{n,\eta,\lambda}=\Phi\left(\argmin_{b\in L^2(\pi)}\left\{
\hat{\cR}_{n,\eta}(\Phi b)
+\lambda\|b\|_{L^2(\pi)}^2\right\}\right).
\end{equation}
We next reformulate \eqref{est_problem feature} as a KRR problem. Writing
$\nu_I=n/n_I$ and $\nu_B=n/n_B$, we have
\[
\hat{\cR}_{n,\eta}(u)
=
\frac1n
\left[
\sum_{i=1}^{n_I}
\frac{
|\sqrt{\nu_I}\cL u(\bm x_i)
-\sqrt{\nu_I}f^*(\bm x_i)|^2
}{q_I(\bm x_i)}
+
\sum_{i=n_I+1}^{n}
\frac{
|\sqrt{\eta\nu_B}\cB u(\bm x_i)
-\sqrt{\eta\nu_B}g^*(\bm x_i)|^2
}{q_B(\bm x_i)}
\right].
\]
Thus we define
$P_{\mathrm{PDE}}:D_{\mathcal R}\to\bar L^2(\bar\Omega)$ by
\begin{equation*}
\label{def pde_mapping}
(P_{\mathrm{PDE}}u)(\bm x)
:=
\begin{cases}
\sqrt{\nu_I}\cL u(\bm x),
&\bm x\in\Omega,\\
\sqrt{\eta\nu_B}\cB u(\bm x),
&\bm x\in\partial\Omega.
\end{cases}
\end{equation*}
The corresponding feature function is
\begin{equation*}
\label{def pde_feature}
\phi_{\mathrm{PDE}}(\bm x,\bm\theta)
:=
\bigl(P_{\mathrm{PDE}}\phi(\cdot,\bm\theta)\bigr)(\bm x)
=
\begin{cases}
\sqrt{\nu_I}\cL\phi(\bm x,\bm\theta),
&\bm x\in\Omega,\\
\sqrt{\eta\nu_B}\cB\phi(\bm x,\bm\theta),
&\bm x\in\partial\Omega.
\end{cases}
\end{equation*}
The transformed features therefore induce the kernel
\[
\begin{aligned}
k_{\mathrm{PDE}}(\bm x,\bm x')
&=
\int_\Theta
\phi_{\mathrm{PDE}}(\bm x,\bm\theta)
\phi_{\mathrm{PDE}}(\bm x',\bm\theta)
\,\dd\pi(\bm\theta)=
\langle
\phi_{\mathrm{PDE}}(\bm x,\cdot),
\phi_{\mathrm{PDE}}(\bm x',\cdot)
\rangle_{L^2(\pi)}.\\
\end{aligned}
\]
Let $\mathcal H_{\mathrm{PDE}}$ denote the associated RKHS,
viewed as continuously embedded in $\bar L^2(\bar\Omega)$.
The feature operator
$\Phi_{\mathrm{PDE}}:L^2(\pi)\to\bar L^2(\bar\Omega)$
associated with $\phi_{\mathrm{PDE}}$ and its adjoint
$\Phi_{\mathrm{PDE}}^*:\bar L^2(\bar\Omega)\to L^2(\pi)$
are given by
\begin{equation*}
\label{def: Phi_PDE}
(\Phi_{\mathrm{PDE}}b)(\bm x)
=
\langle\phi_{\mathrm{PDE}}(\bm x,\cdot),b\rangle_{L^2(\pi)},
\qquad
(\Phi_{\mathrm{PDE}}^*f)(\bm\theta)
=
\langle\phi_{\mathrm{PDE}}(\cdot,\bm\theta),f\rangle_{\bar L^2(\bar\Omega)}.
\end{equation*}
Since Assumption~\ref{assump:feature-operator-compatibility} implies $\Phi_{\mathrm{PDE}}=P_{\mathrm{PDE}}\Phi$,
the minimum-coefficient characterization identifies
$\mathcal H_{\mathrm{PDE}}$ with
$P_{\mathrm{PDE}}\mathcal H_k$, equipped with the norm
\[
\|h\|_{\mathcal H_{\mathrm{PDE}}}
=
\min_{\substack{b\in L^2(\pi)\\
                 \Phi_{\mathrm{PDE}}b=h}}
\|b\|_{L^2(\pi)}
=
\min_{\substack{u\in\mathcal H_k\\
                 P_{\mathrm{PDE}}u=h}}
\|u\|_{\mathcal H_k}.
\]
Thus, $P_{\mathrm{PDE}}$ maps the closed unit ball of
$\mathcal H_k$ onto that of $\mathcal H_{\mathrm{PDE}}$. Let $h^*=P_{\mathrm{PDE}}u^*$. Then \eqref{est_problem feature} is equivalent
to a KRR problem in $\mathcal H_{\mathrm{PDE}}$, with
\[
P_{\mathrm{PDE}}\hat u_{n,\eta,\lambda}
=
\argmin_{h\in\mathcal H_{\mathrm{PDE}}}
\left\{
\frac1n\sum_{i=1}^n
\frac{|h(\bm x_i)-h^*(\bm x_i)|^2}
     {q(\bm x_i)}
+\lambda\|h\|_{\mathcal H_{\mathrm{PDE}}}^2
\right\},
\]
where $q(\bm x):=q_I(\bm x)$ for $\bm x\in\Omega$ and
$q(\bm x):=q_B(\bm x)$ for $\bm x\in\partial\Omega$.
To derive the estimator explicitly, set
\[
(\bm\phi_{n,\mathrm{PDE}}(\bm\theta))_i
:=
\frac{\phi_{\mathrm{PDE}}(\bm x_i,\bm\theta)}
     {\sqrt{nq(\bm x_i)}},
\qquad
(\bm y_{\mathrm{PDE}})_i
:=
\frac{h^*(\bm x_i)}
     {\sqrt{nq(\bm x_i)}},
\quad i=1,\ldots,n.
\]
The corresponding discrete observation operator
$\cS_{n,\mathrm{PDE}}:L^2(\pi)\to\mathbb R^n$
and its adjoint are given by
\[
(\cS_{n,\mathrm{PDE}}b)_i
=
\frac{
\langle\phi_{\mathrm{PDE}}(\bm x_i,\cdot),b\rangle_{L^2(\pi)}
}{\sqrt{nq(\bm x_i)}},
\quad i=1,\ldots,n,
\qquad
\cS_{n,\mathrm{PDE}}^*\bm a
=
\bm\phi_{n,\mathrm{PDE}}(\cdot)^\top\bm a,
\]
for $b\in L^2(\pi)$ and $\bm a\in\mathbb R^n$. The coefficient minimization in \eqref{est_problem feature}
is thus finding
\[
\hat b_{n,\eta,\lambda}=\argmin_{b\in L^2(\pi)}
\left\{
\|\cS_{n,\mathrm{PDE}}b-\bm y_{\mathrm{PDE}}\|_{\mathbb R^n}^2
+\lambda\|b\|_{L^2(\pi)}^2
\right\}.
\]
For every $v\in L^2(\pi)$, the first-order optimality
condition gives
\[
\left\langle
(\cS_{n,\mathrm{PDE}}\hat b_{n,\eta,\lambda}
-\bm y_{\mathrm{PDE}}),
\cS_{n,\mathrm{PDE}}v
\right\rangle_{\mathbb R^n}
+
\lambda
\langle\hat b_{n,\eta,\lambda},v\rangle_{L^2(\pi)}
=0.
\]
Since $v$ is arbitrary, this is equivalent to
\[
\cS_{n,\mathrm{PDE}}^*
\bigl(
\cS_{n,\mathrm{PDE}}\hat b_{n,\eta,\lambda}
-\bm y_{\mathrm{PDE}}
\bigr)
+\lambda\hat b_{n,\eta,\lambda}=0.
\]
Hence
\[
\hat b_{n,\eta,\lambda}
=
(\cS_{n,\mathrm{PDE}}^*\cS_{n,\mathrm{PDE}}+\lambda I)^{-1}
\cS_{n,\mathrm{PDE}}^*\bm y_{\mathrm{PDE}}=\cS_{n,\mathrm{PDE}}^*
( \cS_{n,\mathrm{PDE}}\cS_{n,\mathrm{PDE}}^*+\lambda\bm I_n)^{-1}
\bm y_{\mathrm{PDE}},
\]
where the second equality follows from \eqref{matrix_inverse_lemma}. Since $\bm y_{\mathrm{PDE}}
=\cS_{n,\mathrm{PDE}}b^*$, we obtain
\begin{equation}
\label{best estimation estimator}
\hat u_{n,\eta,\lambda}
=
\Phi\hat b_{n,\eta,\lambda}
=
\Phi\left(
\cS_{n,\mathrm{PDE}}^*
( \cS_{n,\mathrm{PDE}}\cS_{n,\mathrm{PDE}}^*+\lambda\bm I_n)^{-1}
\cS_{n,\mathrm{PDE}}b^*
\right).
\end{equation}

For any $u=\Phi b$, the estimation error in terms of the residual seminorm is determined by:
\[
\begin{aligned}
\|u-u^*\|_{\cR_\eta}^2
&=
\|P_{\cR_\eta}(u-u^*)\|_{\bar L^2(\bar\Omega)}^2=
\|\Phi_{\cR_\eta}(b-b^*)\|_{\bar L^2(\bar\Omega)}^2\\
&=
\left\langle
\widetilde{\cT}_{\cR_\eta}(b-b^*),b-b^*
\right\rangle_{L^2(\pi)}=
\left\|
\widetilde{\cT}_{\cR_\eta}^{1/2}(b-b^*)
\right\|_{L^2(\pi)}^2.\\
\end{aligned}
\]
Here
$\widetilde{\cT}_{\cR_\eta}
=\Phi_{\cR_\eta}^*\Phi_{\cR_\eta}$,
and $\widetilde{\cT}_{\cR_\eta}^{1/2}$ denotes its
positive square root. We therefore relate this operator to
$\widetilde{\cT}_{\mathrm{PDE}}:=\Phi_{\mathrm{PDE}}^*\Phi_{\mathrm{PDE}}$. To express $\widetilde{\cT}_{\mathrm{PDE}}$ and $\widetilde{\cT}_{\cR_\eta}$ using the same
interior and boundary components, we introduce
$\Phi_{\cL}:L^2(\pi)\to L^2(\rho_I)$ and
$\Phi_{\cB}:L^2(\pi)\to L^2(\rho_B)$.
These operators and their adjoints are given by
\begin{equation}
\label{def Phi L operator}
\Phi_{\cL}b
=
\int_\Theta
\cL\phi(\cdot,\bm\theta)b(\bm\theta)
\,\dd\pi(\bm\theta),
\qquad
\Phi_{\cL}^*f
=
\int_\Omega
\cL\phi(\bm x,\cdot)f(\bm x)
\,\dd\rho_I(\bm x),
\end{equation}
and
\begin{equation}
\label{def Phi B operator}
\Phi_{\cB}b
=
\int_\Theta
\cB\phi(\cdot,\bm\theta)b(\bm\theta)
\,\dd\pi(\bm\theta),
\qquad
\Phi_{\cB}^*f
=
\int_{\partial\Omega}
\cB\phi(\bm x,\cdot)f(\bm x)
\,\dd\rho_B(\bm x).
\end{equation}
Writing
$\widetilde{\cT}_{\cL}:=\Phi_{\cL}^*\Phi_{\cL}$
and
$\widetilde{\cT}_{\cB}:=\Phi_{\cB}^*\Phi_{\cB}$,
we obtain
\(
\widetilde{\cT}_{\cR_\eta}
=
\widetilde{\cT}_{\cL}
+\eta\widetilde{\cT}_{\cB},
\)
and
\begin{equation}
\label{pde_int_operator}
\widetilde{\cT}_{\mathrm{PDE}}
:=
\Phi_{\mathrm{PDE}}^*\Phi_{\mathrm{PDE}}
=
\nu_I\widetilde{\cT}_{\cL}
+\eta\nu_B\widetilde{\cT}_{\cB}.
\end{equation}
Thus, setting
$\nu_-:=\min\{\nu_I,\nu_B\}$ and
$\nu_+:=\max\{\nu_I,\nu_B\}$, we have
\begin{equation}
\label{mu_PDE mu_R_eta} 
\nu_-\widetilde{\cT}_{\cR_\eta}
\preccurlyeq
\widetilde{\cT}_{\mathrm{PDE}}
\preccurlyeq
\nu_+\widetilde{\cT}_{\cR_\eta}.
\end{equation}

The corresponding comparison in the function space is
\begin{equation}
\label{eq: pde-residual-comparison}
\nu_-\|v\|_{\cR_\eta}^2
\leq
\|P_{\mathrm{PDE}}v\|_{\bar L^2(\bar\Omega)}^2
\leq
\nu_+\|v\|_{\cR_\eta}^2,
\qquad v\in D_{\mathcal R}.
\end{equation}

To prove Theorem~\ref{thm: Estimation_error upper},
we introduce the relevant degrees of freedom.
Let $(\mu_j^{\cL},e_j^{\cL})$ and
$(\mu_j^{\cB},e_j^{\cB})$ denote the positive eigenpairs
of the integral operators
$\Phi_{\cL}\Phi_{\cL}^*:L^2(\rho_I)\to L^2(\rho_I)$
and
$\Phi_{\cB}\Phi_{\cB}^*:L^2(\rho_B)\to L^2(\rho_B)$,
respectively.
Their kernels $k_{\cL}$ and $k_{\cB}$ are defined in
\eqref{eq: k_L} and \eqref{eq: k_B}. Following \eqref{eq: dof}, for $\lambda>0$, define
\[
N_{\cL}(\lambda)
:=
\sum_{j\geq1}
\frac{\mu_j^{\cL}}{\mu_j^{\cL}+\lambda},
\qquad
N_{\cB}(\lambda)
:=
\sum_{j\geq1}
\frac{\mu_j^{\cB}}{\mu_j^{\cB}+\lambda},
\qquad 
N_{\mathrm{PDE}}(\lambda)
:=
\sum_{j\geq1}
\frac{\mu_j^{\mathrm{PDE}}}
     {\mu_j^{\mathrm{PDE}}+\lambda}.
\]
For the sampling densities $q_I$ and $q_B$, define
\[
F_{\cL}(\lambda)
:=
\operatorname*{ess\,sup}_{\bm x\in\Omega}
\frac{1}{q_I(\bm x)}
\sum_{j\geq1}
\frac{\mu_j^{\cL}}{\mu_j^{\cL}+\lambda}
|e_j^{\cL}(\bm x)|^2,
\]
and
\[
F_{\cB}(\lambda)
:=
\operatorname*{ess\,sup}_{\bm x\in\partial\Omega}
\frac{1}{q_B(\bm x)}
\sum_{j\geq1}
\frac{\mu_j^{\cB}}{\mu_j^{\cB}+\lambda}
|e_j^{\cB}(\bm x)|^2.
\]
The dependence of $F_{\cL}$ and $F_{\cB}$ on
$q_I$ and $q_B$ is suppressed in the notation. To characterize the admissible range of $\lambda$ in Theorem~\ref{thm: Estimation_error upper}, we introduce a concept of critical penalty.
\begin{definition}[Critical penalty]
For $\lambda>0$, set
\[
F_I(\lambda)
:=
\max\{N_{\mathrm{PDE}}(\lambda),F_{\cL}(\lambda)\},
\]
and
\[
F_B(\lambda)
:=
\max\{N_{\mathrm{PDE}}(\lambda),F_{\cB}(\lambda/\eta)\}.
\]
For $n\in\mathbb N^+$ and $\delta\in(0,1)$, define
\begin{equation}
\label{critical penalty L}
\Lambda_I(n,\delta)
:=
\inf\left\{
\lambda>0:
n\geq16F_I(\lambda)
\max\left\{1,\log\frac{14F_I(\lambda)}{\delta}\right\}
\right\},
\end{equation}
and
\begin{equation}
\label{critical penalty B}
\Lambda_B(n,\delta)
:=
\inf\left\{
\lambda>0:
n\geq16F_B(\lambda)
\max\left\{1,\log\frac{14F_B(\lambda)}{\delta}\right\}
\right\}.
\end{equation}
\end{definition}
We define the weighted empirical operators as 
\[
\widetilde{\cT}_{\cL,n_I}
=
\frac1{n_I}\sum_{i=1}^{n_I}
\frac{\cL\phi(\bm x_i,\cdot)\otimes\cL\phi(\bm x_i,\cdot)}
{q_I(\bm x_i)},
\quad
\widetilde{\cT}_{\cB,n_B}
=
\frac1{n_B}\sum_{i=1}^{n_B}
\frac{\cB\phi(\bm x_{n_I+i},\cdot)\otimes\cB\phi(\bm x_{n_I+i},\cdot)}
{q_B(\bm x_{n_I+i})}.
\]
Thus,
\(
A_n
:=
\cS_{n,\mathrm{PDE}}^*\cS_{n,\mathrm{PDE}}
=
\widetilde{\cT}_{\cL,n_I}
+\eta\widetilde{\cT}_{\cB,n_B},
\)
and
\(
\mathbb E[A_n]
=
\widetilde{\cT}_{\cL}+\eta\widetilde{\cT}_{\cB}
=
\widetilde{\cT}_{\cR_\eta}.
\)

We now prove the estimation
upper bound in Theorem~\ref{thm: Estimation_error upper}.
\begin{proof}[Proof of Theorem \ref{thm: Estimation_error upper}]
Let $b^*$ be the minimum-norm coefficient representing
$u^*=\Phi b^*$. Then
\[
\|b^*\|_{L^2(\pi)}=\|u^*\|_{\mathcal H_k}.
\]
By \eqref{best estimation estimator} and \eqref{matrix_inverse_lemma},
\[
\begin{aligned}
\hat b_{n,\eta,\lambda}-b^*
&=
\left[
\cS_{n,\mathrm{PDE}}^*
(\cS_{n,\mathrm{PDE}}\cS_{n,\mathrm{PDE}}^*
+\lambda\bm I_n)^{-1}
\cS_{n,\mathrm{PDE}}
-I
\right]b^*\\
&=
\left[
(A_n+\lambda I)^{-1}A_n-I
\right]b^*\\
&=
-\lambda(A_n+\lambda I)^{-1}b^*.
\end{aligned}
\]
It follows that
\[
\begin{aligned}
\|\hat u_{n,\eta,\lambda}-u^*\|_{\cR_\eta}
&=
\left\|
\widetilde{\cT}_{\cR_\eta}^{1/2}
(\hat b_{n,\eta,\lambda}-b^*)
\right\|_{L^2(\pi)}\\
&\leq
\lambda
\left\|
\widetilde{\cT}_{\cR_\eta}^{1/2}
(A_n+\lambda I)^{-1}
\right\|
\|b^*\|_{L^2(\pi)}.
\end{aligned}
\]
Fix a deterministic $\lambda_n>0$ satisfying \eqref{lambda union condition}, specifically
\begin{equation*}
\lambda_n\geq
\max\{\Lambda_I(n_I,\delta/2),\Lambda_B(n_B,\delta/2)\}.  
\end{equation*}
By Lemmas~\ref{lemma:concentration}
and~\ref{lemma:constant bound}, with probability at least
$1-\delta$,
\[
\begin{aligned}
\left\|
\widetilde{\cT}_{\cR_\eta}^{1/2}
(A_n+\lambda_n I)^{-1/2}
\right\|&\leq
\left\|
\widetilde{\cT}_{\cR_\eta}^{1/2}
(\widetilde{\cT}_{\cR_\eta}+\lambda_n I)^{-1/2}
\right\|
\left\|
(\widetilde{\cT}_{\cR_\eta}+\lambda_n I)^{1/2}
(A_n+\lambda_n I)^{-1/2}
\right\|\\
&\leq2.
\end{aligned}
\]
On this event, for every $0<\lambda\leq\lambda_n$,
Lemma~\ref{lemma: AB-bound} gives
\[
\begin{aligned}
\lambda
\left\|
\widetilde{\cT}_{\cR_\eta}^{1/2}
(A_n+\lambda I)^{-1}
\right\|&\leq
\lambda
\left\|
\widetilde{\cT}_{\cR_\eta}^{1/2}
(A_n+\lambda_n I)^{-1/2}
\right\|\\
&\qquad\quad\times
\left\|
(A_n+\lambda_n I)^{1/2}
(A_n+\lambda I)^{-1/2}
\right\|
\left\|
(A_n+\lambda I)^{-1/2}
\right\|\\
&\leq
2\lambda
\left(\frac{\lambda_n}{\lambda}\right)^{1/2}
\lambda^{-1/2}
=
2\lambda_n^{1/2}.
\end{aligned}
\]
Consequently,
\[
\sup_{\|u^*\|_{\mathcal H_k}\leq1}
\cR_\eta(\hat u_{n,\eta,\lambda})
\leq4\lambda_n,
\qquad 0<\lambda\leq\lambda_n.
\]

We now choose the sampling densities so that $\lambda_n$
can be selected from the PDE eigenvalues.
By Lemma~\ref{lemma: N(gamma)}, there exists a constant
$C>0$ such that
\[
N_{\mathrm{PDE}}(\mu_j^{\mathrm{PDE}})
\leq Cj,
\qquad j\geq 1.
\]
Since the sampling proportions $1/\nu_I$ and $1/\nu_B$
are fixed and positive, choose $c_0>0$ such that
\[
Cc_0
\leq
\frac1{40}\min\{1/\nu_I,1/\nu_B\},
\]
and then choose $n_0>0$ sufficiently large that
$\lfloor c_0n_0\rfloor\geq 1$.
For every $t\geq n_0$, we therefore have
\[
\begin{aligned}
N_{\mathrm{PDE}}\!
\left(\mu_{\lfloor c_0t\rfloor}^{\mathrm{PDE}}\right)\leq C\lfloor c_0t\rfloor\leq Cc_0t\leq
\frac{t}{40}\min\{1/\nu_I,1/\nu_B\}.
\end{aligned}
\]
The constants $c_0$ and $n_0$ are independent of $n$
and $\delta$.
For
\(
\frac{n}{\log(28n/\delta)}\geq n_0,
\)
set
\[
l_n
:=
\left\lfloor
\frac{c_0n}{\log(28n/\delta)}
\right\rfloor,
\qquad
\lambda_n:=\mu_{l_n}^{\mathrm{PDE}}>0.
\]
Choose
\[
q_I(\bm x)
=q^*(k_\cL,\rho_I,\lambda_n;\bm x)=
\frac12+
\frac{1}{2N_{\cL}(\lambda_n)}
\sum_{j\geq1}
\frac{\mu_j^{\cL}}{\mu_j^{\cL}+\lambda_n}
|e_j^{\cL}(\bm x)|^2,
\]
and
\[
q_B(\bm x)
=q^*(k_\cB,\rho_B,\lambda_n/\eta;\bm x)=
\frac12+
\frac{1}{2N_{\cB}(\lambda_n/\eta)}
\sum_{j\geq1}
\frac{\mu_j^{\cB}}{\mu_j^{\cB}+\lambda_n/\eta}
|e_j^{\cB}(\bm x)|^2.
\]
By construction,
\[
F_{\cL}(\lambda_n)\leq2N_{\cL}(\lambda_n),
\qquad
F_{\cB}(\lambda_n/\eta)
\leq2N_{\cB}(\lambda_n/\eta).
\]
Since $\nu_I,\nu_B\geq1$, the decomposition of
$\widetilde{\cT}_{\mathrm{PDE}}$ gives
\[
\widetilde{\cT}_{\cL}
\preccurlyeq\widetilde{\cT}_{\mathrm{PDE}},
\qquad
\eta\widetilde{\cT}_{\cB}
\preccurlyeq\widetilde{\cT}_{\mathrm{PDE}}.
\]
Then
\[
\begin{aligned}
\widetilde{\cT}_{\cL}
(\widetilde{\cT}_{\cL}+\lambda_n I)^{-1}
=
I-\lambda_n
(\widetilde{\cT}_{\cL}+\lambda_n I)^{-1}&\preccurlyeq
I-\lambda_n
(\widetilde{\cT}_{\mathrm{PDE}}+\lambda_n I)^{-1}\\
&=
\widetilde{\cT}_{\mathrm{PDE}}
(\widetilde{\cT}_{\mathrm{PDE}}+\lambda_n I)^{-1}.
\end{aligned}
\]
Taking traces yields
\[
N_{\cL}(\lambda_n)
=
\operatorname{tr}\left[
\widetilde{\cT}_{\cL}
(\widetilde{\cT}_{\cL}+\lambda_n I)^{-1}
\right]
\leq N_{\mathrm{PDE}}(\lambda_n).
\]
Applying the same argument to
$\eta\widetilde{\cT}_{\cB}$ gives
\[
\eta\widetilde{\cT}_{\cB}
(\eta\widetilde{\cT}_{\cB}+\lambda_n I)^{-1}
\preccurlyeq
\widetilde{\cT}_{\mathrm{PDE}}
(\widetilde{\cT}_{\mathrm{PDE}}+\lambda_n I)^{-1}.
\]
Consequently,
\[
\begin{aligned}
N_{\cB}(\lambda_n/\eta)
=
\operatorname{tr}\left[
\widetilde{\cT}_{\cB}
\bigl(\widetilde{\cT}_{\cB}+(\lambda_n/\eta)I\bigr)^{-1}
\right]&=
\operatorname{tr}\left[
\eta\widetilde{\cT}_{\cB}
(\eta\widetilde{\cT}_{\cB}+\lambda_n I)^{-1}
\right]\\
&\leq N_{\mathrm{PDE}}(\lambda_n).
\end{aligned}
\]
Therefore,
\[
F_I(\lambda_n),\,F_B(\lambda_n)
\leq2N_{\mathrm{PDE}}(\lambda_n)
\leq
\frac{\min\{n_I,n_B\}}{20\log(28n/\delta)}.
\]
For $J\in\{I,B\}$, it follows that
\[
16F_J(\lambda_n)
\max\left\{1,\log\frac{28F_J(\lambda_n)}{\delta}\right\}
\leq
\frac45\min\{n_I,n_B\}
\leq n_J.
\]
Thus the chosen $\lambda_n$ satisfies \eqref{lambda union condition}. Therefore with
probability at least $1-\delta$,
\[
\sup_{\|u^*\|_{\mathcal H_k}\leq1}
\|\hat u_{n,\eta,\lambda}-u^*\|^2_{\cR_\eta}
\leq
4\mu_{\left\lfloor c_0n/\log(28n/\delta)\right\rfloor}
^{\mathrm{PDE}},
\qquad 0<\lambda\leq\lambda_n.
\]
\end{proof}
We next establish the lower bound on the worst-case
estimation error stated in
Theorem~\ref{thm: Estimation_error lower}.
\begin{proof}[Proof of Theorem \ref{thm: Estimation_error lower}]
Set $h^*:=P_{\mathrm{PDE}}u^*$. For $h\in\mathcal H_{\mathrm{PDE}}$, write
\[
h(X_n):=
(h(\bm x_1),\ldots,h(\bm x_n))^\top,
\qquad
\bm D_\nu:=
\operatorname{diag}
\bigl(
\sqrt{\nu_I}\bm I_{n_I},
\sqrt{\eta\nu_B}\bm I_{n_B}
\bigr).
\]
Then
\[
h^*(X_n)=\bm D_\nu\bm z_n(X_n,u^*).
\]
Since $\bm D_\nu$ is invertible, the PDE observations and
the pointwise observations of $h^*$ determine each other. Each $G\in\mathcal M_n^{\mathrm{PDE}}$ therefore induces
an estimator
\[
\mathcal Q_Gh
:=
P_{\mathrm{PDE}}
G\bigl(X_n,\bm D_\nu^{-1}h(X_n)\bigr).
\]
By \eqref{eq: pde-residual-comparison},
\[
\|G(S_n(u^*))-u^*\|_{\cR_\eta}^2
\geq
\frac1{\nu_+}
\|\mathcal Q_Gh^*-h^*\|_{\bar L^2(\bar\Omega)}^2.
\]
Let $\mathcal M_n$ denote the class of all estimators based
only on $(X_n,h(X_n))$, with outputs in
$\bar L^2(\bar\Omega)$.
Since
$\{\mathcal Q_G:G\in\mathcal M_n^{\mathrm{PDE}}\}
\subseteq\mathcal M_n$, we obtain
\[
\begin{aligned}
\inf_{G\in\mathcal M_n^{\mathrm{PDE}}}
\sup_{\|u^*\|_{\mathcal H_k}\leq1}
\|G(S_n(u^*))-u^*\|_{\cR_\eta}^2&\geq
\frac1{\nu_+}
\inf_{\mathcal Q\in\mathcal M_n}
\sup_{\|h^*\|_{\mathcal H_{\mathrm{PDE}}}\leq1}
\|\mathcal Qh^*-h^*\|_{\bar L^2(\bar\Omega)}^2\\
&\geq
\frac{\mu_{n+1}^{\mathrm{PDE}}}{\nu_+}.
\end{aligned}
\]
The last inequality follows from
\cite[Theorem~4.1]{long2025optimalratessaturationnoiseless}
with $s=1$ and $p=0$.
\end{proof}
We now consider the simultaneous restrictions imposed by
finitely many features and finitely many collocation observations.
Fix feature parameters
$\bm\theta_1,\ldots,\bm\theta_m\in\Theta$ and collocation points
\[
X_n=(\bm x_1,\ldots,\bm x_n)
\in\Omega^{n_I}\times(\partial\Omega)^{n_B},
\qquad n=n_I+n_B.
\]
Let $\mathcal M_n^{\mathrm{PDE}}$ be the class of estimators
based only on collocation observations, as defined in \eqref{M_n pde}. Let $\mathcal M_{m,n}$ denote the class of estimators
defined in \eqref{M_mn}, which use collocation
observations and take values in $\sigmam$.

To distinguish the restriction imposed by finitely many features from that imposed
by finitely many collocation observations, define
\begin{align*}
E_1
&:=
\sup_{\|u^*\|_{\mathcal H_k}\leq1}
\inf_{v\in\sigmam}
\|v-u^*\|_{\cR_\eta}^2,\\
E_2
&:=
\inf_{G\in\mathcal M_n^{\mathrm{PDE}}}
\sup_{\|u^*\|_{\mathcal H_k}\leq1}
\|G(S_n(u^*))-u^*\|_{\cR_\eta}^2.
\end{align*}
Here, $E_1$ measures the best approximation error in the
random feature model when the target is fully known,
whereas $E_2$ measures the optimal recovery error from
the collocation observations without restricting the
output to $\sigmam$.

\begin{proof}[Proof of Corollary \ref{thm: joint-residual-lower}]
For every $G\in\mathcal M_{m,n}$ and every
$u^*\in\mathcal H_k$, the inclusion
$G(S_n(u^*))\in\sigmam$ implies
\[
\|G(S_n(u^*))-u^*\|_{\cR_\eta}^2
\geq
\inf_{v\in\sigmam}
\|v-u^*\|_{\cR_\eta}^2.
\]
Taking the supremum over the unit ball of $\mathcal H_k$
and then the infimum over $G\in\mathcal M_{m,n}$ gives
$ E_{m,n}\geq E_1$.
Moreover,
$\mathcal M_{m,n}\subseteq\mathcal M_n^{\mathrm{PDE}}$,
so $ E_{m,n}\geq E_2$.

It remains to relate these two errors to the PDE spectrum.
Since
\(
\widetilde{\mathcal T}_{\mathrm{PDE}}
\preccurlyeq
\nu_+\widetilde{\mathcal T}_{\cR_\eta},
\)
\[
\mu_j^{\mathrm{PDE}}
\leq
\nu_+\tilde\mu_j^{(\cR_\eta)},
\qquad j\geq1.
\]
Theorem~\ref{lower_appro_residual} therefore gives
\[
E_1
\geq
\tilde\mu_{m+1}^{(\cR_\eta)}
\geq
\frac{\mu_{m+1}^{\mathrm{PDE}}}{\nu_+}.
\]
On the other hand,
Theorem~\ref{thm: Estimation_error lower} gives
\[
E_2
\geq
\frac{\mu_{n+1}^{\mathrm{PDE}}}{\nu_+}.
\]
Combining these inequalities proves
\eqref{eq: joint-residual-lower}.
\end{proof}

\begin{remark}[A two-stage approximator]
\label{rem: two-stage-approximator}
Fix $\delta\in(0,1)$.
Choose $l_n,q_I,q_B$ according to
Theorem~\ref{thm: Estimation_error upper}
and $l_m,p_m$ according to
Theorem~\ref{upper_appro_residual},
in both cases replacing the failure probability by
$\delta/2$.
Assume that the hypotheses and sample-size requirements
of both theorems are satisfied.
First, sample the interior and boundary collocation points
independently from $q_I\,\dd\rho_I$ and $q_B\,\dd\rho_B$,
respectively, with the two collections independent.
Compute $\hat u_{n,\eta,\lambda}$ according to \eqref{est_problem} with some fixed $0<\lambda\leq\mu_{l_n}^{\mathrm{PDE}}$.
Next, independently of the collocation points, draw
$\bm\theta_1,\ldots,\bm\theta_m$ independently from
$p_m\,\dd\pi$. Taking $\hat u_{n,\eta,\lambda}\in \cH_k$ as a target function, define
\begin{equation}
\label{eq: constructed solution}
\hat u_n^m=\argmin_{u\in\sigmam} \|u-\hat u_{n,\eta,\lambda}\|_{\cR_\eta}.
\end{equation}
Hence,
\[
\|\hat u_n^m-u^*\|^2_{\cR_\eta}
\leq
2\|\hat u_n^m-\hat u_{n,\eta,\lambda}\|^2_{\cR_\eta}
+2\|\hat u_{n,\eta,\lambda}-u^*\|^2_{\cR_\eta}.
\]
Therefore, Theorems~\ref{thm: Estimation_error upper}
and~\ref{upper_appro_residual} imply that, with probability at least $1-\delta$,
\[
\sup_{\|u^*\|_{\mathcal H_k}\leq1}
\|\hat u_n^m-u^*\|^2_{\cR_\eta}
\leq
8\mu_{l_n}^{\mathrm{PDE}}
+
32\tilde\mu_{l_m}^{(\cR_\eta)}.
\]
This construction differs from the empirical least-squares
estimator $\hat u_{n,m,\eta}$ defined in \eqref{eq:RFM solution},
which directly minimizes the weighted empirical residual
in \eqref{discre_loss} over $\sigmam$.
We leave it to future work to determine whether this
estimator attains the same bound under suitable
sampling assumptions.
\end{remark}

\section{Conclusions}
\label{sec:conclusions}
% \section{Numerical Experiments}
% \subsection{Dot-product on the Sphere}
% \subsection{Simple Kernel on [0,1]}

% activation function: 
% \begin{itemize}
%     \item $\text{ReLU}^k$
%     \item $\sin$ or $\cos$
%     \item $\tanh$
% \end{itemize}
In this work, we analyze RFM within
a kernel framework, separating finite-feature approximation
from recovery over the full RKHS using finitely many
collocation observations. Under the stated assumptions, we establish
eigenvalue-based lower bounds and high-probability
upper bounds for worst-case errors over the unit ball
of $\mathcal H_k$. The approximation bounds cover Sobolev norms and the
PDE residual seminorm, while the estimation bounds
are measured in the residual seminorm. 
Combining the approximation and estimation lower bounds,
we obtain a lower bound for reconstruction with
$m$ features and $n$ collocation observations. 
This bound identifies the separate limitations imposed
by the feature space and the available observations. We further investigate how activation regularity
affects the relevant eigenvalues, obtaining algebraic
eigenvalue upper bounds for finitely smooth activations, and
stretched-exponential bounds for analytic activations. The upper bounds in the residual seminorm can also
be combined through the two-stage construction
in \eqref{eq: constructed solution}, which first
performs recovery over the full RKHS and then projects
onto the random feature model using exact residual
inner products.
Establishing a comparable bound for the direct empirical
least-squares RFM estimator in \eqref{eq:RFM solution}
remains an open problem.

\bmhead{Acknowledgements}

% Acknowledgements are not compulsory. Where included they should be brief. Grant or contribution numbers may be acknowledged.
This work is partially supported by the NSFC Major Research Plan - Interpretable and General Purpose Next-generation Artificial Intelligence (Nos. 92570001 and 92370205), NSFC grant 12425113, the Basic Research Program of Jiangsu under Grant BK20253039 and Grant BK20250480, and Key Laboratory of the Ministry of Education for Mathematical Foundations and Applications of Digital Technology, University of Science and Technology of China.

\bmhead{Data Availability} Data will be made available on request.

\section*{Declarations}
\bmhead{Competing interests} The authors declare that they have no known competing financial interests or personal relationships that could have appeared to influence the work reported in this paper.

% Please refer to Journal-level guidance for any specific requirements.

% \section*{Declarations}

% Some journals require declarations to be submitted in a standardised format. Please check the Instructions for Authors of the journal to which you are submitting to see if you need to complete this section. If yes, your manuscript must contain the following sections under the heading `Declarations':

% \begin{itemize}
% \item Conflict of interest/Competing interests (check journal-specific guidelines for which heading to use)
% \item Ethics approval and consent to participate
% \item Consent for publication
% \item Materials availability
% \end{itemize}

% \noindent
% If any of the sections are not relevant to your manuscript, please include the heading and write `Not applicable' for that section. 

% %%===================================================%%
% %% For presentation purpose, we have included        %%
% %% \bigskip command. Please ignore this.             %%
% %%===================================================%%
% \bigskip

% Springer journals and proceedings: \url{https://www.springer.com/gp/editorial-policies}

% Nature Portfolio journals: \url{https://www.nature.com/nature-research/editorial-policies}

% \textit{Scientific Reports}: \url{https://www.nature.com/srep/journal-policies/editorial-policies}

% BMC journals: \url{https://www.biomedcentral.com/getpublished/editorial-policies}
% \end{flushleft}

\begin{appendices}
\numberwithin{equation}{section}

\section{Auxiliary Results}
\label{math_preliminaries}
\subsection{Singular value decomposition (SVD)}
\label{sec: bounded-operator}
% \subsection{Definitions}
% {\bf Bounded linear operator.~}
For a compact operator $A:\cH_1\mapsto\cH_2$, both $AA^*$ and $A^*A$ are compact, self-adjoint, and positive operators. 
Moreover, $AA^*$ and $A^*A$ share the same non-zero eigenvalues, and the singular value decomposition of $A$ is given by
\begin{equation}
\label{eqn: svd origin}
    A  = \sum_{j=1}^{\infty} \mu_j^{\frac{1}{2}} \psi'_j\otimes \psi_j,
\end{equation}
where $\{\mu_j\}_{j=1}^\infty$ are the eigenvalues of $AA^*$ in a decreasing order, $\{\psi_j\}_{j=1}^\infty$ and $\{\psi'_j\}_{j=1}^\infty$ are the corresponding orthonormal eigenfunctions of $A^*A$ and $AA^*$, respectively.  Through the singular value decomposition, it is easy to prove that for any $r \in \RR,\lambda>0$, it holds that
\begin{equation}\label{matrix_inverse_lemma}
    (A^*A + \lambda)^{r}A^* = A^*(AA^*+\lambda)^{r}.
\end{equation}
% , where $\{f_i\}_{i=1}^\infty$ is an orthonormal basis of $\cH$. It can be shown that the trace does not depend on the choice of orthonormal basis and hence if $B$ is compact, $\tr(B) = \sum_{j=1}^\infty \lambda_j$. 
For more information on bounded linear operators in Hilbert spaces, we refer to \cite{kato2013perturbation}. 

\subsection{Preliminary Results on Singular Values of Compact Operators}
\begin{proposition}[{\cite[Proposition 2.1]{volkov2024eigenvalue}}]
\label{prop: minimax principle}
Let $\cH_1, \cH_2$ be two Hilbert spaces. Let $T : \cH_1 \to \cH_2$ be a compact linear operator.
Fix an integer $n \geq 1$ and let $\mathcal{S}_n$ be the set of all subspaces of $\cH_1$ with co-dimension less than or equal to $n - 1$. Then the singular value
\begin{equation}
s_n(T) = \min_{X \in \mathcal{S}_n} \max_{x \in X, \|x\|=1} \|Tx\|.
\end{equation}
\end{proposition}

\begin{proposition}[{\cite[Proposition 2.5]{volkov2024eigenvalue}}]
\label{prop:singular-value-inequality}
Let $\cH_1, \cH_2$ be two Hilbert spaces. Let $S, T : \cH_1 \to \cH_2$ be compact linear operators, and let $p, q \ge 1$ be two integers. Then
\begin{equation}
s_{p+q-1}(S+T) \le s_p(S) + s_q(T).
\end{equation}
\end{proposition}
\begin{proposition}
\label{prop: multi-singular-value-inequality}
Let $\cH_1$ and $\cH_2$ be Hilbert spaces, and let
$T_i:\cH_1\to\cH_2$, $i=1,\ldots,M$, be compact linear
operators. Then, for every $n\geq1$,
\[
s_n\left(\sum_{i=1}^M T_i\right)
\leq
\sum_{i=1}^M s_{\lceil n/M\rceil}(T_i).
\]
If, in addition, there exist
$T_{\max}\in\{T_1,\ldots,T_M\}$ and a constant $C>0$
such that
\[
s_k(T_i)\leq C s_k(T_{\max}),
\qquad i=1,\ldots,M,\quad k\geq1,
\]
then
\[
s_n\left(\sum_{i=1}^M T_i\right)
\leq CM s_{\lceil n/M\rceil}(T_{\max}),
\qquad n\geq1.
\]

% Define $T_{\max}\in\{T_i\}_{i=1}^{M}$ be the operator such that there exist a constant $C$ such that for any $i=1,2,\cdots M$ and $n\in \mathbb{N}^+$,
% \[
% s_n(T_i)\le C s_n(T_{\max}).
% \]
% Then for any $n\in \mathbb{N}^+$, 
% \[
% s_{n}\left(\sum_{i=1}^M T_i\right)  \le \sum_{i=1}^M s_{\lfloor \frac{n}{M} \rfloor}(T_i)\le C M s_{\lfloor \frac{n}{M} \rfloor}(T_{\max}).
% \]
\begin{proof}
Repeated application of
Proposition~\ref{prop:singular-value-inequality} gives
\[
s_{M(k-1)+1}\left(\sum_{i=1}^M T_i\right)
\leq
\sum_{i=1}^M s_k(T_i),
\qquad k\geq1.
\]
For any $n\geq1$, set $k=\lceil n/M\rceil$.
Since $M(k-1)+1\leq n$, monotonicity of the singular
values yields
\[
s_n\left(\sum_{i=1}^M T_i\right)
\leq
s_{M(k-1)+1}\left(\sum_{i=1}^M T_i\right)
\leq
\sum_{i=1}^M s_k(T_i).
\]
The final assertion follows from
$s_k(T_i)\leq C s_k(T_{\max})$.
\end{proof}
\end{proposition}

\begin{lemma}
\label{lemma: multiply by a bdd function}
Let $\Omega$ be a bounded open subset of $\mathbb R^d$
or a compact $d$-dimensional $C^1$ manifold embedded in
$\mathbb R^{\tilde d}$, and let $\Omega'$ be an open subset
of $\mathbb R^{d'}$.
Integration over $\Omega$ is with respect to Lebesgue
measure in the Euclidean case and the induced Riemannian
volume measure in the manifold case; both are denoted
by $\dd\bm y$.
Let $k\in L^2(\Omega'\times\Omega)$ with respect to the
corresponding product measure, and define
\[
K:L^2(\Omega)\to L^2(\Omega'),
\qquad
(K\varphi)(\bm x)
=
\int_\Omega
k(\bm x,\bm y)\varphi(\bm y)\,\dd\bm y.
\]
Let $a\in L^\infty(\Omega)$ satisfy
$\|a\|_{L^\infty(\Omega)}\leq M$ for some $M\geq0$,
and define
\[
\tilde K:L^2(\Omega)\to L^2(\Omega'),
\qquad
(\tilde K\varphi)(\bm x)
=
\int_\Omega
a(\bm y)k(\bm x,\bm y)\varphi(\bm y)\,\dd\bm y.
\]
Then
\[
s_n(\tilde K)\leq Ms_n(K),
\qquad n\geq1.
\]
\end{lemma}

\begin{proof}
Define $M_a:L^2(\Omega)\to L^2(\Omega)$ by
$M_a\varphi=a\varphi$.
Then $\|M_a\|\leq M$ and $\tilde K=KM_a$.
Since $K$ is Hilbert--Schmidt, both $K$ and $\tilde K$
are compact. By the finite-rank approximation characterization of
singular values
\cite[Chapter~VI, Theorem~1.5]{gohberg1990classes},
for every bounded linear operator
$Q:L^2(\Omega)\to L^2(\Omega')$ with
$\operatorname{rank}(Q)<n$, we have
\[
s_n(\tilde K)
\leq \|KM_a-QM_a\|
\leq M\|K-Q\|,
\]
where
$\operatorname{rank}(QM_a)\leq\operatorname{rank}(Q)<n$.
Taking the infimum over such $Q$ yields
$s_n(\tilde K)\leq Ms_n(K)$.
\end{proof}
\begin{theorem}[{\cite[Theorems~1.1 and~1.3]{volkov2024eigenvalue}}]
\label{thm: regularity-eigendecay}
 Let $\Omega$ be a bounded open set of $\mathbb{R}^d$ such that 
$\partial \Omega$ is Lipschitz regular. Let $\Omega'$ be an open set of 
$\mathbb{R}^{d'}$ and $k : \Omega' \times \Omega \to \mathbb{R}$ be in 
$L^2(\Omega') \times H^p(\Omega)$. Define the linear operator
\[
K : L^2(\Omega) \to L^2(\Omega'),
\qquad
Kf(\bm{x}) = \int_{\Omega} k(\bm{x},\bm{y})f(\bm{y})\dd \bm{y}.
\]
Then the singular values satisfy $s_n(K)=o\left(n^{-\frac{p}{d}-\frac{1}{2}}\right)$. In addition, if 
$k$ is analytic in a neighborhood of $\bar{\Omega}$ as a function valued in $L^2\left(\Omega^{\prime}\right)$, there exist $\tau \in (0,1)$ and $C>0$ such that 
\begin{equation*}
\label{ana_act}
 s_n(K)\leq C \tau^{n^{1/d}}.    
\end{equation*}
\end{theorem}

We remark that if $k(\bm{x},\bm{y})$ is real analytic for $(\bm{x},\bm{y})$ in a neighborhood of $\overline{\Omega' \times \Omega}$, it satisfies the analyticity conditions; we refer to \cite[Section 4.4]{volkov2024eigenvalue}
for more details.

\subsection{Useful Operator Inequalities}
\label{sec: useful-inequalties}
\begin{lemma}[{\cite[Lemma A.10]{long2025optimalratessaturationnoiseless}}]\label{lemma: AB-bound}
Let $A$ be a positive self-adjoint compact operator defined over a separable Hilbert space and $\alpha,\beta,\gamma\in (0,\infty)$. If $\beta\leq \alpha$, we have
$\|(A+\alpha )^\gamma (A+\beta)^{-\gamma}\|\leq (\alpha/\beta)^\gamma$.
\end{lemma}

%\begin{proposition}[Proposition 8 in \cite{rudi2017generalization}]\label{pro: operator-app}
%Let $\cH$ be a separable Hilbert space, let $A, B$ two bounded self-adjoint positive linear operators on $\cH$ and $\lambda>0$. Then
%$$
%\left\|(A+\lambda )^{-1 / 2} B^{1 / 2}\right\| \leq\left\|(A+\lambda )^{-1 / 2}(B+\lambda )^{1 / 2}\right\| \leq(1-\alpha)^{-1 / 2}
%$$
%with
%$$
%\alpha=\lambda_{\max }\left[(B+\lambda )^{-1 / 2}(B-A)(B+\lambda )^{-1 / 2}\right] .
%$$
%Note that $\alpha \leq \frac{\lambda_{\max }(B)}{\lambda_{\max }(B)+\lambda}<1$ by definition.
%\end{proposition}
%\wl{it seems that we did not use this proposition.}

The following concentration inequality originates from
\cite[Theorem~3.1]{minsker2017some}.
We use the formulation in
\cite[Theorem~A.2]{long2025optimalratessaturationnoiseless}.
\begin{theorem}[{\cite[Theorem~3.1]{minsker2017some};
\cite[Theorem~A.2]{long2025optimalratessaturationnoiseless}}]
\label{thm: concentration}
Let $X_1, \ldots, X_n $ be a sequence of independent self-adjoint random operators on a Hilbert space $\cH$ such that $\mathbb{E} X_i=0$ for $i=1, \ldots, n$ and $\left\|\sum_{i=1}^n \mathbb{E} X_i^2\right\|\leq \sigma^2$. Assume that $\left\|X_i\right\| \leq U$ almost surely for all $1 \leq i \leq n$ and some positive $U \in \mathbb{R}$. Then, for any $t \ge \frac{1}{6}\left(U+\sqrt{U^2+36 \sigma^2}\right)$,
\begin{equation}\label{eqn: 111}
\mathbb{P}\left(\left\|\sum_{i=1}^n X_i\right\|>t\right) \leq 14 \frac{\tr(\sum_{i=1}^n \mathbb{E} X_i^2)}{\sigma^2} \exp \left[-\frac{t^2 / 2}{\sigma^2+t U / 3}\right].
\end{equation}
\end{theorem}
% We remark that the above theorem is slightly tighter than the version stated in \cite{minsker2017some}*{Theorem 3.1}, where $\frac{\tr(\sum_{i=1}^n \mathbb{E} X_i^2)}{\sigma^2}$ in \eqref{eqn: 111} is replaced by $\frac{\tr(\sum_{i=1}^n \mathbb{E} X_i^2)}{\left\|\sum_{i=1}^n \mathbb{E} X_i^2\right\|}$. However,  \cite{minsker2017some} in fact has proved the tighter presented here; we refer to the paragraphs below Eq.~(3.9) in \cite{minsker2017some} for more details.
% It is worth noting that a larger $\sigma$ imposes a stricter condition on $t$, which limits the size of deviation for which we can provide concentration guarantees. 

% \lipsum[71]

% \begin{lemma}
% Test Lemma.
% \end{lemma}

% \section*{Acknowledgments}
% This work is partially supported by the NSFC Major Research Plan - Interpretable and General Purpose Next-generation Artificial Intelligence (Nos. 92570001 and 92370205), NSFC grant 12425113, the Basic Research Program of Jiangsu under Grant BK20253039 and Grant BK20250480, and Key Laboratory of the Ministry of Education for Mathematical Foundations and Applications of Digital Technology, University of Science and Technology of China.

\subsection{DoF and eigenvalue decay}
The behavior of $N(\gamma)$ can be explicitly characterized using the eigenvalue decay.
\begin{lemma}
\label{lemma: N(gamma)}
Let $\{\mu_j\}_{j=1}^{\infty}$ be a nonincreasing sequence
of positive numbers, and define
\[
N(\gamma)
:=
\sum_{j=1}^{\infty}\frac{\mu_j}{\mu_j+\gamma},
\qquad \gamma>0.
\]
Suppose that $\{\mu_j\}_{j=1}^{\infty}$ satisfies Assumption~\ref{assump:eigenvalue-decay}.
Then, as $\gamma\to 0^+$,
\[
N(\gamma)\asymp
\begin{cases}
\gamma^{-1/t},
& \text{if }\mu_j\asymp j^{-t},\\
[\log(1/\gamma)]^d,
& \text{if }\mu_j\asymp\tau^{-j^{1/d}}.
\end{cases}
\]
Consequently, $N(\mu_k)\lesssim k$ for all $k\geq 1$.
\end{lemma}

\begin{proof}
The asymptotic estimates for polynomial and exponential decay are
established in \cite[Lemma~3.3]{long2025optimalratessaturationnoiseless}.
Here we provide a direct proof for the stretched exponential case
$\mu_j\asymp\tau^{-j^{1/d}}$, noting that the argument below also recovers the exponential decay
regime when $d=1$.

Suppose that $\mu_j\asymp e^{-aj^{1/d}}$ with $a:=\log\tau>0$ and
$d\geq 1$. Then there exist constants $0<c_-\leq c_+$ such that
\[
c_-e^{-aj^{1/d}}
\leq\mu_j
\leq c_+e^{-aj^{1/d}},
\qquad j\geq1.
\]
For sufficiently small $\gamma\in(0, \min\{1/2, c_-e^{-a}\})$, define
\[
J_-:=
\left\lfloor
\left(\frac{\log(c_-/\gamma)}{a}\right)^d
\right\rfloor\geq 1,
\qquad
J_+:=
\left\lceil
\left(\frac{\log(c_+/\gamma)}{a}\right)^d
\right\rceil\geq 1.
\]
For $1\leq j\leq J_-$, monotonicity yields $\mu_j\geq c_-e^{-aJ_-^{1/d}}\geq\gamma$,
so that $\mu_j/(\mu_j+\gamma)\geq 1/2$. Retaining only the first $J_-$
terms in the series gives the lower bound
\[
N(\gamma)
\geq\sum_{j=1}^{J_-}\frac{\mu_j}{\mu_j+\gamma}
\geq\frac{J_-}{2}
\gtrsim(\log(1/\gamma))^d.
\]

For the upper bound, we split the series at $J_+$. Bounding the first
$J_+$ terms by $1$ and applying integral comparison to the tail, we obtain
\begin{align*}
N(\gamma)
&\leq J_++\frac{c_+}{\gamma}
\sum_{j>J_+}e^{-aj^{1/d}}
\leq J_++\frac{c_+}{\gamma}
\int_{J_+}^{\infty}e^{-ax^{1/d}}\,\dd x\\
&=J_++\frac{c_+d}{\gamma}
\int_{J_+^{1/d}}^{\infty}y^{d-1}e^{-ay}\,\dd y
\leq J_++C\frac{c_+}{\gamma}
J_+^{(d-1)/d}e^{-aJ_+^{1/d}}\\
&\leq J_++C J_+^{(d-1)/d}
\lesssim(\log(1/\gamma))^d.
\end{align*}
Here, setting $Y:=J_+^{1/d}\geq 1$, the second-to-last inequality follows from
the exact integration-by-parts identity
\[
\int_Y^{\infty}y^{d-1}e^{-ay}\,\dd y
=
e^{-aY}\sum_{m=0}^{d-1}\frac{(d-1)!}{m!\,a^{d-m}}Y^m
\leq
\left(\sum_{m=0}^{d-1}\frac{(d-1)!}{m!\,a^{d-m}}\right)Y^{d-1}e^{-aY},
\]
with $C:=d\sum_{m=0}^{d-1}\frac{(d-1)!}{m!\,a^{d-m}}>0$ depending only on $a$ and $d$,
while the last inequality uses $c_+e^{-aJ_+^{1/d}}\leq\gamma$. Finally, since $\mu_k\to 0^+$ as $k\to\infty$, substituting $\gamma=\mu_k$
into the above asymptotic bounds on $N(\gamma)$ immediately yields
$N(\mu_k)\lesssim k$ for all $k\geq 1$.
\end{proof}

\subsection{Useful Lemmas in the Proof of Theorem \ref{thm: Estimation_error upper}}

\begin{lemma}
\label{lemma:dof-rf}
For every $\lambda>0$,
\[
F_{\cL}(\lambda)
=
\operatorname*{ess\,sup}_{\bm x\in\Omega}
\frac{
\|(\widetilde{\cT}_{\cL}+\lambda I)^{-1/2}
\cL\phi(\bm x,\cdot)\|_{L^2(\pi)}^2
}{q_I(\bm x)},
\]
and
\[
F_{\cB}(\lambda)
=
\operatorname*{ess\,sup}_{\bm x\in\partial\Omega}
\frac{
\|(\widetilde{\cT}_{\cB}+\lambda I)^{-1/2}
\cB\phi(\bm x,\cdot)\|_{L^2(\pi)}^2
}{q_B(\bm x)}.
\]
\end{lemma}

\begin{proof}
We prove the first identity; the second follows by the same argument.
Write the singular value decompositions as
\begin{equation}
\label{eqn: svd}
\Phi_{\cL}
=
\sum_{j\geq1}
(\mu_j^{\cL})^{1/2}
e_j^{\cL}\otimes\tilde e_j^{\cL},
\qquad
\Phi_{\cL}^*
=
\sum_{j\geq1}
(\mu_j^{\cL})^{1/2}
\tilde e_j^{\cL}\otimes e_j^{\cL},
\end{equation}
where $\{e_j^{\cL}\}_{j\geq1}\subset L^2(\rho_I)$ and
$\{\tilde e_j^{\cL}\}_{j\geq1}\subset L^2(\pi)$ are
orthonormal eigenfunctions of
$\Phi_{\cL}\Phi_{\cL}^*$ and
$\widetilde{\cT}_{\cL}=\Phi_{\cL}^*\Phi_{\cL}$,
respectively.
As in \cite[Lemma~D.2]{long2025optimalratessaturationnoiseless},
for $\rho_I$-almost every $\bm x$,
\[
\cL\phi(\bm x,\cdot)
=
\sum_{j\geq1}
(\mu_j^{\cL})^{1/2}
e_j^{\cL}(\bm x)\tilde e_j^{\cL}(\cdot)
\quad\text{in }L^2(\pi).
\]
Consequently,
\[
\begin{aligned}
\|(\widetilde{\cT}_{\cL}+\lambda I)^{-1/2}
\cL\phi(\bm x,\cdot)\|_{L^2(\pi)}^2&=
\left\|
\sum_{j\geq1}
\left(\frac{\mu_j^{\cL}}{\mu_j^{\cL}+\lambda}\right)^{1/2}
e_j^{\cL}(\bm x)\tilde e_j^{\cL}(\cdot)
\right\|_{L^2(\pi)}^2\\
&=
\sum_{j\geq1}
\frac{\mu_j^{\cL}}{\mu_j^{\cL}+\lambda}
|e_j^{\cL}(\bm x)|^2.
\end{aligned}
\]
Dividing by $q_I(\bm x)$ and taking the essential supremum
proves the claim.
\end{proof}

\begin{lemma}[Concentration of the empirical operator]
\label{lemma:concentration}
Let $\delta\in(0,1)$ and $\lambda>0$. If
\begin{equation}
\label{lambda union condition}
\lambda\geq
\max\{\Lambda_I(n_I,\delta/2),\Lambda_B(n_B,\delta/2)\}, 
\end{equation}
then, with probability at least $1-\delta$,
\begin{equation}
\label{equ: concentration_operator}
\left\|
(\widetilde{\cT}_{\cR_\eta}+\lambda I)^{-1/2}
(A_n-\widetilde{\cT}_{\cR_\eta})
(\widetilde{\cT}_{\cR_\eta}+\lambda I)^{-1/2}
\right\|
\leq\frac34.
\end{equation}
\end{lemma}

\begin{proof}
Let $\Delta_n = (\widetilde{\cT}_{\cR_\eta}+\lambda I)^{-1/2}
(A_n-\widetilde{\cT}_{\cR_\eta})
(\widetilde{\cT}_{\cR_\eta}+\lambda I)^{-1/2}$. By the triangle
inequality,
\[
\begin{aligned}
\|\Delta_n\|
&\leq
\left\|
(\widetilde{\cT}_{\cR_\eta}+\lambda I)^{-1/2}
(\widetilde{\cT}_{\cL,n_I}-\widetilde{\cT}_{\cL})
(\widetilde{\cT}_{\cR_\eta}+\lambda I)^{-1/2}
\right\|\\
&\quad+
\left\|
(\widetilde{\cT}_{\cR_\eta}+\lambda I)^{-1/2}
\eta(\widetilde{\cT}_{\cB,n_B}-\widetilde{\cT}_{\cB})
(\widetilde{\cT}_{\cR_\eta}+\lambda I)^{-1/2}
\right\|\\
&=:M_1+M_2.
\end{aligned}
\]

\emph{Bounding $M_1$.}
Define
\[
h_i^I
=
\frac{
(\widetilde{\cT}_{\cR_\eta}+\lambda I)^{-1/2}
\cL\phi(\bm x_i,\cdot)
}{\sqrt{q_I(\bm x_i)}},
\]
and
\[
Z_i^I
:=
\frac1{n_I}\Big[
h_i^I\otimes h_i^I
-
(\widetilde{\cT}_{\cR_\eta}+\lambda I)^{-1/2}
\widetilde{\cT}_{\cL}
(\widetilde{\cT}_{\cR_\eta}+\lambda I)^{-1/2}
\Big].
\]
Note that
\[
\mathbb E[h_i^I\otimes h_i^I]
=
(\widetilde{\cT}_{\cR_\eta}+\lambda I)^{-1/2}
\widetilde{\cT}_{\cL}
(\widetilde{\cT}_{\cR_\eta}+\lambda I)^{-1/2},
\]
the $Z_i^I$ are independent, centered, self-adjoint
operators, and
\(
M_1=\left\|\sum_{i=1}^{n_I}Z_i^I\right\|.
\)

Since
$\widetilde{\cT}_{\cL}\preccurlyeq\widetilde{\cT}_{\cR_\eta}$, $
(\widetilde{\cT}_{\cR_\eta}+\lambda I)^{-1}
\preccurlyeq
(\widetilde{\cT}_{\cL}+\lambda I)^{-1}$. 
Combining this inequality with Lemma~\ref{lemma:dof-rf} gives
\[
\begin{aligned}
\|h_i^I\|_{L^2(\pi)}^2
&=
\frac{
\langle
\cL\phi(\bm x_i,\cdot),
(\widetilde{\cT}_{\cR_\eta}+\lambda I)^{-1}
\cL\phi(\bm x_i,\cdot)
\rangle_{L^2(\pi)}
}{q_I(\bm x_i)}\\
&\leq
\frac{
\|(\widetilde{\cT}_{\cL}+\lambda I)^{-1/2}
\cL\phi(\bm x_i,\cdot)\|_{L^2(\pi)}^2
}{q_I(\bm x_i)}\\
&\leq F_{\cL}(\lambda)
\leq F_I(\lambda),
\end{aligned}
\]
almost surely.
Hence $\|h_i^I\otimes h_i^I\|\leq F_I(\lambda)$ and $\|\mathbb E[h_i^I\otimes h_i^I]\|\leq F_I(\lambda)$. Thus
\[
\|Z_i^I\|\leq\frac{F_I(\lambda)}{n_I}.
\]
Moreover,
\[
\begin{aligned}
\mathbb E[(Z_i^I)^2]
&\preccurlyeq
\frac1{n_I^2}
\mathbb E\left[
\|h_i^I\|_{L^2(\pi)}^2
h_i^I\otimes h_i^I
\right]\\
&\preccurlyeq
\frac{F_I(\lambda)}{n_I^2}
(\widetilde{\cT}_{\cR_\eta}+\lambda I)^{-1/2}
\widetilde{\cT}_{\cL}
(\widetilde{\cT}_{\cR_\eta}+\lambda I)^{-1/2}\\
&\preccurlyeq
\frac{F_I(\lambda)}{n_I^2}
\widetilde{\cT}_{\cR_\eta}
(\widetilde{\cT}_{\cR_\eta}+\lambda I)^{-1}.
\end{aligned}
\]
Since $\|\widetilde{\cT}_{\cR_\eta}
(\widetilde{\cT}_{\cR_\eta}+\lambda I)^{-1}\|\leq 1$, 
\[
\left\|
\sum_{i=1}^{n_I}\mathbb E[(Z_i^I)^2]
\right\|
\leq\frac{F_I(\lambda)}{n_I}.
\]
Furthermore,
$\widetilde{\cT}_{\cR_\eta}
\preccurlyeq\frac{1}{\nu_-}\widetilde{\cT}_{\mathrm{PDE}}$
implies
\[
\widetilde{\cT}_{\cR_\eta}
(\widetilde{\cT}_{\cR_\eta}+\lambda I)^{-1}
\preccurlyeq
\widetilde{\cT}_{\mathrm{PDE}}
(\widetilde{\cT}_{\mathrm{PDE}}+\nu_-\lambda I)^{-1}.
\]
Taking traces yields
\[
\operatorname{tr}\left[
\widetilde{\cT}_{\cR_\eta}
(\widetilde{\cT}_{\cR_\eta}+\lambda I)^{-1}
\right]
\leq
N_{\mathrm{PDE}}(\nu_-\lambda)
\leq
N_{\mathrm{PDE}}(\lambda),
\]
where the last inequality follows from $\nu_-\geq1$
and the monotonicity of $N_{\mathrm{PDE}}$.
Thus
\[
\operatorname{tr}\left(
\sum_{i=1}^{n_I}\mathbb E[(Z_i^I)^2]
\right)
\leq
\frac{F_I(\lambda)N_{\mathrm{PDE}}(\lambda)}{n_I}
\leq\frac{F_I(\lambda)^2}{n_I}.
\]
Applying Theorem~\ref{thm: concentration} with
$U=\sigma^2=F_I(\lambda)/n_I$ yields
\[
\mathbb P(M_1>t)
\leq
14F_I(\lambda)
\exp\left(
-\frac{3n_It^2}{F_I(\lambda)(6+2t)}
\right),
\]
provided
\[
3n_It^2\geq3F_I(\lambda)+F_I(\lambda)t.
\]
For $t=3/8$, this condition is $n_I\geq8F_I(\lambda)$,
and $
\exp\left[
-3n_It^2/(F_I(\lambda)(6+2t))
\right]=\exp[-n_I/(16F_I(\lambda))]$.
Thus,
\begin{equation}
\label{N_I condition}
n_I\geq16F_I(\lambda)
\max\left\{1,\log\frac{28F_I(\lambda)}{\delta}\right\}
\end{equation}
implies
$\mathbb P(M_1\leq3/8)\geq1-\delta/2$.

\emph{Bounding $M_2$.}
Define
\[
h_i^B
=
\sqrt{\frac{\eta}{q_B(\bm x_{n_I+i})}}
(\widetilde{\cT}_{\cR_\eta}+\lambda I)^{-1/2}
\cB\phi(\bm x_{n_I+i},\cdot),
\]
and
\[
Z_i^B
:=
\frac1{n_B}\Big[
h_i^B\otimes h_i^B
-
\eta
(\widetilde{\cT}_{\cR_\eta}+\lambda I)^{-1/2}
\widetilde{\cT}_{\cB}
(\widetilde{\cT}_{\cR_\eta}+\lambda I)^{-1/2}
\Big].
\]
These operators are independent and centered, and
\(
M_2=\left\|\sum_{i=1}^{n_B}Z_i^B\right\|.
\)

Since
$\eta\widetilde{\cT}_{\cB}
\preccurlyeq\widetilde{\cT}_{\cR_\eta}$,
\[
\begin{aligned}
\|h_i^B\|_{L^2(\pi)}^2
&\leq
\frac{\eta}{q_B(\bm x_{n_I+i})}
\left\langle
\cB\phi(\bm x_{n_I+i},\cdot),
(\eta\widetilde{\cT}_{\cB}+\lambda I)^{-1}
\cB\phi(\bm x_{n_I+i},\cdot)
\right\rangle_{L^2(\pi)}\\
&=
\frac{
\|(\widetilde{\cT}_{\cB}+(\lambda/\eta)I)^{-1/2}
\cB\phi(\bm x_{n_I+i},\cdot)\|_{L^2(\pi)}^2
}{q_B(\bm x_{n_I+i})}\\
&\leq F_{\cB}(\lambda/\eta)
\leq F_B(\lambda),
\end{aligned}
\]
almost surely.
The same moment calculation gives
\[
\|Z_i^B\|\leq\frac{F_B(\lambda)}{n_B},
\]
and
\[
\sum_{i=1}^{n_B}\mathbb E[(Z_i^B)^2]
\preccurlyeq
\frac{F_B(\lambda)}{n_B}
\widetilde{\cT}_{\cR_\eta}
(\widetilde{\cT}_{\cR_\eta}+\lambda I)^{-1}.
\]
Therefore,
\[
\left\|
\sum_{i=1}^{n_B}\mathbb E[(Z_i^B)^2]
\right\|
\leq\frac{F_B(\lambda)}{n_B},
\qquad
\operatorname{tr}\left(
\sum_{i=1}^{n_B}\mathbb E[(Z_i^B)^2]
\right)
\leq\frac{F_B(\lambda)^2}{n_B}.
\]
Another application of Theorem~\ref{thm: concentration}
with $t=3/8$ shows that
\begin{equation}
\label{n_B condition}
n_B\geq16F_B(\lambda)
\max\left\{1,\log\frac{28F_B(\lambda)}{\delta}\right\}
\end{equation}
implies
$\mathbb P(M_2\leq3/8)\geq1-\delta/2$.

The condition \eqref{lambda union condition}
ensures that \eqref{N_I condition} and \eqref{n_B condition} hold. The union bound now gives
\[
\mathbb P\left(\|\Delta_n\|\leq\frac34\right)
\geq
\mathbb P\left(M_1\leq\frac38,\ M_2\leq\frac38\right)
\geq1-\delta.
\]
\end{proof}

\begin{lemma}
\label{lemma:constant bound}
On the event \eqref{equ: concentration_operator},
\[
\left\|
(\widetilde{\cT}_{\cR_\eta}+\lambda I)^{1/2}
(A_n+\lambda I)^{-1/2}
\right\|
\leq2.
\]
\end{lemma}

\begin{proof}
Since the operator in \eqref{equ: concentration_operator}
is self-adjoint, its norm bound implies
\[
(\widetilde{\cT}_{\cR_\eta}+\lambda I)^{-1/2}
(\widetilde{\cT}_{\cR_\eta}-A_n)
(\widetilde{\cT}_{\cR_\eta}+\lambda I)^{-1/2}
\preccurlyeq\frac34 I.
\]
Multiplying both sides of this inequality on the left
and right by
$(\widetilde{\cT}_{\cR_\eta}+\lambda I)^{1/2}$ gives
\[
\widetilde{\cT}_{\cR_\eta}-A_n
\preccurlyeq
\frac34(\widetilde{\cT}_{\cR_\eta}+\lambda I).
\]
Hence
\[
A_n+\lambda I
\succcurlyeq
\frac14(\widetilde{\cT}_{\cR_\eta}+\lambda I).
\]
Multiplying both sides of this inequality on the left
and right by
$(A_n+\lambda I)^{-1/2}$ gives
\[
\begin{aligned}
&(A_n+\lambda I)^{-1/2}
(\widetilde{\cT}_{\cR_\eta}+\lambda I)
(A_n+\lambda I)^{-1/2}\\
&\qquad\preccurlyeq4I.
\end{aligned}
\]
Taking operator norms proves the claim.
\end{proof}

\section{Proofs of Activation-Dependent Eigenvalue Bounds}
\label{proof_of_thm}

\subsection{Proof of Theorem \ref{activation_fun_rate}}
\label{proof: activation_fun_rate}
We first introduce a criterion that converts finite-rank
approximation estimates in the Hilbert--Schmidt norm
into singular-value decay bounds.
\begin{lemma}[Spectral decay from finite-rank approximation]
\label{lemma: finite-rank-spectral-decay}
Let $T:H_1\to H_2$ be Hilbert--Schmidt between separable
Hilbert spaces, and let $d\geq1$ and $s>0$.
Suppose that for all sufficiently small dyadic $h=2^{-\ell}$
there is a bounded linear operator $T_h:H_1\to H_2$ such that
\[
\operatorname{rank}T_h\leq C_0h^{-d},
\qquad
\|T-T_h\|_{\mathrm{HS}}\leq C_1h^s,
\]
with constants independent of $h$.
Then $s_j(T)=O(j^{-1/2-s/d})$.
If instead $\|T-T_h\|_{\mathrm{HS}}=o(h^s)$,
then $s_j(T)=o(j^{-1/2-s/d})$.
\end{lemma}
\begin{proof}
Set $N_h:=\lceil\max\{C_0,1\}h^{-d}\rceil$.
Since $s_{N_h+1}(T_h)=0$,
Proposition~\ref{prop:singular-value-inequality} gives
\[
s_{2N_h}(T)
\leq s_{N_h}(T-T_h)
\leq N_h^{-1/2}\|T-T_h\|_{\mathrm{HS}}
\leq C h^{s+d/2}.
\]
The second inequality follows from
$s_{N_h}(T-T_h)^2\leq1/N_h\sum_j s_j(T-T_h)^2$.
For every sufficiently large $j$, choose a dyadic $h$ with
$2N_h\leq j<2N_{h/2}$. Then $j\asymp h^{-d}$, and
$s_j(T)\leq s_{2N_h}(T)\leq C h^{s+d/2}=O(j^{-1/2-s/d}) $. The same argument
preserves little-$o$ estimates.
\end{proof}

To estimate the singular value decay of the boundary feature operator $\Phi_\cB$
without imposing additional smoothness on $\partial\Omega$, we introduce the following two lemmas.

\begin{lemma}[Ambient finite regularity]
\label{lemma: ambient-finite-decay}
Let $d\geq1$, let $q\geq1$ be an integer, and let
$X\subset\mathbb R^d$ be compact with a finite positive Borel
measure $\rho$. Let $E$ be a separable real Hilbert space.
Suppose $a:X\to E$ admits an $E$-valued $C^q$ extension to
an open neighborhood of $X$. Then
\[
K:L^2(X,\rho;\mathbb R)\to E,
\qquad Kf=\int_X a(\bm x)f(\bm x)\,\dd\rho(\bm x),
\]
is compact and $s_j(K)=o(j^{-1/2-q/d})$.
\end{lemma}
\begin{proof}
If $\rho(X)=0$, the assertion is immediate. Otherwise choose
a compact neighborhood of $X$ contained in the extension domain.
For sufficiently small dyadic $h$, a grid partitions $X$ into at
most $Ch^{-d}$ Borel sets of diameter at most $h$.
On each nonempty set use the degree-$q$ Taylor polynomial of
the extension about a point of that set. All line segments used
in Taylor's formula lie in the chosen compact neighborhood.
Uniform continuity of the $q$th derivatives gives the uniform
remainder bound
\[
\sup_{\bm x\in X}\|a(\bm x)-a_h(\bm x)\|_E=o(h^q).
\]
The integral operator $K_h$ with kernel $a_h$ consequently satisfies
\[
\operatorname{rank}K_h\leq C_qh^{-d},
\qquad
\|K-K_h\|_{\mathrm{HS}}^2
=\int_X\|a(\bm x)-a_h(\bm x)\|_E^2\,\dd\rho(\bm x)
=o(h^{2q}).
\]
The boundedness of $a$ and finiteness of $\rho$ make $K$
Hilbert--Schmidt. Applying 
Lemma~\ref{lemma: finite-rank-spectral-decay} with $s=q$ proves the lemma.
\end{proof}

\begin{lemma}[Ambient analyticity]
\label{lemma: ambient-analytic-decay}
Let $d\geq1$, let $X\subset\mathbb R^d$ be compact, and let
$\rho$ be a finite positive Borel measure on $X$.
Let $E$ be a separable real Hilbert space and $E_{\mathbb C}$
its complexification. Suppose $a:X\to E$ has a norm-holomorphic
extension $A:U\to E_{\mathbb C}$ on an open neighborhood
$U\subset\mathbb C^d$ of $X$.
Then the Bochner integral operator
\[
K:L^2(X,\rho;\mathbb R)\to E,
\qquad Kf=\int_X a(\bm x)f(\bm x)\,\dd\rho(\bm x),
\]
is compact, and there exist $C,c>0$, independent of $j$, such that
\begin{equation}
\label{eq: ambient-analytic-decay}
s_j(K)\leq C\exp(-cj^{1/d}),\qquad j\geq1.
\end{equation}
\end{lemma}
\begin{proof}
The case $\rho(X)=0$ is immediate. Write
$P_R(\bm x)=\{\bm z\in\mathbb C^d:
\max_i|z_i-x_i|<R\}$.
By compactness, choose $R>0$, centers $\bm x_1,\ldots,\bm x_N\in X$,
and a disjoint Borel partition $X=\bigcup_{\ell=1}^N F_\ell$
such that $F_\ell\subset P_{R/2}(\bm x_\ell)$ and
$\overline{P_R(\bm x_\ell)}\subset U$.
Set
\[
M:=\max_{1\leq\ell\leq N}
\sup_{\bm z\in\overline{P_R(\bm x_\ell)}}
\|A(\bm z)\|_{E_{\mathbb C}}<\infty.
\]
In each polydisc, the Taylor expansion and the vector-valued
Cauchy formula \cite[Theorem~5.1]{Kruse2020} give
\[
A(\bm z)=\sum_{\bm\alpha\in\mathbb N_0^d}
a_{\ell,\bm\alpha}(\bm z-\bm x_\ell)^{\bm\alpha},
\qquad
\|a_{\ell,\bm\alpha}\|_{E_{\mathbb C}}
\leq MR^{-|\bm\alpha|}.
\]
For an integer $m\geq1$ and $\bm x\in F_\ell$, define
\[
a_m(\bm x)
:=\sum_{0\leq\alpha_1,\ldots,\alpha_d<m}
\operatorname{Re}(a_{\ell,\bm\alpha})
(\bm x-\bm x_\ell)^{\bm\alpha}.
\]
Since $A(\bm x)=a(\bm x)\in E$ and real-part projection is
contractive, absolute convergence and the Cauchy estimate imply
\[
\begin{aligned}
\|a(\bm x)-a_m(\bm x)\|_E
&\leq\sum_{\max_i\alpha_i\geq m}
\|a_{\ell,\bm\alpha}\|_{E_{\mathbb C}}
\prod_i|x_i-x_{\ell,i}|^{\alpha_i}\\
&\leq M\sum_{\max_i\alpha_i\geq m}2^{-|\bm\alpha|}\\
&\leq Md\left(\sum_{q=m}^\infty2^{-q}\right)
\left(\sum_{q=0}^\infty2^{-q}\right)^{d-1}
=Md\,2^d2^{-m}.
\end{aligned}
\]
All multi-index sums in this display are over $\mathbb N_0^d$.
Let $K_m$ be the integral operator with kernel $a_m$.
Its range is contained in the span of the $Nm^d$ retained
coefficient vectors, and Cauchy--Schwarz gives
\[
\operatorname{rank}K_m\leq Nm^d,
\qquad
\|K-K_m\|\leq Md\,2^d\rho(X)^{1/2}2^{-m}.
\]
Thus $K$ is compact and
$s_{Nm^d+1}(K)\leq\|K-K_m\|$.
For $j>N$, take $m=\lfloor((j-1)/N)^{1/d}\rfloor$.
Increasing $C$ to cover the remaining indices proves
\eqref{eq: ambient-analytic-decay}.
\end{proof}
\begin{proof}[Proof of Theorem \ref{activation_fun_rate}]
Set $E:=L^2(\pi)$ and, for each derivative occurring below, write
\[
\phi_{\bm\beta}(\bm x)
:=\mathrm D_{\bm x}^{\bm\beta}\phi(\bm x,\cdot).
\]
Define
\[
K_{\bm\beta}f
:=\int_\Omega\phi_{\bm\beta}(\bm x)f(\bm x)\,\dd\bm x,
\qquad
K_{\bm\beta,\rho_I}f
:=\int_\Omega\phi_{\bm\beta}(\bm x)f(\bm x)\,\dd\rho_I(\bm x),
\]
with domains $L^2(\Omega,\dd\bm x)$ and $L^2(\rho_I)$,
respectively, and common codomain $E$. Write $M_g$ for multiplication by a real-valued function $g$ on the relevant $L^2$ space.
% the relevant $L^2$ space.
% The regularity estimates below ensure that these operators are
% Hilbert--Schmidt. Since $E$ is separable, it embeds isometrically
% into $L^2(0,1)$. Composition with this embedding preserves singular
% values and Bochner Sobolev regularity, so
% Theorem~\ref{thm: regularity-eigendecay} applies to these
% $E$-valued kernels.

% Write $M_g$ for multiplication by a real-valued function $g$ on
% the relevant $L^2$ space. For $g\in L^\infty(\Omega)$,
% \[
% (K_{\bm\beta,\rho_I}M_g)(K_{\bm\beta,\rho_I}M_g)^*
% =(K_{\bm\beta}M_{g\sqrt{w_I}})
%  (K_{\bm\beta}M_{g\sqrt{w_I}})^*.
% \]
% Indeed, both sides map $v\in E$ to
% \[
% \int_\Omega g(\bm x)^2w_I(\bm x)\phi_{\bm\beta}(\bm x)
% \langle\phi_{\bm\beta}(\bm x),v\rangle_E\,\dd\bm x.
% \]
% The two operators therefore have the same singular values.
Since $\dd \rho_I = w_I \dd \bm x$, Lemma~\ref{lemma: multiply by a bdd function} yields
\begin{equation}
\label{eq: activation-weighted-singular-values}
s_j(K_{\bm\beta,\rho_I}M_g)
\leq\|g\|_{L^\infty(\Omega)}
\|w_I\|_{L^\infty(\Omega)}^{1/2}s_j(K_{\bm\beta}).
\end{equation}
Let $\bm e_i$ denote the $i$th coordinate vector. Then
\[
\widetilde{\mathcal T}_r
=\sum_{|\bm\beta|\leq r}
K_{\bm\beta}K_{\bm\beta}^*,
\]
\[
\Phi_{\cL}^*
=-\sum_{i,\ell=1}^d
 K_{\bm e_i+\bm e_\ell,\rho_I}M_{a^{i\ell}}
 +\sum_{i=1}^dK_{\bm e_i,\rho_I}M_{b^i}
 +K_{\bm0,\rho_I}M_c,
\qquad
\Phi_{\cB}^*f
=\int_{\partial\Omega}\phi_{\bm0}(\bm x)f(\bm x)
\,\dd\rho_B(\bm x).
\]

\emph{Finite regularity.}
For $|\bm\alpha|=p$ and $\bm k\ne\bm0$,
the weak chain rule and integration along lines parallel to
$\bm k$ give
\begin{equation}
\label{eq: activation-slicing}
\int_\Omega
|\mathrm D_{\bm x}^{\bm\alpha}\phi(\bm x,\bm\theta)|^2
\,\dd\bm x
\leq C|\bm k|^{2p-1}
\|\sigma^{(p)}\|_{L^2(-M,M)}^2.
\end{equation}
% To see the factor $|\bm k|^{2p-1}$, write
% $\bm x=\bm y+t\bm k/|\bm k|$ with $\bm y\perp\bm k$.
% The chain rule contributes $|\bm k^{\bm\alpha}|^2
% \leq|\bm k|^{2p}$, while the substitution
% $s=|\bm k|t+b$ contributes $|\bm k|^{-1}$.
% The measure of the orthogonal projection of $\Omega$ is uniformly
% bounded because $\Omega$ is bounded.
For $\bm k=\bm0$, all positive-order spatial derivatives vanish.
Compactness of $\Theta$ and continuity of the lower derivatives
of $\sigma$ consequently imply
\[
\sup_{\bm\theta\in\Theta}
\|\phi(\cdot,\bm\theta)\|_{H^p(\Omega,\dd\bm x)}<\infty.
\]
Fubini's theorem, applied also to the weak derivatives, gives
\[
\phi_{\bm\beta}\in H^{p-|\bm\beta|}(\Omega;E),
\qquad |\bm\beta|\leq p,
\]
where the Sobolev space uses Lebesgue measure.
Theorem~\ref{thm: regularity-eigendecay} therefore gives
\[
s_j(K_{\bm\beta})
=o\!\left(j^{-\frac{p-|\bm\beta|}{d}-\frac12}\right).
\]
For $|\bm\beta|=p$, this follows instead from
$\sum_j s_j(K_{\bm\beta})^2<\infty$ and monotonicity.
Combining \eqref{eq: activation-weighted-singular-values}
with Proposition~\ref{prop: multi-singular-value-inequality} yields
\[
s_j(\Phi_{\cL}^*)
=o\!\left(j^{-\frac{p-2}{d}-\frac12}\right).
\]

The strict choice of $M$ and compactness of $\Theta$ imply that
$\bm x\mapsto\phi_{\bm0}(\bm x)$ has an $E$-valued
$C^{p-1}$ extension to a neighborhood of $\overline\Omega$.
Lemma~\ref{lemma: ambient-finite-decay}, applied with
$X=\partial\Omega$ and $\rho=\rho_B$, gives
\[
s_j(\Phi_{\cB}^*)
=o\!\left(j^{-\frac{p-1}{d}-\frac12}\right).
\]
Let $J:=\#\{\bm\beta\in\mathbb N_0^d:|\bm\beta|\leq r\}$.
Applying Proposition~\ref{prop: multi-singular-value-inequality},

% to the positive covariance operators, and using
% $s_j(T)=s_j(T^*)$, gives

\begin{equation}
\label{eq: activation-covariance-sums}
\begin{aligned}
\tilde\mu_j^{(r)}
&\leq\sum_{|\bm\beta|\leq r}
s_{\lceil j/J\rceil}(K_{\bm\beta})^2,\\
\mu_j^{\mathrm{PDE}}
&\leq\nu_I s_{\lceil j/2\rceil}(\Phi_{\cL}^*)^2
+\eta\nu_B s_{\lceil j/2\rceil}(\Phi_{\cB}^*)^2.
\end{aligned}
\end{equation}
Thus
\[
\tilde\mu_j^{(r)}
=o\!\left(j^{-\frac{2(p-r)}d-1}\right),
\qquad
\mu_j^{\mathrm{PDE}}
=o\!\left(j^{-\frac{2(p-2)}d-1}\right).
\]

\emph{Analytic regularity.}
Real analyticity on an open interval containing $[-M,M]$
provides a holomorphic extension of $\sigma$ to a complex
neighborhood of this interval segment.
By compactness of $\overline\Omega$ and $\Theta$, there is a
common complex neighborhood $U$ of $\overline\Omega$ such that
for every $|\bm\beta|\leq\max\{r,2\}$,
\[
\bm z\longmapsto
\left[\bm\theta\longmapsto
\bm k^{\bm\beta}\sigma^{(|\bm\beta|)}(\bm k^\top\bm z+b)
\right]
\]
is a norm-holomorphic $L^2(\pi;\mathbb C)$-valued extension of
$\phi_{\bm\beta}$. Indeed, after shrinking $U$, all values
$\bm k^\top\bm z+b$, with $\bm z\in U$ and
$\bm\theta=(\bm k,b)\in\Theta$, lie in a fixed compact
subset of the holomorphic extension domain of $\sigma$.
Cauchy estimates and compactness of $\Theta$ then
provide uniform bounds for the derivatives with
respect to $\bm z$, which justify differentiation
in $L^2(\pi;\mathbb C)$.

Apply Lemma~\ref{lemma: ambient-analytic-decay} to
$X=\overline\Omega$ with Lebesgue measure restricted to $\Omega$,
and to $X=\partial\Omega$ with measure $\rho_B$ and
$a=\phi_{\bm0}$. Taking common constants over the finitely
many derivatives gives
\[
s_j(K_{\bm\beta})\leq C_0e^{-c_0j^{1/d}},
\qquad
s_j(\Phi_{\cB}^*)\leq C_0e^{-c_0j^{1/d}}.
\]
Equation~\eqref{eq: activation-weighted-singular-values}
and Proposition~\ref{prop: multi-singular-value-inequality}
then give
$s_j(\Phi_{\cL}^*)\leq C_1e^{-c_1j^{1/d}}$.
Using \eqref{eq: activation-covariance-sums} and adjusting
the constants proves
\[
\tilde\mu_j^{(r)}\leq C\tau^{2j^{1/d}},
\qquad
\mu_j^{\mathrm{PDE}}\leq C\tau^{2j^{1/d}},
\qquad j\geq1,
\]
for some $C>0$ and $\tau\in(0,1)$ independent of $j$.
\end{proof}

\subsection{Proof of Theorem \ref{thm: reluk-spectral-decay}}
\label{proof: thm: reluk-spectral-decay}
\begin{proof}
Let $Q_0$ be a cube containing $\overline\Omega$.
For sufficiently small dyadic $h$, a uniform grid of side length
$h$ meets $Q_0$ in $O(h^{-d})$ cubes.
Let $V_h\subset L^2(\Omega)$ consist of restrictions to $\Omega$
of functions that are polynomials of total degree at most $k$
on each grid cube, and let $P_h$ be the orthogonal projection
onto $V_h$. Then $\dim V_h\leq Ch^{-d}$.

Fix $\bm\theta=(\bm k,b)$ and a multi-index $\bm\beta$ with
$q:=|\bm\beta|\leq k$.
If $\bm k\ne\bm0$, the weak derivative is, almost everywhere,
\[
\mathrm D_{\bm x}^{\bm\beta}\phi(\bm x,\bm\theta)
=\frac{k!}{(k-q)!}\bm k^{\bm\beta}
(\bm k^\top\bm x+b)_+^{k-q},
\]
where $t_+^0:=\mathbf1_{\{t>0\}}$.
Call a grid cube cut if its closure meets
$\{\bm x:\bm k^\top\bm x+b=0\}$, and let $S_{h,\bm\theta}$
be the union of the cut cubes meeting $Q_0$.
Within a fixed bounded enlargement of $Q_0$, this union lies
in a slab of thickness $2\sqrt d\,h$; hence
$|S_{h,\bm\theta}|\leq Ch$, uniformly in $\bm\theta$.
Outside this union the derivative is a polynomial of degree
at most $k-q$. On the cut cubes,
\[
|\mathrm D_{\bm x}^{\bm\beta}\phi(\bm x,\bm\theta)|
\leq C|\bm k|^k h^{k-q}\leq Ch^{k-q}.
\]
Approximate the derivative by itself on the uncut cubes and
by zero on the cut cubes. This defines a function in $V_h$.
The best-approximation
property of $P_h$ gives
\begin{equation}
\label{eq: reluk-projection-error}
\|(I-P_h)\mathrm D_{\bm x}^{\bm\beta}
\phi(\cdot,\bm\theta)\|_{L^2(\Omega)}^2
\leq Ch^{2(k-q)+1}.
\end{equation}
If $\bm k=\bm0$, the feature is constant and all its
positive-order derivatives vanish, so the same estimate holds.
All constants are uniform over the compact set $\Theta$.

Use the operators $K_{\bm\beta}$ defined in the proof
of Theorem~\ref{activation_fun_rate}, then
\[
\begin{aligned}
\operatorname{rank}(K_{\bm\beta}P_h)&\leq Ch^{-d},\\
\|K_{\bm\beta}-K_{\bm\beta}P_h\|_{\mathrm{HS}}^2
&=\int_\Theta
\|(I-P_h)\mathrm D_{\bm x}^{\bm\beta}
\phi(\cdot,\bm\theta)\|_{L^2(\Omega)}^2
\,\dd\pi(\bm\theta)\\
&\leq Ch^{2(k-q)+1}.
\end{aligned}
\]
Lemma~\ref{lemma: finite-rank-spectral-decay}, with
$s=k-q+\tfrac12$, yields
\[
s_j(K_{\bm\beta})
\leq Cj^{-1/2-(k-q+1/2)/d}.
\]
The first inequality in
\eqref{eq: activation-covariance-sums}
then proves \eqref{eq: reluk-eigenvalue-decay}.

Moreover, \eqref{eq: activation-weighted-singular-values}
with $g\equiv1$ gives
\[
s_j(K_{\bm\beta,\rho_I})
\leq \|w_I\|_{L^\infty(\Omega)}^{1/2}
s_j(K_{\bm\beta})
\leq Cj^{-1/2-(k-q+1/2)/d}.
\]
For the PDE statement, let $k\geq2$.
Lemma~\ref{lemma: multiply by a bdd function} and
Proposition~\ref{prop: multi-singular-value-inequality}, applied
to the decomposition of $\Phi_{\cL}^*$, give
\[
s_j(\Phi_{\cL}^*)
\leq Cj^{-1/2-(k-3/2)/d}.
\]
Since $t_+^k\in C^{k-1}(\mathbb R)$,
Lemma~\ref{lemma: ambient-finite-decay} gives
\[
s_j(\Phi_{\cB}^*)=o(j^{-1/2-(k-1)/d}).
\]
The second inequality in \eqref{eq: activation-covariance-sums}
then proves the claimed PDE eigenvalue bound.
\end{proof}

\end{appendices}


\bibliography{refabbr}% common bib file

\end{document}